\documentclass[reqno]{amsart}
\usepackage{amsmath,amssymb,stmaryrd}
\newtheorem{thm}{Theorem}[section]
\newtheorem{lem}{Lemma}[section]
\newtheorem{prop}{Proposition}[section]

\def \eps {\epsilon}
\def \rr {\mathbb{R}}
\def \nn {\mathbb{N}}
\def \rn {\mathbb{R}^n}
\def \rnm {\mathbb{R}^n_-}

\def \ne {\nu_\eps}
\def \ze {z_\eps}
\def \re {r_\eps}

\def \pe {p_\eps}
\def \tge {\tilde{g}_\eps}

\def \xe {x_\eps}
\def \ye {y_\eps}
\def \fe {f_\eps}
\def \ae {a_\eps}
\def \ue {u_\eps}

\def \ke {k_\eps}
\def \be {\beta_\eps}
\def \ve {v_\eps}
\def \elle {\ell_\eps}
\def \tv {\tilde{v}}
\def \tu {\tilde{u}}
\def \tve {\tilde{v}_\eps}
\def \tye {\tilde{y}_\eps}
\def \tue {\tilde{u}_\eps}
\def \tuei {\tilde{u}_{\eps,i}}
\def \tueun {\tilde{u}_{\eps,1}}
\def \tui {\tilde{u}_i}

\def \tze {\tilde{z}_\eps}

\def \bye {\bar{y}_\eps}
\def \mei {\mu_{\eps,i}}
\def \kei {k_{\eps,i}}
\def \meun {\mu_{\eps,1}}
\def \keun {k_{\eps,1}}
\def \xeun {x_{\eps,1}}
\def \meN {\mu_{\eps,N}}
\def \keN {k_{\eps,N}}
\def \crit {2^{\star}}

\def \huno {H_{1,0}^2(\Omega)}
\def \hunrn {H_{1,0}^2(\rn)}
\def \hunrnm {H_{1,0}^2(\rnm)}

\def \ds {\displaystyle}
\def \beq {\begin{eqnarray*}}
\def \eeq {\end{eqnarray*}}
\def \beqn {\begin{eqnarray}}
\def \eeqn {\end{eqnarray}}
\def \bequa {\begin{equation}}
\def \eequa {\end{equation}}

\begin{document}

\title[Concentration estimates
for Emden-Fowler equations]{Concentration estimates
for Emden-Fowler equations with boundary singularities and critical 
growth}
\author{N. Ghoussoub}
\address{Nassif Ghoussoub, Department of Mathematics, University of
British Columbia, Vancouver, Canada}
 \email{nassif@math.ubc.ca}
\author{F. Robert}
\address{Fr\'ed\'eric Robert, Laboratoire J.A.Dieudonn\'e, Universit\'e de
Nice Sophia-Antipolis,
Parc Valrose, 06108 Nice cedex 2, France}
\email{frobert@math.unice.fr}
\thanks{Both authors gratefully acknowledge the support of the Natural
Sciences and Engineering Research Council of Canada and the hospitality of
the University of British Columbia where this work was initiated.}

\date{March, 13th 2005}
\begin{abstract} We establish --among other things-- existence and
multiplicity of solutions for the Dirichlet problem
$\sum_i\partial_{ii}u+\frac{|u|^{\crit-2}u}{|x|^s}=0$ on smooth bounded
domains $\Omega$ of $ \rn$ ($n\geq 3$) involving the critical
Hardy-Sobolev exponent $\crit =\frac{2(n-s)}{n-2}$ where $0<s<2$, and in
the case where zero (the point of singularity) is on the boundary
$\partial \Omega$. Just as in the Yamabe-type non-singular framework
(i.e., when $s=0$), there is no nontrivial solution under global convexity
assumption (e.g., when $\Omega$ is star-shaped around $0$). However, in
contrast to the non-satisfactory situation of the non-singular case, we
show the existence of an infinite number of solutions under an assumption
of local strict concavity of $\partial \Omega$ at $0$ in at least one
direction.
More precisely, we need  the principal curvatures of $\partial
\Omega$ at $0$ to be  non-positive but not all vanishing. We also show 
that
the best constant in the Hardy-Sobolev inequality
   is attained as long as the mean curvature of $\partial \Omega$ at $0$ is 
negative, extending  the results of
\cite{gk} and completing our result of \cite{gr1} to include dimension 
$3$.  The key ingredients in our proof are refined  concentration
estimates which yield compactness for certain Palais-Smale sequences which
do not hold in the non-singular case.
\end{abstract}
\maketitle

\section{Introduction and statement of the results}
We address the problem of existence and multiplicity of possibly 
sign-changing solutions of the following Emden-Fowler boundary value 
problem on a smooth domain $\Omega$ of $\rn$, $n\geq 3$:

\bequa\label{def:equa}
\left\{\begin{array}{ll}
\Delta u=\frac{|u|^{\crit-2}u}{|x|^s}& \hbox{ in }{\mathcal D}'(\Omega)\\
 \quad u=0&\hbox{ on }\partial\Omega.
\end{array}\right.
\eequa
where here and throughout the paper, $\Delta=-\sum_i\partial_{ii}$ is the
Laplacian with minus sign convention, and 
$\crit:=2^*(s)=\frac{2(n-s)}{n-2}$ with
$s\in [0,2]$. The non-singular case, i.e., when $s=0$, is the Euclidean 
version of the celebrated
Yamabe problem considered first by Brezis and Nirenberg \cite{bn} followed 
by a large number of authors.
Here again the situation is interesting since we are dealing
with the corresponding critical exponent  in the Hardy-Sobolev embedding 
$\huno \to L^p(\Omega; |x|^{-s} dx)$ which is not compact when $p=2^*(s)$. 
We recall that $\huno$ is the completion of $C_c^\infty(\Omega)$, the set 
of smooth
functions compactly supported in $\Omega$, for the norm
$\Vert u\Vert_{\huno}=\sqrt{\int_{\Omega}|\nabla u|^2\, dx}$, and that the 
above embedding follows from the
Hardy-Sobolev inequality (\cite{ca}, \cite{cw}, \cite{gy}) which states 
that
  the constant
defined as
 \bequa\label{def:mus}
   \mu_s (\Omega) := \inf \left \{ \int_{\Omega}| \nabla u|^2
dx;\, u \in \huno \hbox{ and }  \int_{\Omega} \frac {|u|^{\crit}}{|x|^s}\,
dx =1\right\}
\eequa
satisfies $0<\mu_s(\Omega) <+\infty$. This in turn allows for a
variational approach for  the problem of finding solutions in $\huno\cap 
C^0(\overline{\Omega})$ for the Dirichlet problem (\ref{def:equa}).

Now the story of the state of the art in the non-singular case is quite 
extensive (see for instance Struwe \cite{st}), but for our purpose we 
single out the following highlights:

1) For any domain  $\Omega$, the best constant $\mu_0(\Omega) $ is the 
same as $\mu_0(\rn)$  and it  is never attained unless $\Omega$ is 
essentially $\rn$ (i.e., ${\rm cap}(\rn\setminus \Omega) =0$), in which 
case there is an infinite number
of sign-changing solutions for
\bequa\label{def:equa.orig}
\left\{\begin{array}{ll}
\Delta u=|u|^{\crit-2}u & \hbox{ in }{\mathcal D}'(\Omega)\\
 \quad u=0&\hbox{ on }\partial\Omega.
\end{array}\right.
\eequa
Moreover, there are no solution for (\ref{def:equa.orig})  whenever 
$\Omega$ is bounded convex or star-shaped. On the other hand, there
are solutions if $\Omega$ is not contractible (in dimension 3) and an 
infinite number of them \cite{bc}, if the domain $\Omega$ has non-trivial 
homology (i.e., $H_d(\Omega, {\bf Z}_2)\neq 0$ for some $d>0$). 
Unfortunately, these topological conditions are far from being optimal and 
no geometric condition that would guarantee the existence of one or more 
solutions, have so far been isolated.

2) On the other hand, the addition of a linear term to the equation, such 
as
\bequa\label{def:equa.linear}
\left\{\begin{array}{ll}
\Delta u=|u|^{\crit-2}u +\lambda u& \hbox{ in }{\mathcal D}'(\Omega)\\
\quad u=0&\hbox{ on }\partial\Omega.
\end{array}\right.
\eequa
improves the situation dramatically, especially when $0<\lambda 
<\lambda_1$, since there is then a positive solution  for any smooth 
bounded domain $\Omega$ in $\rn$ as long as $n\geq 4$ (See 
Brezis-Nirenberg \cite{bn}). The case $n=3$ is more delicate and was dealt 
with by Druet \cite{d2}. Most relevant to our work, are the recent 
results by Devillanova and Solimini who managed in a remarkable paper 
\cite{ds}, to establish  the existence of an infinite number of solutions 
for (\ref{def:equa.linear}) in dimension $n\geq 7$. \\

The situation for the Emden-Fowler equations (i.e., when $s>0$) turned out 
to be at least as interesting, and somewhat more satisfactory. Actually, 
the case when $0$ belongs to the interior of the domain $\Omega$ is almost 
identical to the non-singular case  \cite{gy} as one can prove essentially 
the same results with a suitable adaptation of  the same techniques.
However, the situation is much different when $0\in \partial \Omega$.

1) Indeed,  Egnell showed in  \cite{eg} that for open
cones of the form ${C}=\{x\in \rn; x=r\theta, \theta\in D \ {\rm and}\
r>0\}$ where the base $D$ is a connected domain of the unit sphere
$S^{n-1}$ of $\rn$, the best constant $\mu_{s}(C)$ is attained for
$0<s<2$ even when $\bar {C} \neq \rn$.
 The case where $\partial \Omega$ is smooth at $0$ was tackled in \cite{gk}
and it turned out to be also quite interesting since the curvature of the
boundary at $0$ gets to play an important role.
It was shown there
that in dimension $n\geq 4$, the negativity of all principal
curvatures\footnote{In our context, we specify the orientation of
$\partial\Omega$ in such a way that the normal vectors of $\partial
\Omega$ are pointing outward from the domain $\Omega$.}  at $0$ --which is
essentially a condition of {\it ``local strict concavity''} -- leads to
attainability of the best constant for problems with Dirichlet boundary
conditions, while the Neumann problems required the positivity of the mean
curvature at $0$.

More recently, we show in  \cite{gr1} that for dimension $n\geq 4$, the 
negativity of the mean
curvature of $\partial \Omega$ at $0$ is sufficient to ensure the 
attainability of $\mu_s (\Omega)$. This result is quite satisfactory, 
since standard Pohozaev type arguments show  non-attainability in the case 
where $\Omega$ is convex or star-shaped at  $0$. One of the results  of 
this paper is  the extension of this attainability result  to  cover all 
dimensions (greater than $3$)  including the more subtle context of 
dimension $3$. We shall establish the following

\begin{thm}\label{th:intro}
Let $\Omega$ be a smooth bounded oriented domain of $\rn$, $n\geq 3$, such
that $0\in\partial\Omega$ and assume $s\in (0,2)$.  If the mean curvature 
of $\partial\Omega$ at $0$ is negative, then $\mu_{s}(\Omega)$ is achieved 
by a
positive function which  is --a positive multiple of-- a solution for
\bequa\label{def:equapositive}
\left\{\begin{array}{ll}
\Delta u=\frac{|u|^{\crit-2}u}{|x|^s}& \hbox{ in }{\mathcal D}'(\Omega)\\
u>0&\hbox{ in }\Omega\\
u=0&\hbox{ on }\partial\Omega.
\end{array}\right.
\eequa
\end{thm}

  2)  As to the question of multiplicity of solutions for (\ref{def:equa}), 
we note that Ghoussoub-Kang had shown in \cite{gk} the existence of  two 
solutions under the assumption that all principal 
curvatures at $0$ are negative.  More precisely, assuming that the principal curvatures
$\alpha_1,...,\alpha_{n-1}$ of $\partial \Omega$ at $0$ are finite,  the
oriented boundary $\partial\Omega$ near the origin can then be represented
(up to rotating the coordinates if necessary) by $x_1 = \varphi_0 (x') =
-\frac{1}{2} \sum_{i=2}^{n}\alpha_{i-1} x^2_i +o(|x'|^2)$, where $x'=(x_2,
...,x_{n}) \in B_\delta(0) \cap \{ x_1 =0 \}$ for some $\delta >0$ where
$B_\delta(0)$ is the ball in $\rn$ centered at $0$ with radius $\delta$.
If the principal curvatures at $0$ are all negative, i.e.,  if
\begin{equation}
\max_{1 \leq i \leq n-1}\alpha_i <0,
\end{equation}
then the sectional curvature at $0$ is negative and therefore $\partial
\Omega$ --viewed as an $(n-1)$-Riemannian submanifold of $\rn$-- is
strictly convex at $0$ (see for instance \cite{ghl}). The latter property
means that there exists a neighborhood $U$ of $0$ in $\partial \Omega$,
such that the whole of $U$ lies on one side of a hyperplane $H$ that is
tangent to $\partial \Omega$ at $0$ and $U\cap H=\{0\}$, and so does the
complement $\rn\setminus\Omega$, at least locally. In other words, the 
above curvature
condition then amounts to a notion of strict local convexity of $\rn
\setminus \Omega$ at $0$.
In this paper, we complete and extend these results  in many  ways, since 
we establish the existence
of infinitely many solutions under the following much weaker assumption:
\begin{equation}
\max_{1 \leq i \leq n-1}\alpha_i \leq 0 \quad {\rm and} \quad  \min_{1
\leq i \leq n-1}\alpha_i <0.
\end{equation}
which  is a condition of {\it ``local concavity at $0$''} that is 
``strict"  in at least one direction.

\begin{thm}\label{th:multi}
Let $\Omega$ be a smooth bounded oriented domain of $\rn$, $n\geq 3$, such
that $0\in\partial\Omega$. Let $s\in (0,2)$ and  $a\in
C^1(\overline{\Omega})$ be such that the operator $\Delta+a$ is coercive 
in
$\Omega$. If the principal curvatures of $\partial\Omega$ at $0$ are
non-positive, but  not all vanishing, then there exists an infinite number
of solutions  $u\in \huno\cap C^1(\overline{\Omega})$ for
$$\left\{\begin{array}{ll}
\Delta u+au=\frac{|u|^{\crit-2}u}{|x|^s}& \hbox{ in }{\mathcal
D}'(\Omega)\\
u=0&\hbox{ on }\partial\Omega.
\end{array}\right.$$
   \end{thm}
\noindent We do not know if the negativity of the mean curvature at 
$0$ is sufficient for the above result, however
it is a remarkably satisfactory once compared to what is known in the 
nonsingular case and  since --as mentioned above-- we have no solution 
when $\Omega$ is convex or star-shaped at $0$.

3)  All these results rely on blow-up analysis techniques where the 
limiting
spaces (i.e., on which  the blown-up solutions of corresponding
Euler-Lagrange equations eventually live) play an important role. In the
non-singular case, the limiting space is $\rn$ while in our framework,
the limiting cases occur on half-spaces of the form  $\rnm=\{x\in\rnm/\, 
x_1<0\}$, where $x_1$ denotes
the first coordinate of a generic point $x\in\rn$ in the canonical basis
of $\rn$.
  The above theorem is a corollary of a more powerful result
established below about the asymptotic behaviour of a family of solutions
to elliptic pde's,  which are not necessarily minimizing sequences. We
actually study families of solutions to related  subcritical problems, and
we completely describe their asymptotic behaviour --potentially developing
a singularity at zero-- as we approach the critical exponent.

More
precisely, we say  that a function is in $C^1(\overline{\Omega})$ if it 
can
be extended to a $C^1-$function in a open neighborhood of
$\overline{\Omega}$, and consider a family $(\ae)_{\eps>0}\in 
C^1(\overline{\Omega})$ and a function $a\in C^1(\overline{\Omega})$ such 
that there exists an open subset ${\mathcal U}\subset\rn$ such that $\ae, 
a$ can be extended to ${\mathcal U}$ by $C^1-$functions that we still 
denote by $\ae,a$. We assume that they satisfy

\bequa\label{hyp:ae}
\overline{\Omega}\subset\subset{\mathcal U}\hbox{ and }\lim_{\eps\to
0}\ae=a\hbox{ in }C^1({\mathcal U}).
\eequa
Here is the main result of this paper.

\begin{thm}\label{th:cpct}
Let $\Omega$ be a smooth bounded oriented domain of $\rn$, $n\geq 3$,
such that $0\in\partial\Omega$. Assume $s\in (0,2)$ and consider
$(\ae)_{\eps>0}\in C^1(\overline{\Omega})$ such that (\ref{hyp:ae}) hold. 
We let $(\pe)_{\eps>0}$ such that $\pe\in
[0,\crit-2)$ for all $\eps>0$ and $\lim_{\eps\to 0}\pe=0$. We assume that
the principal curvatures of $\partial\Omega$ at $0$ are non-positive but
do not all vanish. We consider a family of functions
$(\ue)_{\eps>0}$ that  is uniformly bounded in $\huno$ and satisfying
$$\left\{\begin{array}{ll}
\Delta \ue+\ae\ue=\frac{|\ue|^{\crit-2-\pe}}{|x|^s}\ue& \hbox{ in
}{\mathcal D}'(\Omega)\\
\hskip 40pt \ue=0&\hbox{ on }\partial\Omega.
\end{array}\right.$$

1)   The family $(\ue)_{\eps>0}$ is then pre-compact in the $C^1-$topology. In
particular, there exists
$u_0\in \huno$ such that, up to a subsequence, we have that $\lim_{\eps\to
0}\ue=u_0$ in $C^1(\overline{\Omega})$.

2)  Moreover, if the $\ue$'s are nonnegative for all $\eps>0$, then the 
same conclusion holds under the sole hypothesis that the mean curvature of 
$\partial \Omega$ at $0$ is negative.
\end{thm}
The proof of this last theorem uses the machinery developed in 
Druet-Hebey-Robert \cite{dhr} and is in the spirit of Druet \cite{d3}, 
where the concentration analysis is studied in the intricate Riemannian 
setting. The study of the asymptotic for elliptic nonlinear pde's
was initiated by Atkinson-Peletier \cite{ap}, see also Br\'ezis-Peletier \cite{bp}. In the Riemannian context, the asymptotics have first been studied by Schoen \cite{s} and Hebey-Vaugon \cite{hv1}. This tool turned out to be very powerful
in the study of best constant problems in Sobolev inequalities, see
for instance Druet \cite{d1}, Hebey-Vaugon \cite{hv1}, \cite{hv2} and
Robert \cite{r2}). We also mention the study of the asymptotics for
solutions to nonlinear pde's (Han \cite{ha}, Hebey \cite{he}, 
Druet-Robert \cite{dr} and Robert \cite{r1}). In the case of arbitrary large energies, the compactness issues become quite intricate, especially in 
the Riemannian context, see for instance the pioneer work of Schoen 
\cite{s}. We also refer to the recent work of Druet \cite{d3} and Marques \cite{marques}. One can also find compactness results for fourth order 
equations  in the work of Hebey-Robert \cite{hr} and Hebey-Robert-Wen \cite{hrw}.
In a forthcoming paper \cite{gr1}, we tackle similar questions for
various critical equations involving a whole affine subspace of
singularities on the boundary.

The paper is organized as follows. In Section \ref{sec:prem}, we state
general facts and two lemmae that will be useful throughout the paper. In
Section \ref{sec:exh}, we construct the different scales of blow-up. In
Sections \ref{sec:fe:1} and \ref{sec:fe:2}, we prove strong pointwise
estimates for sequences of solutions to our problem. In Section
\ref{sec:poho}, we use the Pohozaev identity to describe precisely the
asymptotics related to our problem and we prove theorem \ref{th:cpct}.
Section \ref{sec:proof:th} contains the proofs of Theorems 1.1 and 1.2. 
Finally, we give in
the Appendix a regularity result for solutions to a critical PDE, some 
useful properties of the Green's function and a symmetry property of 
solutions to some nonlinear elliptic equations on the half-plane.

\section{Basic facts and preliminary Lemmae}\label{sec:prem}

Throughout the paper,  $\Omega$ will be a smooth bounded domain of $\rn$, $n\geq 
3$, such that $0\in\partial\Omega$. For $s\in (0,2)$, we write $\crit =\crit (s):=\frac{2(n-s)}{n-2}$
and for each $\eps>0$, we consider $\pe\in [0,\crit-2)$ such that
\bequa\label{lim:pe}
\lim_{\eps\to 0}\pe=0.
\eequa
We let $a\in C^1(\overline{\Omega})$ and a family $(\ae)_{\eps>0}\in 
C^1(\overline{\Omega})$ such that (\ref{hyp:ae}) holds. For any $\eps>0$, 
we consider $\ue\in \huno$ to be a solution to the
system
$$\left\{\begin{array}{ll}
\Delta\ue+\ae\ue=\frac{|\ue|^{\crit-2-\pe}}{|x|^s}\ue& \hbox{ in 
}{\mathcal D}'(\Omega)\\
\ue=0 &\hbox{ on }\partial\Omega
\end{array}\right.\eqno{(E_\eps)}$$
for all $\eps>0$. Note that it follows from Proposition \ref{prop:app} of
the Appendix that
$$\ue\in C^{1,\theta}(\overline{\Omega})\cap 
C^2(\overline{\Omega}\setminus\{0\})$$
for all $\theta\in (0,\min\{1,\crit-s\})$. In addition, we assume that
there exists $\Lambda>0$ such that
\bequa\label{bnd:ue}
\Vert\ue\Vert_{\huno}\leq\Lambda
\eequa
for all $\eps>0$. It then follows from the weak compactness of the
unit ball of $\huno$ that there exists $u_0\in\huno$ such that
\bequa\label{weak:lim:ue}
\ue\rightharpoonup u_0
\eequa
weakly in $\huno$ when $\eps\to 0$. Note that $u_0$ verifies

$$\Delta u_0+a u_0=\frac{|u_0|^{\crit-2}}{|x|^s}u_0 \hbox{ in }{\mathcal 
D}'(\Omega).$$
It follows from the Appendix that
$$u_0\in C^{1,\theta}(\overline{\Omega})\cap 
C^2(\overline{\Omega}\setminus\{0\})$$
for all $\theta\in (0,\min\{1,\crit-s\})$. The following Proposition 
addresses the case when $\ue$ is uniformly bounded in $L^\infty$. Note that here and in the sequel, all the convergence results are up to the extraction of a subsequence.

\begin{prop}\label{prop:bound} Let $\Omega$ be a smooth bounded domain of
$\rn$, $n\geq 3$, such that $0\in\partial\Omega$. We let $(\ue)$, $(\ae)$
and $(\pe)$ such that $(E_\eps)$, (\ref{hyp:ae}) and
(\ref{lim:pe}) hold. We assume that there exists $C>0$ such that
$|\ue(x)|\leq C$ for all $x\in\Omega$. Then up to a subsequence,
$\lim_{\eps\to 0}\ue=u_0$ in $C^1(\overline{\Omega})$, where $u_0$ is as 
in (\ref{weak:lim:ue}).
\end{prop}

\noindent{\it Proof:} It follows from the
proof of Proposition \ref{prop:app} of the Appendix that for any
$\theta\in (0,\min\{1,\crit-s\})$, there exists $C>0$ such that
$\Vert\ue\Vert_{C^{1,\theta}(\overline{\Omega})}\leq C$ for all $\eps>0$.
The conclusion of the Proposition then follows. We refer to the Appendix
for the details.\hfill$\Box$

\medskip\noindent From now on, we assume that

\bequa\label{hyp:blowup}
\lim_{\eps\to 0}\Vert\ue\Vert_{L^\infty(\Omega)}=+\infty.
\eequa
Throughout the paper, we shall say that blow-up occurs whenever 
(\ref{hyp:blowup})
holds. We define

$$\rnm=\{x\in\rn/\, x_1<0\}$$
where $x_1$ is the first coordinate of a generic point of $\rn$. This
space will be the limit space after blow-up. In the sequel of this
section, we give some useful tools for the blow-up analysis. We let
$y_0\in\partial\Omega$. Since $\partial\Omega$ is smooth and $y_0\in
\partial\Omega$, there exist $U,V$ open subsets of $\rn$, there exists $I$
an open intervall of $\rr$, there exists $U'$ an open subset of
$\rr^{n-1}$ such that $0\in U=I\times U'$ and $y_0\in V$. There exist
$\varphi\in C^\infty(U,V)$ and $\varphi_0\in C^\infty(U')$ such that, up
to rotating the coordinates if necessary,

\bequa\begin{array}{ll}

(i) & \varphi: U\to V\hbox{ is a }C^\infty-\hbox{diffeomorphism}\\

(ii) & \varphi(0)=y_0\\

(iii) & D_{0}\varphi=Id_{\rn}\\

(iv) & \varphi(U\cap\{x_1<0\})=\varphi(U)\cap\Omega\hbox{ and
}\varphi(U\cap\{x_1=0\})=\varphi(U)\cap\partial\Omega.\\

(v) & \varphi(x_1,y)=y_0+(x_1+\varphi_0(y),y)\hbox{ for all }(x_1,y)\in
I\times U'=U\\

(vi) & \varphi_0(0)=0\hbox{ and }\nabla\varphi_0(0)=0.

\end{array}\label{def:vphi}
\eequa
Here $D_x\varphi$ denotes the differential of $\varphi$ at $x$. This chart
will be useful throughout all the paper.

\medskip\noindent We prove two useful blow-up lemmae:

\begin{lem}\label{lem:blowup:1} We let $\Omega$ be a smooth bounded domain
of $\rn$, $n\geq 3$. We assume that $0\in\partial\Omega$. We let $(\ue)$,
$(\ae)$ and $(\pe)$ such that $(E_\eps)$, (\ref{hyp:ae}), (\ref{lim:pe}) 
and (\ref{bnd:ue}) hold. We let $(\ye)_{\eps>0}\in\Omega$. Let

$$\ne:=|\ue(\ye)|^{-\frac{2}{n-2}}\hbox{ and
}\be:=|\ye|^{\frac{s}{2}}|\ue(\ye)|^{\frac{2+\pe-\crit}{2}}.$$

We assume that $\lim_{\eps\to 0}\ne=0$. In particular, $\lim_{\eps\to
0}\be=0$. We assume that for any $R>0$, there exists $C(R)>0$ such that

\bequa\label{bnd:lem:1}
|\ue(x)|\leq C(R) |\ue(\ye)|
\eequa
for all $x\in B_{R\be}(\ye)\cap\Omega$ and all $\eps>0$. Then we have that

$$\ye=O\left(\ne^{1-\frac{\pe}{\crit-2}}\right)$$
when $\eps\to 0$. In particular, $\lim_{\eps\to 0}\ye=0$.

\end{lem}

\noindent{\it Proof of Lemma \ref{lem:blowup:1}:} We proceed by
contradiction and assume that

\bequa\label{lim:step:1}
\lim_{\eps\to 0}\frac{|\ye|}{\elle}=+\infty.
\eequa
where $\elle:=\ne^{1-\frac{\pe}{\crit-2}}$ for all $\eps>0$. In
particular, it follows from the definition of $\be$ and (\ref{lim:step:1})
that

\bequa\label{ppty:be}
\lim_{\eps\to 0}\be=0,\; \lim_{\eps\to 0}\frac{\be}{\elle}=+\infty\hbox{
and }\lim_{\eps\to 0}\frac{\be}{|\ye|}=0.
\eequa

\medskip\noindent{\it Case 1:} We assume that there exists $\rho>0$ such
that

$$\frac{d(\ye,\partial\Omega)}{\be}\geq 3\rho$$
for all $\eps>0$. For $x\in B_{2\rho}(0)$ and $\eps>0$, we define

$$\ve(x):=\frac{\ue(\ye+\be x)}{\ue(\ye)}.$$
Note that this is well defined since $\ye+\be x\in\Omega$ for all $x\in
B_{2\rho}(0)$. It follows from (\ref{bnd:lem:1}) that there exists 
$C(\rho)>0$ such that

\bequa\label{bnd:lem:1:bis}
|\ve(x)|\leq C(\rho)
\eequa
for all $\eps>0$ and all $x\in B_{2\rho}(0)$. As easily checked, we have
that

$$\Delta\ve+\be^2\ae(\ye+\be
x)\ve=\frac{|\ve|^{\crit-2-\pe}\ve}{\left|\frac{\ye}{|\ye|}+\frac{\be}{|\ye|}
x\right|^s}$$
weakly in $B_{2\rho}(0)$. Since (\ref{ppty:be}) holds, we have that

\bequa\label{eq:lem:1}
\Delta\ve+\be^2\ae(\ye+\be x)\ve=(1+o(1))|\ve|^{\crit-2-\pe}\ve
\eequa
weakly in $B_{2\rho}(0)$, where $\lim_{\eps\to 0}o(1)=0$ in
$C^0_{loc}(B_{2\rho}(0))$. It follows from (\ref{bnd:lem:1:bis}),
(\ref{eq:lem:1}) and standard elliptic theory that there exists $v\in
C^1(B_{2\rho}(0))$ such that

$$\ve\to v$$
in $C^1_{loc}(B_{2\rho}(0))$ when $\eps\to 0$. In particular,

\bequa\label{lim:v:case1:nonzero}
v(0)=\lim_{\eps\to 0}\ve(0)=1
\eequa
and $v\not\equiv 0$. With a change of variables and the definition of
$\be$, we get that

\beq
&&\int_{\Omega\cap B_{\rho\be}(\ye)}\frac{|\ue|^{\crit-\pe}}{|x|^s}\,
dx=\frac{|\ue(\ye)|^{\crit-\pe}\be^n}{|\ye|^s}\int_{B_{\rho}(0)}\frac{|\ve|^{\crit-\pe}}{\left|\frac{\ye}{|\ye|}+\frac{\be}{|\ye|}
x\right|^s}\, dx\\
&&\geq
\left(\frac{|\ye|}{\elle}\right)^{s\frac{n-2}{2}}\int_{B_{\rho}(0)}\frac{|\ve|^{\crit-\pe}}{\left|\frac{\ye}{|\ye|}+\frac{\be}{|\ye|}x\right|^s}\,
dx.
\eeq
Using the equation $(E_\eps)$, (\ref{bnd:ue}), (\ref{lim:step:1}) and 
(\ref{ppty:be}) and passing to the limit
$\eps\to 0$, we get that
$$\int_{B_{\rho}(0)}|v|^{\crit}\, dx=0,$$
and then $v\equiv 0$ in $B_\rho(0)$. A contradiction with
(\ref{lim:v:case1:nonzero}). Then (\ref{lim:step:1}) does not hold in Case
1.

\medskip\noindent{\it Case 2:} We assume that, up to a subsequence,
\bequa\label{lim:d:be:0}
\lim_{\eps\to 0}\frac{d(\ye,\partial\Omega)}{\be}=0.
\eequa
In this case,
$$\lim_{\eps\to 0}\ye=y_0\in\partial\Omega.$$
Since $y_0\in\partial\Omega$, we let $\varphi:U\to V$ as in
(\ref{def:vphi}), where $U,V$ are open neighborhoods of $0$ and $y_0$
respectively. We let $\tue=\ue\circ\varphi$, which is defined on $U\cap
\{x_1\leq 0\}$. For any $i,j=1,...,n$, we let
$g_{ij}=(\partial_i\varphi,\partial_j\varphi)$, where $(\cdot,\cdot)$
denotes the Euclidean scalar product on $\rn$, and we consider $g$ as a
metric on $\rn$. We let $\Delta_g=-div_g(\nabla)$ the Laplace-Beltrami
operator with respect to the metric $g$. In our basis, we have that

$$\Delta_g=-g^{ij}\left(\partial_{ij}-\Gamma_{ij}^k\partial_k\right),$$
where $g^{ij}=(g^{-1})_{ij}$ are the coordinates of the inverse of the
tensor $g$ and the $\Gamma_{ij}^k$'s are the Christoffel symbols of the
metric $g$. As easily checked, we have that

$$\Delta_g\tue+\ae\circ\varphi(x)\cdot
\tue=\frac{|\tue|^{\crit-2-\pe}\tue}{|\varphi(x)|^s}$$
weakly in $U\cap\{x_1<0\}$. We let $\ze\in\partial\Omega$ such that

\bequa\label{def:ze}
|\ze-\ye|=d(\ye,\partial\Omega).
\eequa
We let $\tye,\tze\in U$ such that

\bequa\label{def:tye:tze}
\varphi(\tye)=\ye\hbox{ and }\varphi(\tze)=\ze.
\eequa
It follows from the properties (\ref{def:vphi}) of $\varphi$ that

\bequa\label{ppty:tye:tze}
\lim_{\eps\to 0}\tye=\lim_{\eps\to 0}\tze=0,\; (\tye)_1<0\hbox{ and
}(\tze)_1=0.\eequa
At last, we let

$$\tve(x):=\frac{\tue(\tze+\be x)}{\tue(\tye)}$$
for all $x\in \frac{U-\tze}{\be}\cap \{x_1<0\}$. With
(\ref{ppty:tye:tze}), we get that $\tve$ is defined on
$B_R(0)\cap\{x_1<0\}$ for all $R>0$, as soon as $\eps$ is small enough. It
follows from (\ref{bnd:lem:1}) that there exists $C'(R)>0$ such that

\bequa\label{bnd:lem:1:ter}
|\tve(x)|\leq C'(R)
\eequa
for all $\eps>0$ and all $x\in B_{R}(0)\cap\{x_1\leq 0\}$. The function
$\tve$ verifies

$$\Delta_{\tge}\tve+\be^2\ae\circ\varphi(\tze+\be
x)\tve=\frac{|\tve|^{\crit-2-\pe}\tve}{\left|\frac{\varphi(\tze+\be
x)}{|\ye|}\right|^s}$$
weakly in $B_R(0)\cap\{x_1<0\}$. In this expression, $\tge=g(\tze+\be x)$
and $\Delta_{\tge}$ is the Laplace-Beltrami operator with respect to the
metric $\tge$. With (\ref{lim:d:be:0}), (\ref{def:ze}) and
(\ref{def:tye:tze}), we get that

$$\varphi(\tze+\be x)=\ye+O_R(1)\be,$$
for all $x\in B_{R}(0)\cap\{x_1\leq 0\}$ and all $\eps>0$, where there
exists $C_R>0$ such that $|O_R(1)|\leq C_R$ for all $x\in
B_{R}(0)\cap\{x_1<0\}$. With (\ref{ppty:be}), we then get that

$$\lim_{\eps\to 0}\frac{|\varphi(\tze+\be x)|}{|\ye|}=1$$
in $C^0(B_{R}(0)\cap\{x_1\leq 0\})$. It then follows that

$$\Delta_{\tge}\tve+\be^2\ae\circ\varphi(\tze+\be
x)\tve=(1+o(1))|\tve|^{\crit-2-\pe}\tve$$
weakly in $B_R(0)\cap\{x_1<0\}$, where $\lim_{\eps\to 0}o(1)=0$ in
$C^0(B_R(0)\cap\{x_1\leq 0\})$. Since $\tve$ vanishes on
$B_R(0)\cap\{x_1=0\}$ and (\ref{bnd:lem:1:ter}) holds, it follows
from standard elliptic theory that there exists $\tv\in 
C^1(B_R(0)\cap\{x_1\leq 0\})$ such that
$$\lim_{\eps\to 0}\tve=\tv$$
in $C^{0}(B_{\frac{R}{2}}(0)\cap \{x_1\leq 0\})$. In particular,

\bequa\label{eq:v:vanish}
\tv\equiv 0\hbox{ on }B_{\frac{R}{2}}(0)\cap \{x_1=0\}.
\eequa
Moreover, it follows from (\ref{lim:d:be:0}), (\ref{def:ze}) and
(\ref{def:tye:tze}) that

$$\tve\left(\frac{\tye-\tze}{\be}\right)=1\hbox{ and }\lim_{\eps\to
0}\frac{\tye-\tze}{\be}=0.$$
In particular, $\tv(0)=1$. A contradiction with (\ref{eq:v:vanish}). Then
(\ref{lim:step:1}) does not hold in Case 2.

\smallskip\noindent In both cases, we have contradicted
(\ref{lim:step:1}). This proves that $\ye=O(\elle)$ when $\eps\to 0$,
which proves the Lemma.\hfill$\Box$

\begin{lem}\label{lem:blowup:2}  We let $\Omega$ be a smooth bounded
domain of $\rn$, $n\geq 3$. We assume that $0\in\partial\Omega$. We let
$(\ue)$, $(\ae)$ and $(\pe)$ such that $(E_\eps)$, (\ref{hyp:ae}), 
(\ref{lim:pe}) and (\ref{bnd:ue}) hold. We let $(\ne)_{\eps>0}$ and
$(\elle)_{\eps>0}$ such that $\ne,\elle>0$ for all $\eps>0$ and

$$\elle=\ne^{1-\frac{\pe}{\crit-2}}\hbox{ and }\lim_{\eps\to 0}\ne=0.$$
Since $0\in\partial\Omega$, we let $\varphi:U\to V$ as in (\ref{def:vphi})
with $y_0=0$, where $U,V$ are open neighborhoods of $0$. We let

$$\tue(x):=\ne^{\frac{n-2}{2}}\ue\circ\varphi(\elle x)$$
for all $x\in \frac{U}{\elle}\cap \{x_1\leq 0\}$ and all $\eps>0$. We
assume that either

\smallskip (L1) for all $R>0$, there exists $C(R)>0$ such that
$$|\tue(x)|\leq C(R)$$
for all $x\in B_R(0)\cap \{x_1<0\}$, or\par

\smallskip (L2) for all $R>\delta>0$, there exists $C(R,\delta)>0$ such
that

$$|\tue(x)|\leq C(R,\delta)$$
for all $x\in \left(B_R(0)\setminus \overline{B}_\delta(0)\right)\cap
\{x_1<0\}$.\par

\smallskip\noindent Then there exists $\tu\in\hunrnm\cap
C^{1}(\overline{\rnm})$ such that

$$\Delta\tu=\frac{|\tu|^{\crit-2}\tu}{|x|^s}\hbox{ in }{\mathcal
D}'(\rnm)$$
and
$$\lim_{\eps\to 0}\tue=\tu\hbox{ in }\left\{\begin{array}{ll}
C^{1}_{loc}(\overline{\rnm})&\hbox{if (L1) holds}\\
C^{1}_{loc}(\overline{\rnm}\setminus\{0\})&\hbox{if (L2) holds}
\end{array}\right.$$
\end{lem}

\noindent{\it Proof of Lemma \ref{lem:blowup:2}:} Let $\eta\in
C^\infty(\rn)$. As easily checked, we have that

$$\eta\tue\in \hunrnm$$
for all $\eps>0$ small enough, and

$$\nabla(\eta\tue)(x)=\tue\nabla\eta+\eta
\elle\ne^{\frac{n-2}{2}}D_{(\elle
x)}\varphi[(\nabla\ue)(\varphi(\elle x))],$$
for all $\eps>0$ and all $x\in\rnm$. In this expression, $D_x\varphi$ is
the differential of the function $\varphi$ at $x$. We get that

\beq
&&\int_{\rnm}|\nabla(\eta\tue)|^2\, dx\leq
2\int_{\rnm}|\nabla\eta|^2\tue^2\, dx\\
&&+2 \elle^2\ne^{n-2}\int_{\rnm\cap \hbox{Supp
}\eta}|D_{(\elle x)}\varphi[(\nabla\ue)(\varphi(\elle x))]|^2\, dx.
\eeq
With H\"older's inequality and a change of variables, we get that

\beqn
&&\int_{\rnm} |\nabla(\eta\tue)|^2\, dx\leq
2\left(\int_{\rnm}|\nabla\eta|^n\,
dx\right)^{\frac{2}{n}}\cdot\left(\int_{\rnm\cap\hbox{Supp
}\nabla\eta}|\tue|^{\frac{2n}{n-2}}\, dx\right)^{\frac{n-2}{n}}\nonumber\\
&&+4\elle^n\left(\frac{\ne}{\elle}\right)^{n-2}\int_{\rnm\cap\hbox{Supp
}\eta}|\nabla\ue|^2(\varphi(\elle x))\, dx\nonumber\\
&&\leq
2\Vert\nabla\eta\Vert_{n}^2\Vert\tue\Vert_{L^{\frac{2n}{n-2}}(\hbox{Supp
}\nabla\eta)}^{2}\nonumber\\
&&+C \ne^{\frac{\pe(n-2)}{\crit-2}}\int_\Omega|\nabla\ue|^2\,
dx\label{ineq:ve:1}
\eeqn
With another change of variables, we get that

\beqn
&&\int_{\rnm} |\nabla(\eta\tue)|^2\, dx\leq C
\ne^{\frac{(n-2)\pe}{\crit-2}}\Vert\nabla\eta\Vert_{n}^2\Vert\ue\Vert_{L^{\frac{2n}{n-2}}(\Omega)}^{2}\nonumber\\
&&+C \ne^{\frac{\pe(n-2)}{\crit-2}}\int_\Omega|\nabla\ue|^2\,
dx\label{ineq:ve:2}
\eeqn
for all $\eps>0$, where $C$ is independant of $\eps$. With (\ref{bnd:ue}), 
Sobolev's inequality and since $\ne^{\pe}\leq 1$ for
all $\eps>0$ small enough, we get with (\ref{ineq:ve:2}) that

$$\Vert\eta\tue\Vert_{\hunrnm}=O(1)$$
when $\eps\to 0$. It then follows that there exists
$\tilde{u}_\eta\in\hunrnm$ such that, up to a subsequence,

$$\eta\tue\rightharpoonup \tilde{u}_\eta$$
weakly in $\hunrnm$ when $\eps\to 0$. We let $\eta_1\in C_c^\infty(\rn)$
such that $\eta_1\equiv 1$ in $B_1(0)$ and $\eta_1\equiv 0$ in
$\rn\setminus \overline{B}_2(0)$. For any $R\in\nn^\star$, we let
$\eta_R(x)=\eta_1(\frac{x}{R})$ for all $x\in\rn$. With a diagonal
argument, we can assume that, up to a subsequence, for any $R>0$, there
exists $\tu_R\in\hunrnm$ such that

$$\eta_R\tue\rightharpoonup \tu_R$$
weakly in $\hunrnm$ when $\eps\to 0$. Letting $\eps\to 0$ in
(\ref{ineq:ve:2}), with (\ref{bnd:ue}), Sobolev's inequality and since
$\ne^{\pe}\leq 1$ for all $\eps>0$ small enough, we get that there exists
a constant $C>0$ independant of $R$ such that

$$\int_{\rnm} |\nabla \tu_R|^2\, dx\leq  C\Vert\nabla\eta_R\Vert_{n}^2+ C$$
for all $R>0$. Since
$\Vert\nabla\eta_R\Vert_{n}^2=\Vert\nabla\eta_1\Vert_{n}^2$ for all $R>0$,
we get that there exists $C>0$ independant of $R$ such that

$$\int_{\rnm} |\nabla \tu_R|^2\, dx\leq C$$
for all $R>0$. It then follows that there exists $\tu\in \hunrnm$ such
that $\tu_R\rightharpoonup \tu$ weakly in $\hunrnm$ when $R\to +\infty$. As easily checked, we then obtain that $\tu_\eta=\eta \tu$ (we omit the proof
of this fact. It is straightforward).

\medskip\noindent For any $i,j=1,...,n$, we let
$(\tge)_{ij}=(\partial_i\varphi(\elle x),\partial_j\varphi(\elle x))$,
where $(\cdot,\cdot)$ denotes the Euclidean scalar product on $\rn$. We
consider $\tge$ as a metric on $\rn$. We let

$$\Delta_{\tge}=-{\tge}^{ij}\left(\partial_{ij}-\Gamma_{ij}^k(\tge)\partial_k\right),$$
where $\tge^{ij}:=(\tge^{-1})_{ij}$ are the coordinates of the inverse of
the tensor $\tge$ and the $\Gamma_{ij}^k(\tge)$'s are the Christoffel
symbols of the metric $\tge$. With a change of variable, equation
$(E_\eps)$ rewrites as

\bequa\label{eq:ve:sec3}
\Delta_{\tge}(\eta_R\tue)+\elle^2\ae\circ\varphi(\elle
x)\eta_R\tue=\frac{|\eta_R\tue|^{\crit-2-\pe}\eta_R\tue}{\left|\frac{\varphi(\elle
x)}{\elle}\right|^s}\hbox{ in }{\mathcal D}'(B_R(0)\cap\{x_1<0\})
\eequa
for all $\eps>0$. Passing to the weak limit $\eps\to 0$ and then $R\to 
+\infty$ in this equation, we get that

$$\Delta\tu=\frac{|\tu|^{\crit-2}\tu}{|x|^s}\hbox{ in }{\mathcal
D}'(\rnm).$$
Since $\tu\in \hunrnm$, it follows from Proposition \ref{prop:app} of the Appendix that $\tu\in 
C^{1,\theta}(\overline{\rnm})$ for all $\theta\in(0,\min\{1,\crit-s\})$.

\medskip\noindent We deal with case (L1). Since $s\in (0,2)$, (L1) and
(\ref{eq:ve:sec3}) hold and $\tue\equiv 0$ on $\{x_1=0\}$, it follows from 
arguments similar to the ones
developed in the Appendix that for any $\theta\in (0,\min\{1,\crit-s\})$
and any $R>0$, there exists $C(\theta,R)>0$ independant of $\eps>0$ small
such that

$$\Vert\tue\Vert_{C^{1,\theta}(B_R(0)\cap\{x_1\leq 0\})}\leq C(\theta,R)$$
for all $\eps>0$ small. It then follows from Ascoli's theorem that for any
$\theta\in (0,\min\{1,\crit-s\})$,

$$\lim_{\eps\to 0}\tue=\tu$$
in $C^{1,\theta}_{loc}(\overline{\rnm})$. The proof proceeds similarly in
Case (L2). This ends the proof of the Lemma.\hfill$\Box$

\section{Construction and exhaustion of the blow-up scales}\label{sec:exh}
This section is devoted to the proof of the following proposition:
\begin{prop}\label{prop:exhaust} We let $\Omega$ be a smooth bounded
domain of $\rn$, $n\geq 3$. We assume that $0\in\partial\Omega$. We let
$(\ue)$, $(\ae)$ and $(\pe)$ such that $(E_\eps)$, (\ref{hyp:ae}), 
(\ref{lim:pe}) and (\ref{bnd:ue}) hold. We assume that blow-up occurs,
that is

$$\lim_{\eps\to 0}\Vert\ue\Vert_{L^\infty(\Omega)}=+\infty.$$
Then there exists $N\in\nn^\star$, there exists $N$ families of points 
$(\mei)_{\eps>0}$ such that we have that\par

\smallskip\noindent {\bf(A1)} $\lim_{\eps\to 0}\ue=u_0$ in
$C^2_{loc}(\overline{\Omega}\setminus\{0\})$ where $u_0$ is as in
(\ref{weak:lim:ue}),

\smallskip\noindent {\bf(A2)} $0< \meun< ...<\meN$ for all $\eps>0$,\par

\smallskip\noindent {\bf(A3)} $$\lim_{\eps\to 0}\meN=0\hbox{ and

}\lim_{\eps\to 0}\frac{\mu_{\eps,i+1}}{\mei}=+\infty\hbox{ for all

}i=1...N-1$$

\smallskip\noindent {\bf(A4)} For all $i=1...N$, there exists $\tui\in
\hunrnm\cap C^{1}(\overline{\rnm})\setminus\{0\}$ such that

$$\Delta \tui=\frac{|\tui|^{\crit-2}\tui}{|x|^s}\hbox{ in }{\mathcal
D}'(\rnm)$$
and

$$\lim_{\eps\to 0}\tuei=\tui$$
in $C^{1}_{loc}(\overline{\rnm}\setminus\{0\})$, where

$$\tuei(x):=\mei^{\frac{n-2}{2}}\ue(\varphi(\kei x))$$
for all $x\in \frac{U}{\kei}\cap\{x_1\leq 0\}$ and 
$\kei=\mei^{1-\frac{\pe}{\crit-2}}$. Moreover, $\lim_{\eps\to
0}\tueun=\tu_1$ in $C^{1}_{loc}(\overline{\rnm})$.\par

\smallskip\noindent {\bf(A5)}

$$\lim_{R\to +\infty}\lim_{\eps\to 0}\sup_{|x|\geq R
\keN}|x|^{\frac{n-2}{2}}|\ue(x)-u_0(x)|^{1-\frac{\pe}{\crit-2}}=0$$

\smallskip\noindent {\bf(A6)} For any $\delta>0$ and any $i=1...N-1$, we have
that

$$\lim_{R\to +\infty}\lim_{\eps\to 0}\sup_{\delta k_{\eps,i+1}\geq |x|\geq R
\kei}|x|^{\frac{n-2}{2}}\left|\ue(x)-\mu_{\eps,i+1}^{-\frac{n-2}{2}}u_{i+1}\left(\frac{\varphi^{-1}(x)}{k_{\eps,i+1}}\right)\right|^{1-\frac{\pe}{\crit-2}}=0.$$
\smallskip\noindent {\bf(A7)} For any $i\in\{1,...,N\}$, there exists
$\alpha_i\in (0,1]$ such that

$$\lim_{\eps\to 0}\mei^{\pe}=\alpha_i.$$
\end{prop}

\noindent The proof of this proposition proceeds in seven steps.

\medskip\noindent{\bf Step \ref{sec:exh}.1:} We let $\xeun\in\Omega$ and
$\meun,\keun>0$ such that

\bequa\label{def:me:xe}
\max_\Omega|\ue|=|\ue(\xeun)|=\meun^{-\frac{n-2}{2}}\hbox{ and
}\keun=\meun^{1-\frac{\pe}{\crit-2}}.
\eequa
We claim that
\bequa\label{lim:xe:ke}
|\xeun|=O(\keun)
\eequa
when $\eps\to 0$, and in particular that $\lim_{\eps\to 0}\xe=0$. Indeed,
we use Lemma \ref{lem:blowup:1} with $\ye=\xeun$, $\ne=\meun$ and
$C(R)=1$. We then immediately get that $|\xeun|=O(\keun)$ when $\eps\to
0$.

\medskip\noindent From now on, we let $\varphi:U\to V$ as in
(\ref{def:vphi}) with $y_0=0$ and $U,V$ are open neighborhoods of $0$ in
$\rn$. We then let

\bequa\label{def:ae:be}
\xeun=\varphi(a_\eps,b_\eps),
\eequa
where $a_\eps\in \{x_1<0\}$, $b_\eps\in\rr^{n-1}$ and $(a_\eps,b_\eps)\in
U$. Note that $\lim_{\eps\to 0}(a_\eps,b_\eps)=(0,0)$.

\medskip\noindent{\bf Step \ref{sec:exh}.2:} We claim that

\bequa\label{lim:d:ze}
d(\xeun,\partial\Omega)=(1+o(1))|a_\eps|=O(\keun)
\eequa
where $\lim_{\eps\to 0}o(1)=0$.\par

\smallskip\noindent{\it Proof of the Claim:} Indeed, since 
$0\in\partial\Omega$, we get with (\ref{lim:xe:ke}) that

\bequa\label{lim:prop32:1}
d(\xeun,\partial\Omega)\leq |\xeun-0|=O(\keun)
\eequa
when $\eps\to 0$. We first remark that

$$d(\xeun,\partial\Omega)\leq d(\xeun,\varphi(0,b_\eps))=|a_\eps|.$$
We let $\gamma_\eps\in\rr^{n-1}$ such that $(0,\gamma_\eps)\in
U\cap\{x_1=0\}$ and $Y_\eps=\varphi(0,\gamma_\eps)\in\partial\Omega$ such
that $d(\xeun,\partial\Omega)=|\xeun-Y_\eps|$. Since
$d(\xeun,\partial\Omega)\leq |a_\eps|$, we get that

$$b_\eps-\gamma_\eps=O(|a_\eps|),$$
when $\eps\to 0$. Since $\nabla\varphi_0(0)=0$ (where $\varphi_0$ is as in
(\ref{def:vphi})), we get that

$$\varphi_0(b_\eps)=\varphi_0(\gamma_\eps)+o(|b_\eps-\gamma_\eps|)=\varphi_0(\gamma_\eps)+o(|a_\eps|)$$
when $\eps\to 0$. Moreover,

\beq
d(\xeun,\partial\Omega)&=&|\xeun-Y_\eps|\\
&=& |(a_\eps+\varphi_0(b_\eps)-\varphi_0(\gamma_\eps),
b_\eps-\gamma_\eps)|\\
&=& |(a_\eps+o(a_\eps), b_\eps-\gamma_\eps)|\leq |a_\eps|
\eeq
when $\eps\to 0$. It then follows that $b_\eps-\gamma_\eps=o(|a_\eps|)$
and $d(\xeun,\partial\Omega)=(1+o(1))|a_\eps|$ when $\eps\to 0$. This
prove (\ref{lim:d:ze}).\hfill$\Box$

\medskip\noindent The classical Hardy-Sobolev inequality asserts that
there exists $C>0$ such that

\bequa\label{ineq:HS:rn}
\left(\int_{\rn}\frac{|u|^{\crit}}{|x|^s}\,
dx\right)^{\frac{2}{\crit}}\leq C \int_{\rn}|\nabla u|^2\, dx
\eequa
for all $u\in \hunrn$. We define

\bequa\label{def:best:cst}
\mu_s(\rnm):=\inf\frac{\int_{\rnm}|\nabla u|^2\,
dx}{\left(\int_{\rnm}\frac{|u|^{\crit}}{|x|^s}\,
dx\right)^{\frac{2}{\crit}}}\
\eequa
where the infimum is taken over functions $u\in\hunrnm\setminus\{0\}$. The
existence of $\mu_s(\rnm)>0$ is a consequence of (\ref{ineq:HS:rn}).

\medskip\noindent{\bf Step \ref{sec:exh}.3:} The construction of the
$(\mei)$'s proceeds by induction. This step is the initiation.

\begin{lem}\label{lem:init} We let
$$\tueun(x):=\meun^{\frac{n-2}{2}}\ue\circ\varphi(\keun x)$$
for all $\eps>0$ and all $x\in\frac{U}{\keun}\cap \{x_1\leq 0\}$. Then,
there exists $\tu_1\in \hunrnm\cap C^1(\overline{\rnm})$ such that

\noindent{\bf(B1)} $\lim_{\eps\to 0}\tueun=\tu_1$ in
$C^1_{loc}(\overline{\rnm})$,\par

\noindent{\bf(B2)}
$$\Delta\tu_1=\frac{|\tu_1|^{\crit-2}\tu_1}{|x|^s}\hbox{ in }{\mathcal
D}'(\rnm),$$

\noindent{\bf(B3)} $$\int_{\rnm}|\nabla \tu_1|^2\, dx\geq
\mu_s(\rnm)^{\frac{\crit}{\crit-2}}.$$
Moreover, there exists $\alpha_1\in (0,1]$ such that $\lim_{\eps\to
0}\meun^{\pe}=\alpha_1$.
\end{lem}

\noindent{\it Proof of Lemma \ref{lem:init}:} Indeed, since
$|\tueun(x)|\leq 1$ for all $x\in\frac{U}{\keun}\cap \{x_1\leq 0\}$,
hypothesis (L1) of Lemma \ref{lem:blowup:2} is satisfied and it follows
from Lemma \ref{lem:blowup:2} that points (B1) and (B2) hold. We let
$\lambda_\eps=-\frac{a_\eps}{\keun}>0$ and
$\theta_\eps=\frac{b_\eps}{\keun}\in\rr^{n-1}$, where $a_\eps,b_\eps$ are
defined in (\ref{def:ae:be}). It follows from Steps \ref{sec:exh}.1 and
\ref{sec:exh}.2 that there exists $\lambda_0\geq 0$ and
$\theta_0\in\rr^{n-1}$ such that $\lim_{\eps\to
0}(\lambda_\eps,\theta_\eps)=(\lambda_0,\theta_0)$. It then follows from
the definition of $\tueun$ and (\ref{def:me:xe}) that

$$|\tueun(-\lambda_\eps,\theta_\eps)|=1.$$
Passing to the limit $\eps\to 0$ and using point (B1), we get that
$|\tu_1(-\lambda_0,\theta_0)|=1$. In particular $\tu_1\not\equiv 0$ and 
$\lambda_0\neq 0$. Multiplying (B2) by $\tu_1$ and integrating by parts 
over $\rnm$, we get
that

$$\int_{\rnm}|\nabla\tu_1|^2\, dx=\int_{\rnm}\frac{|u|^{\crit}}{|x|^s}\,
dx.$$
Using the Hardy-Sobolev inequality (\ref{def:best:cst}) and that
$\tu_1\not\equiv 0$, we get (B3). At last, with (\ref{bnd:ue}),
(\ref{ineq:ve:2}) and Sobolev's inequality, we get that for any $\eta\in
C^\infty_c(\rn)$, there exists $C>0$ such that

$$\int_{\rnm} |\nabla(\eta\tueun)|^2\, dx\leq C
\meun^{\frac{(n-2)\pe}{\crit-2}}$$
for all $\eps>0$. Letting $\eps\to 0$ and using that $\tu_1\not\equiv 0$,
we get that $\lim_{\eps\to 0}\meun^{\pe}>0$.\hfill$\Box$

\medskip\noindent{\bf Step \ref{sec:exh}.4:} We claim that there exists
$C>0$ such that
\bequa\label{ineq:est:1}
|x|^{\frac{n-2}{2}}|\ue(x)|^{1-\frac{\pe}{\crit-2}}\leq C
\eequa
for all $\eps>0$ and all $x\in\Omega$.\par

\smallskip\noindent{\it Proof of the Claim:} We argue by contradiction and
we let $(\ye)_{\eps>0}\in\Omega$ such that

\bequa\label{hyp:step34}
\sup_{x\in\Omega}|x|^{\frac{n-2}{2}}|\ue(x)|^{1-\frac{\pe}{\crit-2}}=|\ye|^{\frac{n-2}{2}}|\ue(\ye)|^{1-\frac{\pe}{\crit-2}}\to
+\infty
\eequa
when $\eps\to 0$. We let

$$\ne:=|\ue(\ye)|^{-\frac{2}{n-2}}\hbox{ and
}\elle:=\ne^{1-\frac{\pe}{\crit-2}}$$
for all $\eps>0$. It follows from (\ref{hyp:step34}) that

\bequa\label{lim:infty:34}
\lim_{\eps\to 0}\frac{|\ye|}{\elle}=+\infty\hbox{ and }\lim_{\eps\to
0}\ne=0.
\eequa
We let
$$\be:=|\ye|^\frac{s}{2}|\ue(\ye)|^{\frac{2+\pe-\crit}{2}}.$$
It follows from (\ref{hyp:step34}) that

\bequa\label{lim:proof:34}
\lim_{\eps\to 0}\frac{\be}{|\ye|}=0.
\eequa
We let $R>0$. We let $x\in B_R(0)$ such that $\ye+\be x\in\Omega$. It
follows from the definition (\ref{hyp:step34}) of $\ye$ that

$$|\ye+\be x|^{\frac{n-2}{2}}|\ue(\ye+\be x)|\leq
|\ye|^{\frac{n-2}{2}}|\ue(\ye)|,$$
and then

$$\left(\frac{|\ue(\ye+\be
x)|}{|\ue(\ye)|}\right)^{1-\frac{\pe}{\crit-2}}\leq
\left(\frac{1}{1-\frac{\be}{|\ye|}R}\right)^{\frac{n-2}{2}}$$
for all $\eps>0$ and all $x\in B_R(0)$ such that $\ye+\be x\in\Omega$.
With (\ref{lim:proof:34}), we get that there exists $\eps(R)>0$ such that

$$|\ue(\ye+\be x)|\leq 2 |\ue(\ye)|$$
for all $x\in B_R(0)$ such that $\ye+\be x\in\Omega$ and all
$0<\eps<\eps(R)$. It then follows from Lemma \ref{lem:blowup:1} that
$\ye=O(\elle)$ when $\eps\to 0$. A contradiction with
(\ref{lim:infty:34}). This proves (\ref{ineq:est:1}).\hfill$\Box$

\medskip\noindent As a remark, it follows from $(E_\eps)$,
(\ref{weak:lim:ue}), (\ref{ineq:est:1}) and standard elliptic theory that

\bequa\label{cv:out:0}
\lim_{\eps\to 0}\ue=u_0\hbox{ in
}C^2_{loc}(\overline{\Omega}\setminus\{0\}).
\eequa

\medskip\noindent We let $p\in\mathbb{N}^\star$. We consider the following
assertions:

\smallskip\noindent{\bf(C1)} $0< \meun< ...<\mu_{\eps,p}$\par

\smallskip\noindent{\bf(C2)} $$\lim_{\eps\to 0}\mu_{\eps,p}=0\hbox{ and

}\lim_{\eps\to 0}\frac{\mu_{\eps,i+1}}{\mei}=+\infty\hbox{ for all

}i=1...p-1$$

\smallskip\noindent{\bf(C3)} For all $i=1...p$, there exists $\tui\in
\hunrnm\cap C^{1}(\overline{\rnm})\setminus\{0\}$ such that

$$\Delta \tui=\frac{|\tui|^{\crit-2}\tui}{|x|^s}\hbox{ in }{\mathcal
D}'(\rnm), \qquad\int_{\rnm}|\nabla \tui|^2\, dx\geq
\mu_s(\rnm)^{\frac{\crit}{\crit-2}}$$
and

$$\lim_{\eps\to 0}\tuei=\tui$$
in $C^{1}_{loc}(\overline{\rnm}\setminus\{0\})$, where

$$\tuei(x):=\mei^{\frac{n-2}{2}}\ue(\varphi(\kei x))$$
for all $x\in \frac{U}{\kei}\cap\{x_1\leq 0\}$ and 
$\kei:=\mei^{1-\frac{\pe}{\crit-2}}$.\par

\smallskip\noindent{\bf (C4)} For any $i\in\{1,...,p\}$, there exists
$\alpha_i\in (0,1]$ such that

$$\lim_{\eps\to 0}\mei^{\pe}=\alpha_i.$$
We say that ${\mathcal H}_p$ holds if there exists $p$ families of points
$(\mei)_{\eps>0}$, $i=1,...,p$ such that $(\meun)_{\eps>0}$ is as in
(\ref{def:me:xe}) and points (C1), (C2) (C3) and (C4) hold. Note that it
follows from Step \ref{sec:exh}.4 that ${\mathcal H}_1$ holds with the
improvement that the convergence in (C3) holds in
$C^{1}_{loc}(\overline{\rnm})$.

\medskip\noindent{\bf Step \ref{sec:exh}.5:} We prove the following
proposition:

\begin{prop}\label{prop:HN}  Let $\Omega$ be a smooth bounded domain of
$\rn$, $n\geq 3$, such that $0\in\partial\Omega$. We let $(\ue)$, $(\ae)$
and $(\pe)$ such that $(E_\eps)$, (\ref{hyp:ae}), (\ref{lim:pe}) and 
(\ref{bnd:ue}) hold. Let $p\geq 1$. We assume that ${\mathcal H}_p$ holds.
Then either

$$\lim_{R\to +\infty}\lim_{\eps\to 0}\sup_{|x|\geq R 
k_{\eps,p}}|x|^{\frac{n-2}{2}}|\ue(x)-u_0(x)|^{1-\frac{\pe}{\crit-2}}=0$$
or ${\mathcal H}_{p+1}$ holds.
\end{prop}

\noindent{\it Proof of Proposition \ref{prop:HN}:} We assume that

$$\lim_{R\to +\infty}\lim_{\eps\to 0}\sup_{|x|\geq R
k_{\eps,p}}|x|^{\frac{n-2}{2}}|\ue(x)-u_0(x)|^{1-\frac{\pe}{\crit-2}}\neq 
0.$$
It then follows that there exists a family $(\ye)_{\eps>0}\in\Omega$ such
that

\bequa\label{hyp:lim:ye:1}
\lim_{\eps\to 0}\frac{|\ye|}{k_{\eps,p}}=+\infty\hbox{ and }\lim_{\eps\to 
0}|\ye|^{\frac{n-2}{2}}|\ue(\ye)-u_0(\ye)|^{1-\frac{\pe}{\crit-2}}=\alpha>0.
\eequa
We claim that $\lim_{\eps\to 0}\ye=0$. Otherwise, it follows from
(\ref{cv:out:0}) that $\lim_{\eps\to 0}|\ue(\ye)-u_0(\ye)|=0$. A 
contradiction.

\medskip\noindent Since $u_0\in C^0(\overline{\Omega})$ and $\lim_{\eps\to
0}\ye=0$, we get that

\bequa\label{hyp:lim:ye}
\lim_{\eps\to 
0}|\ye|^{\frac{n-2}{2}}|\ue(\ye)|^{1-\frac{\pe}{\crit-2}}=\alpha>0.
\eequa
In particular, $\lim_{\eps\to 0}|\ue(\ye)|=+\infty$. We let

$$\mu_{\eps,p+1}:=|\ue(\ye)|^{-\frac{2}{n-2}}\hbox{ and 
}k_{\eps,p+1}:=\mu_{\eps,p+1}^{1-\frac{\pe}{\crit-2}}.$$
As a consequence, $\lim_{\eps\to 0}\mu_{\eps,p+1}=0$. We define

$$\tu_{\eps,p+1}(x):=\mu_{\eps,p+1}^{\frac{n-2}{2}}\ue(\varphi(k_{\eps,p+1}
x))$$
for all $x\in\frac{U}{k_{\eps,p+1}}\cap\{x_1\leq 0\}$. It follows from 
(\ref{ineq:est:1}) that
$$|\varphi(k_{\eps,p+1} x)|^{\frac{n-2}{2}}|\ue(\varphi(k_{\eps,p+1}
x))|^{1-\frac{\pe}{\crit-2}}\leq C$$
for all $x\in\frac{U}{k_{\eps,p+1}}\cap\{x_1\leq 0\}$. With the definition 
of $\tu_{\eps,p+1}$ and the properties (\ref{def:vphi}) of $\varphi$, we
get that there exists $C>0$ such that

$$|x|^{\frac{n-2}{2}}|\tu_{\eps,p+1}(x)|^{1-\frac{\pe}{\crit-2}}\leq C$$
for all $x\in\frac{U}{k_{\eps,p+1}}\cap\{x_1\leq 0\}$. It then follows
that hypothesis (L2) of Lemma \ref{lem:blowup:2} is satisfied. It then
follows from Lemma \ref{lem:blowup:2} that there exists $\tu_{p+1}\in
\hunrnm\cap C^{1}(\overline{\rnm})$ such that

$$\Delta \tu_{p+1}=\frac{|\tu_{p+1}|^{\crit-2}\tu_{p+1}}{|x|^s}\hbox{ in
}{\mathcal D}'(\rnm), $$
and

\bequa\label{cv:tup1}
\lim_{\eps\to 0}\tu_{\eps,p+1}=\tu_{\eps,p+1}
\eequa
in $C^{1}_{loc}(\overline{\rnm}\setminus\{0\})$. It follows from
(\ref{hyp:lim:ye}) and the definition of $k_{\eps,p+1}$ that

$$\lim_{\eps\to 0}\frac{|\ye|}{k_{\eps,p+1}}=\alpha>0.$$
We let $\tye\in \{x_1<0\}$ such that $\ye=\varphi(k_{\eps,p+1}\tye)$. It
then exists $\tilde{y}_0\in\overline{\rnm}$ such that $\lim_{\eps\to
0}\tye=\tilde{y}_0\neq 0$. It then follows from (\ref{cv:tup1}) that

$$|\tu_{p+1}(y_0)|=\lim_{\eps\to 0}|\tu_{\eps,p+1}(\tye)|=1,$$
and then $\tu_{p+1}\not\equiv 0$. With arguments similar to the ones
developed in the proof of Lemma \ref{lem:init}, we then get that

$$\int_{\rnm}|\nabla \tu_{p+1}|^2\, dx\geq
\mu_s(\rnm)^{\frac{\crit}{\crit-2}}$$
and there exists $\alpha_{p+1}\in (0,1]$ such that $\lim_{\eps\to
0}\mu_{\eps,p+1}^{\pe}=\alpha_{p+1}$. Moreover, it follows from
(\ref{hyp:lim:ye}), (\ref{hyp:lim:ye:1}) and the definition of
$\mu_{\eps,p+1}$ that

$$\lim_{\eps\to 0}\frac{\mu_{\eps,p+1}}{\mu_{\eps,p}}=+\infty\hbox{ and
}\lim_{\eps\to 0}\mu_{\eps,p+1}=0.$$
As easily checked, the families $(\mei)_{\eps>0}$, $i\in\{1,...,p+1\}$
satisfy ${\mathcal H}_{p+1}$.\hfill$\Box$

\medskip\noindent{\bf Step \ref{sec:exh}.6:} Next proposition is the
equivalent of Proposition \ref{prop:HN} at smaller scales.

\begin{prop}\label{prop:HP} Let $\Omega$ be a smooth bounded domain of
$\rn$, $n\geq 3$, such that $0\in\partial\Omega$. We let $(\ue)$, $(\ae)$
and $(\pe)$ such that $(E_\eps)$, (\ref{hyp:ae}), (\ref{lim:pe}) and (\ref{bnd:ue}) hold. Let $p\geq 1$. We assume that ${\mathcal H}_p$ holds. 
Then either for any
$i\in \{1,...,p-1\}$ and for any $\delta>0$

$$\lim_{R\to +\infty}\lim_{\eps\to 0}\sup_{x\in B_{\delta
k_{\eps,i+1}}(0)\setminus \overline{B}_{R 
k_{\eps,i}(0)}}|x|^{\frac{n-2}{2}}\left|\ue(x)-\mu_{\eps,i+1}^{-\frac{n-2}{2}}\tu_{i+1}\left(\frac{\varphi^{-1}(x)}{k_{\eps,i+1}}\right)\right|^{1-\frac{\pe}{\crit-2}}=0$$
or ${\mathcal H}_{p+1}$ holds.
\end{prop}

\noindent{\it Proof of Proposition \ref{prop:HP}:} We assume that there
exist $i\leq p-1$, $\delta>0$ such that

$$\lim_{R\to +\infty}\lim_{\eps\to 0}\sup_{x\in B_{\delta
k_{\eps,i+1}}(0)\setminus \overline{B}_{R 
k_{\eps,i}(0)}}|x|^{\frac{n-2}{2}}\left|\ue(x)-\mu_{\eps,i+1}^{-\frac{n-2}{2}}\tu_{i+1}\left(\frac{\varphi^{-1}(x)}{k_{\eps,i+1}}\right)\right|^{1-\frac{\pe}{\crit-2}}>0.$$
It then follows that there exists a family $(\ye)_{\eps>0}\in\Omega$ such
that

\beqn
&&\lim_{\eps\to 0}\frac{|\ye|}{k_{\eps,i}}=+\infty,\qquad |\ye|\leq \delta
k_{\eps,i+1}\hbox{ for all }\eps>0\label{hyp:lim:ye:1:bis}\\
&&\lim_{\eps\to
0}|\ye|^{\frac{n-2}{2}}\left|\ue(\ye)-\mu_{\eps,i+1}^{-\frac{n-2}{2}}\tu_{i+1}\left(\frac{\varphi^{-1}(\ye)}{k_{\eps,i+1}}\right)\right|^{1-\frac{\pe}{\crit-2}}=\alpha>0.\label{hyp:lim:ye:1:ter}
\eeqn
We let  $\tye\in\rnm$ such that $\ye=\varphi(k_{\eps,i+1}\tye)$. It
follows from (\ref{hyp:lim:ye:1:bis}) that $|\tye|\leq 2\delta$ for all
$\eps>0$. We claim that $\lim_{\eps\to 0}\tye=0$. Indeed, we rewrite
(\ref{hyp:lim:ye:1:ter}) as

$$\lim_{\eps\to 0}
|\tye|^{\frac{n-2}{2}}\left|\tu_{\eps,i+1}(\tye)-\tu_{i+1}(\tye)\right|^{1-\frac{\pe}{\crit-2}}=\alpha>0.$$
A contradiction with point (C3) of ${\mathcal H}_{p}$ in case $\tye\not\to
0$ when $\eps\to 0$.
\medskip\noindent Since $\tu_{i+1}\in C^0(\overline{\rnm})$, we then get
that

$$|\ye|^{\frac{n-2}{2}}\left|\mu_{\eps,i+1}^{-\frac{n-2}{2}}\tu_{i+1}\left(\frac{\varphi^{-1}(\ye)}{k_{\eps,i+1}}\right)\right|^{1-\frac{\pe}{\crit-2}}=O\left(\frac{|\ye|}{k_{\eps,i+1}}\right)^{\frac{n-2}{2}}=o(1)$$
when $\eps\to 0$. We rewrite (\ref{hyp:lim:ye:1:ter}) as

\bequa\label{hyp:lim:ye:bis}
\lim_{\eps\to
0}|\ye|^{\frac{n-2}{2}}|\ue(\ye)|^{1-\frac{\pe}{\crit-2}}=\alpha>0.
\eequa
We let

$$\ne:=|\ue(\ye)|^{-\frac{2}{n-2}}\hbox{ and 
}\elle:=\ne^{1-\frac{\pe}{\crit-2}}.$$
We define

$$\tue(x):=\ne^{\frac{n-2}{2}}\ue(\varphi(\elle x))$$
for all $x\in\frac{U}{\elle}\cap\{x_1\leq 0\}$. It follows from 
(\ref{ineq:est:1}) that
$$|\varphi(\elle x)|^{\frac{n-2}{2}}|\ue(\varphi(\elle
x))|^{1-\frac{\pe}{\crit-2}}\leq C$$
for all $x\in\frac{U}{\elle}\cap\{x_1\leq 0\}$. With the definition of 
$\tue$ and the properties (\ref{def:vphi}) of $\varphi$, we get that
there exists $C>0$ such that

$$|x|^{\frac{n-2}{2}}|\tue(x)|^{1-\frac{\pe}{\crit-2}}\leq C$$
for all $x\in\frac{U}{\elle}\cap\{x_1\leq 0\}$. It then follows that
hypothesis (L2) of Lemma \ref{lem:blowup:2} is satisfied. It then follows
from Lemma \ref{lem:blowup:2} that there exists $\tu\in \hunrnm\cap
C^{1}(\overline{\rnm})$ such that

$$\Delta \tu=\frac{|\tu|^{\crit-2}\tu}{|x|^s}\hbox{ in }{\mathcal
D}'(\rnm), $$
and
\bequa\label{cv:tup1:bis}
\lim_{\eps\to 0}\tue=\tu
\eequa
in $C^{1}_{loc}(\overline{\rnm}\setminus\{0\})$. It follows from
(\ref{hyp:lim:ye:bis}) and the definition of $\elle$ that

$$\lim_{\eps\to 0}\frac{|\ye|}{\elle}=\alpha>0.$$
We let $\bye\in \{x_1<0\}$ such that $\ye=\varphi(\elle\bye)$. It then
exists $\bar{y}_0\in\overline{\rnm}$ such that $\lim_{\eps\to
0}\bye=\bar{y}_0\neq 0$. It follows from (\ref{cv:tup1:bis}) and the
definition of $\tue$ and $\tye$ that

$$|\tu(\bar{y}_0)|=\lim_{\eps\to 0}|\tue(\bye)|=1,$$
and then $\tu\not\equiv 0$. With arguments similar to the ones developed
in the proof of Lemma \ref{lem:init}, we then get that

$$\int_{\rnm}|\nabla \tu|^2\, dx\geq \mu_s(\rnm)^{\frac{\crit}{\crit-2}}$$
and there exists $\alpha\in (0,1]$ such that $\lim_{\eps\to
0}\ne^{\pe}=\alpha$. Moreover, it follows from (\ref{hyp:lim:ye:bis}),
(\ref{hyp:lim:ye:1:bis}) and the definition of $\ne$ that

$$\lim_{\eps\to 0}\frac{\ne}{\mu_{\eps,i}}=+\infty\hbox{ and
}\lim_{\eps\to 0}\frac{\mu_{\eps,i+1}}{\ne}=+\infty.$$
As easily checked, the families $(\meun)$,..., $(\mei)$, $(\ne)$,
$(\mu_{\eps,i+1})$,..., $(\meN)_{\eps>0}$ satisfy ${\mathcal
H}_{p+1}$.\hfill$\Box$

\medskip\noindent{\bf Step \ref{sec:exh}.7:} This last Step is the proof
of Proposition \ref{prop:exhaust}.

\begin{prop}\label{prop:sec.7} Let $\Omega$ be a smooth bounded domain of
$\rn$, $n\geq 3$, such that $0\in\partial\Omega$. We let $(\ue)$, $(\ae)$
and $(\pe)$ such that $(E_\eps)$, (\ref{hyp:ae}), (\ref{lim:pe}) and 
(\ref{bnd:ue}) hold. We let $N_0=\max\{p/\, {\mathcal H}_p\hbox{
holds}\}$. Then $N_0<+\infty$ and the conclusion of Proposition
\ref{prop:exhaust} holds with $N=N_0$.
\end{prop}

\noindent{\it Proof of Proposition \ref{prop:sec.7}:} Indeed, assume that
${\mathcal H}_p$ holds. Let $\delta,R>0$. Since $\mei=o(\mu_{\eps,i+1})$
for all $i\in\{1,...,N-1\}$, we then get with a change of variable and the
definition of $\tuei$ (see (C3)) that

\beq
\int_{\Omega}|\nabla\ue|^2\, dx&\geq &\sum_{i=1}^N\int_{\varphi(B_{R
\kei}(0)\setminus \overline{B}_{\delta\kei}(0))}|\nabla\ue|^2\, dx\\
&\geq & \sum_{i=1}^N\mei^{-\frac{n-2}{\crit-2}\pe}\int_{B_{R}(0)\setminus
\overline{B}_{\delta}(0)}|\nabla\tuei|_{g_{\eps,i}}^2\, dv_{g_{\eps,i}}\\
&\geq & \sum_{i=1}^N\int_{B_{R}(0)\setminus
\overline{B}_{\delta}(0)}|\nabla\tuei|_{g_{\eps,i}}^2\, dv_{g_{\eps,i}}
\eeq
where $g_{\eps,i}$ is the metric such that
$(g_{\eps,i})_{qr}=(\partial_q\varphi(\kei x),\partial_r\varphi(\kei x))$
for all $q,r\in\{1,...,p\}$. Passing to the limit $\eps\to 0$ and using
point (C3) of ${\mathcal H}_p$, we get that

$$\int_{\Omega}|\nabla\ue|^2\, dx\geq p
\mu_s(\rnm)^{\frac{\crit}{\crit-2}}+o(1)$$
when $\eps\to 0$. With (\ref{bnd:ue}), we get that there exists $C>0$ such
that $$p\leq \Lambda^2 \mu_s(\rnm)^{-\frac{\crit}{\crit-2}}.$$ It then
follows that $N_0<+\infty$ exists.

\medskip\noindent We let families $(\meun)_{\eps>0}$,...,
$(\mu_{\eps,N_0})_{\eps>0}$ such that ${\mathcal H}_{N_0}$ holds. We argue
by contradiction and assume that the conclusion of Proposition
\ref{prop:exhaust} does not hold with $N=N_0$. Assertions (A1), (A2), (A3)
(A4) and (A7) hold. Assume that (A5) or (A6) does not hold. It then
follows from Propositions \ref{prop:HN} and \ref{prop:HP} that ${\mathcal
H}_{N+1}$ holds. A contradiction with the choice of $N=N_0$, and the
proposition is proved.\hfill$\Box$

\section{Strong pointwise estimates, Part 1}\label{sec:fe:1}
The objective of this section is the proof of the following strong 
pointwise estimate:

\begin{prop}\label{prop:fund:est:part1} Let $\Omega$ be a smooth bounded
domain of $\rn$, $n\geq 3$. We let $s\in (0,2)$. We let $(\pe)_{\eps>0}$ 
such that $\pe\in [0,\crit-2)$ for all $\eps>0$ and (\ref{lim:pe}) holds. 
We consider $(\ue)_{\eps>0}\in \huno$ such that (\ref{hyp:ae}), $(E_\eps)$ 
and (\ref{bnd:ue}) hold. We assume that blow-up occurs, that is
$$\lim_{\eps\to 0}\Vert\ue\Vert_{L^\infty(\Omega)}=+\infty.$$
We let $\mu_{\eps,1},...,\mu_{\eps,N}$ as in Proposition 
\ref{prop:exhaust}. Then, there exists $C>0$ such that

\bequa\label{eq:fund:est:part1}
|\ue(x)|\leq
C\sum_{i=1}^N\frac{\mei^{\frac{n}{2}}|x|}{\left(\mei^2+|x|^2\right)^{\frac{n}{2}}}+C |x|
\eequa
for all $\eps>0$ and all $x\in\Omega$.
\end{prop}

\noindent The proof of this estimate goes through seven steps. We let 
$s\in (0,2)$. We let $(\pe)_{\eps>0}$ such that $\pe\in [0,\crit-2)$ for 
all $\eps>0$ and (\ref{lim:pe}) holds. We consider $(\ue)_{\eps>0}\in 
\huno$ that satisfies the hypothesis of Proposition 
\ref{prop:fund:est:part1}. We let $\mu_{\eps,1},...,\mu_{\eps,N}$ as in 
Proposition \ref{prop:exhaust}.

\medskip\noindent{\bf Step \ref{sec:fe:1}.1:} We claim that for any 
$\nu\in (0,1)$ and any $R>0$, there exists $C(\nu,R)>0$ such that
\bequa\label{estim:nu:1}
|\ue(x)|\leq 
C(\nu,R)\cdot\left(\frac{\meN^{\frac{n}{2}-\nu(n-1)}d(x,\partial\Omega)^{1-\nu}}{|x|^{n(1-\nu)}}+d(x,\partial\Omega)^{1-\nu}\right)
\eequa
for all $x\in\Omega\setminus \overline{B}_{R k_{\eps,N}}(0)$ and all
$\eps>0$.\par

\smallskip\noindent{\it Proof of the Claim:} Since $\Delta$ is coercive
on $\Omega$, we let $G$ be the Green's function for $\Delta$ in $\Omega$ 
with Dirichlet boundary condition. We let

$$H(x)=-\partial_\nu G(x,0)$$
for all $x\in\overline{\Omega}\setminus\{0\}$. Here $\nu$ denotes the 
outward normal vector at $\partial\Omega$. It follows from Theorem
\ref{th:green:2} of the Appendix that $H\in
C^2(\overline{\Omega}\setminus\{0\})$, that
\bequa\label{eq:Ge}
\Delta H=0
\eequa
in $\Omega$ and that there exist $\delta_1,C_1>0$ such that

\bequa\label{ineq:Ge:1}
\frac{d(x,\partial\Omega)}{C_1|x|^{n}}\leq H(x)\leq \frac{C_1
d(x,\partial\Omega)}{|x|^{n}}
\eequa
and
\bequa\label{ineq:Ge:2}
\frac{|\nabla H(x)|}{H(x)}\geq \frac{1}{C_1d(x,\partial\Omega)}\geq 
\frac{1}{C_1 |x|}
\eequa
for all $x\in \Omega\cap B_{2\delta_1}(0)$.

\medskip\noindent Since $\Delta$ is coercive, we let $\lambda_{1}>0$ be 
the first eigenvalue of $\Delta$ on $\Omega$, and we let
$\psi\in C^2(\overline{\Omega})$ be the unique eigenfunction such
that

$$\left\{\begin{array}{ll}
\Delta\psi=\lambda_{1}\psi & \hbox{ in
}\Omega\\
\psi>0 & \hbox{ in }\Omega\\
\psi=0 & \hbox{ on }\partial\Omega\\
\int_\Omega\psi^2\, dx=1
\end{array}\right\}$$
It follows from standard elliptic theory and Hopf's maximum principle that 
there exists $C'_2,\delta_2>0$ such that

$$\frac{1}{C'_2} d(x,\partial\Omega)\leq \psi(x)\leq C'_2
d(x,\partial\Omega)\hbox{ and }\frac{1}{C'_2}\leq |\nabla\psi(x)|\leq
C'_2
$$
for all $x\in \Omega\cap B_{2\delta_2}(0)$. Consequently, there exists 
$C_2>0$ such that
\bequa\label{ppty:phi}
\frac{1}{C_2}d(x,\partial\Omega)\leq \psi(x)\leq C_2
d(x,\partial\Omega)\hbox{ and
}\frac{|\nabla\psi(x)|}{\psi(x)}\geq\frac{1}{C_2
d(x,\partial\Omega)}\geq\frac{1}{C_2 |x|}
\eequa
for all $x\in \Omega\cap B_{2\delta_2}(0)$. We let the operator

$$L_\eps=\Delta +\left(\ae-\frac{|\ue|^{\crit-2-\pe}}{|x|^s}\right).$$

\medskip\noindent{\it Step \ref{sec:fe:1}.1.1:} We claim that there exist
$\delta_0>0$ and $R_0>0$ such that for any $\nu\in (0,1)$ and any $R>R_0$, 
$\delta\in (0,\delta_0)$, we have that

\bequa\label{ineq:LG}
L_\eps H^{1-\nu}>0,\hbox{ and }L_\eps\psi^{1-\nu}>0
\eequa
for all $x\in \Omega\cap B_\delta(0)\setminus \overline{B}_{R\keN}(0)$ and
for all $\eps>0$ sufficiently small. Indeed, with (\ref{eq:Ge}), we get 
that

\bequa\label{eq:Ge:2}
\frac{L_\eps H^{1-\nu}}{H^{1-\nu}}(x)=\ae(x)+\nu(1-\nu)\frac{|\nabla
H|^2}{H^2}(x)-\frac{|\ue(x)|^{\crit-2-\pe}}{|x|^s}
\eequa
for all $x\in \Omega\setminus\{0\}$ and all $\eps>0$. We let
$0<\delta_0\leq\min\{\delta_1,\delta_2\}$ such that

\bequa\label{def:delta0}
\left\{\begin{array}{l}
2\delta_0^2\sup_{\Omega}|\ae|\leq \frac{\nu(1-\nu)}{2\cdot
\max\{C_1^2,C_2^2\}}\\
\\
2^{\crit+1}\delta_0^{2-s}\Vert
u_0\Vert_{L^\infty(\Omega)}^{\crit-2}<\frac{\nu(1-\nu)}{4\cdot
\max\{C_1^2,C_2^2\}}\end{array}\right\}
\eequa
for all $\eps>0$. This choice is possible thanks to (\ref{hyp:ae}). It
follows from point (A5) of Proposition \ref{prop:exhaust} that there
exists $R_0>0$ such that for any $R>R_0$, we have that

$$|x|^{\frac{n-2}{2}}|\ue(x)-u_0(x)|^{1-\frac{\pe}{\crit-2}}\leq
\left(\frac{\nu(1-\nu)}{2^{\crit+1}\max\{C_1^2,C_2^2\}}\right)^{\frac{1}{\crit-2}}$$
for all $x\in \Omega\setminus \overline{B}_{R\keN}(0)$ and all $\eps>0$.
We then get that

\beq
|x|^{2-s}|\ue(x)|^{\crit-2-\pe}&\leq&
2^{\crit-1-\pe}|x|^{2-s}|\ue(x)-u_0(x)|^{\crit-2-\pe}\\
&&+2^{\crit-1-\pe}|x|^{2-s}
|u_0(x)|^{\crit-2-\pe}\\
&\leq & 2^{-\pe}\frac{\nu(1-\nu)}{4\cdot
\max\{C_1^2,C_2^2\}}\\
&&+2^{\crit-1-\pe}\delta^{2-s}\Vert 
u_0\Vert_{L^\infty(\Omega)}^{\crit-2-\pe}
\eeq
for all $x\in \Omega\setminus \overline{B}_{R\keN}(0)$ and all $\eps>0$. 
We get with the
choice (\ref{def:delta0}) of $\delta_0$ that for any $\delta\in
(0,\delta_0)$ and all $R>R_0$

\beq
|x|^{2-s}|\ue(x)|^{\crit-2-\pe}&\leq& \frac{\nu(1-\nu)}{4\cdot
\max\{C_1^2,C_2^2\}}+2^{\crit-1}\delta_0^{2-s} \Vert
u_0\Vert_{L^\infty(\Omega)}^{\crit-2}\\
&<& \frac{\nu(1-\nu)}{2\cdot \max\{C_1^2,C_2^2\}}
\eeq
for all $x\in (B_\delta(0)\setminus \overline{B}_{R\keN}(0))\cap\Omega$
and all $\eps>0$ small enough. With (\ref{eq:Ge:2}) and
(\ref{def:delta0}), we get that

\beq
\frac{L_\eps H^{1-\nu}}{H^{1-\nu}}(x)&\geq& \frac{\nu(1-\nu)}{
C_1^2|x|^2}+\ae(x)-\frac{\nu(1-\nu)}{2C_1^2|x|^2}\\
&\geq& \frac{\nu(1-\nu)-2C_1^2|x|^2|\ae(x)|}{2C_1^2|x|^2}>0
\eeq
for all $x\in (B_\delta(0)\setminus \overline{B}_{R\keN}(0))\cap\Omega$
and all $\eps>0$ small enough. We deal with the second inequality of
(\ref{ineq:LG}). We have that

$$\frac{L_\eps
\psi^{1-\nu}}{\psi^{1-\nu}}(x)=\ae(x)+(1-\nu)\lambda_{1}+\nu(1-\nu)\frac{|\nabla
\psi|^2}{\psi^2}(x)-\frac{|\ue(x)|^{\crit-2-\pe}}{|x|^s}$$
for all $x\in\Omega$. With (\ref{ppty:phi}) and (\ref{def:delta0}) we get
that
$$\frac{L_\eps \psi^{1-\nu}}{\psi^{1-\nu}}(x)\geq
\frac{\nu(1-\nu)-2 C_2^2|\ae(x)|\delta^2+2(1-\nu)\lambda_{1}|x|^2
C_2^2}{2C_2^2|x|^2}>0$$
for all $x\in (B_\delta(0)\setminus \overline{B}_{R\keN}(0))\cap\Omega$
and all $\eps>0$. This proves the last inequality of (\ref{ineq:LG}).

\medskip\noindent{\it Step \ref{sec:fe:1}.1.2:} It follows from point (A4)
of Proposition \ref{prop:exhaust} that there exists $C_1(R)>0$ such that

\bequa\label{sup:ue:boundary:2}
|\ue(x)|\leq C_1(R) \meN^{-\frac{n}{2}}d(x,\partial\Omega)
\eequa
for all $x\in\Omega\cap\partial B_{R\keN}(0)$ and all $\eps>0$. It follows
from point (A1) of Proposition \ref{prop:exhaust} that there exists
$C_2(\delta)>0$ such that

\bequa\label{sup:ue:boundary:1}
|\ue(x)|\leq C_2(\delta)d(x,\partial\Omega)
\eequa
for all $x\in\Omega\cap\partial B_{\delta}(0)$ and all $\eps>0$. We let
$$D_{\eps,R,\delta}:=(B_\delta(0)\setminus
\overline{B}_{R\keN}(0))\cap\Omega.$$
We let
$$\alpha_\eps:=2C_1(R)C_1^{1-\nu}R^{n-(n-1)\nu}\alpha_N^{-\frac{n-(n-1)\nu}{\crit-2}}\mu_{\eps,N}^{\frac{n}{2}-\nu(n-1)}$$
and
$$\beta_\eps:=2\delta^\nu C_2(\delta)C_2^{1-\nu},$$
and
$$\varphi_\eps(x)=\alpha_\eps H^{1-\nu}(x)+\beta_\eps \psi^{1-\nu}(x)$$
for all $x\in \overline{D_{\eps,R,\delta}}$ and all $\eps>0$. Here,
$\alpha_N$ is as in point (A7) of Proposition \ref{prop:exhaust}. We claim
that
\bequa\label{control:ue:bord}
|\ue(x)|\leq \varphi_\eps(x)
\eequa
for all $\eps>0$ and all $x\in \partial D_{\eps,R,\delta}$. Indeed, with
inequalities (\ref{ineq:Ge:1}) and (\ref{sup:ue:boundary:2}), we get that
for any $x\in \Omega\cap\partial B_{R\keN}(0)$,

\beq
&&\frac{|\ue(x)|}{\alpha_\eps H(x)^{1-\nu}}\leq
\frac{\meN^{\nu(n-1)-n}d(x,\partial\Omega)^\nu|x|^{n-n\nu}}{2 
R^{n-(n-1)\nu}\alpha_N^{-\frac{n-(n-1)\nu}{\crit-2}}}\leq
\frac{\meN^{-\frac{n-(n-1)\nu}{\crit-2}\pe}}{2\alpha_N^{-\frac{n-(n-1)\nu}{\crit-2}}}\leq 
1
\eeq
when $\eps\to 0$ with point (A7) of Proposition \ref{prop:exhaust}.
Similarly, we have with (\ref{ppty:phi}) and (\ref{sup:ue:boundary:1})
that

$$\frac{|\ue(x)|}{\beta_\eps \psi(x)^{1-\nu}}\leq
\frac{d(x,\partial\Omega)^\nu}{2\delta^\nu}<1$$
for all $x\in \Omega\cap\partial B_{\delta}(0)$ and all $\eps>0$. On
$\partial\Omega\cap (B_\delta(0)\setminus \overline{B}_{R\keN}(0))$, we
clearly have $\varphi_\eps(x)>|\ue(x)|=0$. As easily checked, these
assertions prove (\ref{control:ue:bord}).

\medskip\noindent{\it Step \ref{sec:fe:1}.1.3:} We claim that $L_\eps$
verifies the following comparison maximum: if $\varphi\in 
C^2(D_{\eps,R,\delta})\cap C^0(\overline{D_{\eps,R,\delta}})$, then
$$\left\{\begin{array}{ll}
L_\eps\varphi\geq 0 & \hbox{ in }D_{\eps,R,\delta}\\
\varphi\geq 0& \hbox{ on }\partial D_{\eps,R,\delta}
\end{array}\right\}\Rightarrow\;\varphi\geq 0 \hbox{ in
}D_{\eps,R,\delta}.$$
Indeed, we let $U_0$ be an open subset of $\rn$ such that 
$\overline{\Omega}\subset\subset U_0$. Since the operator $\Delta$ is 
coercive in $U_0$ (with boundary Dirichlet condition), we let 
$\tilde{G}\in C^2(U_0\times U_0\setminus\{(x,x)/x\in U_0\})$ be the 
Green's function for $\Delta$ with Dirichlet condition in $U_0$. In other 
words, $\tilde{G}$ satisfies $$\Delta \tilde{G}(x,\cdot)=\delta_x$$
weakly in ${\mathcal D}(U_0)$. For the existence, we refer to Theorem 
\ref{th:green:1} of Appendix B. Moreover, since $0\in U_0$ is in the 
interior of the domain, there exists $\hat{\delta}_0>0$ and $C_0>0$ such 
that $$\frac{1}{C_0|x|^{n-2}}\leq \tilde{G}(0,x)\leq \frac{C_0}{|x|^{n-2}}$$
and
$$\frac{|\nabla \tilde{G}(0,x)|}{\tilde{G}(0,x)}\geq
\frac{C_0'}{|x|}$$
for all $\eps>0$ and all $x\in B_{2\hat{\delta}_0}(0)\setminus\{0\}$. The
proof of these estimates goes as in the proof of points (G9) and (G10) of
Theorem \ref{th:green:2}. We refer to \cite{dhr} for the details. With the
same techniques as in Step \ref{sec:fe:1}.1.1, we get that for $R>0$ large 
enough and $\delta>0$ small enough, then

$$\tilde{G}^{1-\nu}>0\hbox{ and }L_\eps\tilde{G}^{1-\nu}>0$$
in
$\overline{D_{\eps,R,\delta}}$ for all $\eps>0$. It then follows from 
\cite{bnv} that $L_\eps$ verifies the above mentioned comparison 
principle.

\medskip\noindent{\it Step \ref{sec:fe:1}.1.4:} It follows from
(\ref{ineq:LG}) and (\ref{control:ue:bord}) that

$$\left\{\begin{array}{ll}
L_\eps\varphi_\eps\geq 0=L_\eps\ue & \hbox{ in }D_{\eps,R,\delta}\\
\varphi_\eps\geq 0=\ue& \hbox{ on }\partial D_{\eps,R,\delta}\\
L_\eps\varphi_\eps\geq 0=-L_\eps\ue & \hbox{ in }D_{\eps,R,\delta}\\
\varphi_\eps\geq 0=-\ue& \hbox{ on }\partial D_{\eps,R,\delta}
\end{array}\right\}.$$
It follows from the above comparison principle that

$$|\ue(x)|\leq \varphi_\eps(x)$$
for all $x\in D_{\eps,R,\delta}$. With (\ref{ineq:Ge:1}), we then get that
(\ref{estim:nu:1}) holds on $D_{\eps,R,\delta}=(B_\delta(0)\setminus 
\overline{B}_{R k_{\eps,N}}(0))\cap\Omega$ for $R$ large and $\delta$ 
small. It follows
from this last assertion, (\ref{ineq:Ge:1}) and points (A1) and (A4) of
Proposition \ref{prop:exhaust} that (\ref{estim:nu:1}) holds on 
$\Omega\setminus \overline{B}_{R k_{\eps,N}}(0)$ for all 
$R>0$.\hfill$\Box$

\medskip\noindent{\bf Step \ref{sec:fe:1}.2:} Let $i\in \{1,...,N-1\}$. We
claim that for any $\nu\in (0,1)$ and any $R,\rho>0$, there exists 
$C(\nu,R,\rho)>0$ such that

\bequa\label{estim:nu:1:bis}
|\ue(x)|\leq C(\nu,R,\rho)\cdot\left(
\frac{\mu_{\eps,i}^{\frac{n}{2}-\nu(n-1)}d(x,\partial\Omega)^{1-\nu}}{|x|^{n(1-\nu)}}+
\mu_{i+1,\eps}^{\nu-\frac{n}{2}} d(x,\partial\Omega)^{1-\nu}\right)
\eequa
for all $x\in B_{R k_{\eps,i+1}}(0)\setminus \overline{B}_{\rho
k_{\eps,i}}(0)$ and all $\eps>0$.\par

\smallskip\noindent{\it Proof of the Claim:} We let $i\in \{1,...,N-1\}$.
We follow the lines of the proof of Step \ref{sec:fe:1}.1. We let $H$
and $\psi$ as in Step \ref{sec:fe:1}.1. Recall that we then get that
there exists $\delta_1>0$ and $C_1>0$ such that

\bequa\label{ineq:Ge:1:43}
\frac{d(x,\partial\Omega)}{C_1|x|^{n}}\leq H(x)\leq \frac{C_1
d(x,\partial\Omega)}{|x|^{n}}
\eequa
and
\bequa\label{ineq:Ge:2:43}
\frac{|\nabla H(x)|}{H(x)}\geq \frac{1}{C_1d(x,\partial\Omega)}\geq
\frac{1}{C_1|x|}
\eequa
for all $x\in B_{2\delta_1}(0)\setminus\{0\}$. Moreover
there exists $C_2,\delta_2>0$ such that $\psi$ verifies

\bequa\label{ppty:phi:43}
\frac{1}{C_2}d(x,\partial\Omega)\leq \psi(x)\leq C_2
d(x,\partial\Omega)\hbox{ and
}\frac{|\nabla\psi(x)|}{\psi(x)}\geq\frac{1}{C_2|x|}
\eequa
for all $x\in \Omega\cap B_{2\delta_2}(0)$. We let
the operator

$$L_\eps=\Delta +\left(\ae-\frac{|\ue|^{\crit-2-\pe}}{|x|^s}\right).$$

\medskip\noindent{\it Step \ref{sec:fe:1}.2.1:} We claim that there exist
$\rho_0>0$ and $R_0>0$ such that for any $\nu\in (0,1)$ and any $R>R_0$, 
$\rho\in (0,\rho_0)$, we have that

\bequa\label{ineq:LG:43}
L_\eps H^{1-\nu}>0\hbox{ and }L_\eps \psi^{1-\nu}>0
\eequa
for all $x\in \Omega\cap(B_{\rho k_{\eps,i+1}}(0)\setminus \overline{B}_{R
k_{\eps,i}}(0))$ and for all $\eps>0$ sufficiently small. Indeed, as in 
Step \ref{sec:fe:1}.1, we get that

\bequa\label{eq:Ge:2:43}
\frac{L_\eps H^{1-\nu}}{H^{1-\nu}}(x)=\ae(x)+\nu(1-\nu)\frac{|\nabla
  H|^2}{H^2}(x)-\frac{|\ue(x)|^{\crit-2-\pe}}{|x|^s}
\eequa
for all $x\in \Omega\setminus\{0\}$ and all $\eps>0$. We let $0<\rho_0<1$
such that

\bequa\label{def:delta0:43}
2^{\crit+1}\rho_0^{2-s}\Vert
\tu_{i+1}\Vert_{L^\infty(B_2(0)\cap\rnm)}^{\crit-2}<\frac{\nu(1-\nu)}{\max\{C_1^2,C_2^2\}}
\eequa
for all $\eps>0$. It follows from point (A6) of Proposition
\ref{prop:exhaust} that there exists $R_0>0$ such that for any $R>R_0$

$$|x|^{\frac{n-2}{2}}\left|\ue(x)-\mu_{\eps,i+1}^{\frac{n-2}{2}}\tu_{i+1}\left(\frac{\varphi^{-1}(x)}{k_{\eps,i+1}}\right)\right|^{1-\frac{\pe}{\crit-2}}\leq 
\left(\frac{\nu(1-\nu)}{2^{\crit+1}\max\{C_1^2,C_2^2\}}\right)^{\frac{1}{\crit-2}}$$
for all $x\in \Omega\cap\left(B_{k_{\eps,i+1}}(0)\setminus
\overline{B}_{R\kei}(0)\right)$ and all $\eps>0$. We then get that 
\beq
&&|x|^{2-s}|\ue(x)|^{\crit-2-\pe}\\
&&\leq
2^{\crit-1-\pe}|x|^{2-s}\left|\ue(x)-\mu_{\eps,i+1}^{\frac{n-2}{2}}\tu_{i+1}\left(\frac{\varphi^{-1}(x)}{k_{\eps,i+1}}\right)\right|^{\crit-2-\pe}\\
&&+2^{\crit-1-\pe}|x|^{2-s}\mu_{\eps,i+1}^{-\frac{n-2}{2}(\crit-2)\cdot(1-\frac{\pe}{\crit-2})}
\sup_{B_2(0)\cap\rnm}|\tu_{i+1}|^{\crit-2-\pe}\\
&&\leq 
2^{-\pe}\frac{\nu(1-\nu)}{4\cdot\max\{C_1^2,C_2^2\}}+2^{\crit-1-\pe}\left(\frac{|x|}{k_{\eps,i+1}}\right)^{2-s}\Vert 
\tu_{i+1}\Vert_{L^\infty(B_2(0)\cap\rnm)}^{\crit-2-\pe}
\eeq
for all $x\in \Omega\cap \left(B_{k_{\eps,i+1}}(0)\setminus
\overline{B}_{R\kei}(0)\right)$ and all $\eps>0$. We then get with the choice
(\ref{def:delta0:43}) of $\rho_0$ that for any $\rho\in (0,\rho_0)$ and
all $R>R_0$

\beq
|x|^{2-s}|\ue(x)|^{\crit-2-\pe}&\leq&
\frac{\nu(1-\nu)}{4\cdot\max\{C_1^2,C_2^2\}}+2^{\crit-1}\rho_0^{2-s}\Vert\tu_{i+1}\Vert_{L^\infty(B_2(0)\cap\rnm)}^{\crit-2}+o(1)\\
&<& \frac{\nu(1-\nu)}{2\cdot\max\{C_1^2,C_2^2\}}
\eeq
for all $x\in (B_{\rho k_{\eps,i+1}}(0)\setminus
\overline{B}_{R\kei}(0))\cap\Omega$ and all $\eps>0$ small enough. Since 
(\ref{hyp:ae}) holds, we get with (\ref{eq:Ge:2:43}) that

\beq
\frac{L_\eps H^{1-\nu}}{H^{1-\nu}}(x)&\geq&
\frac{\nu(1-\nu)}{C_1^2|x|^2}+\ae(x)-\frac{\nu(1-\nu)}{2C_1^2|x|^2}\\
&\geq& \frac{\nu(1-\nu)+2|x|^2C_1^2\ae(x)}{2C_1^2|x|^2}>0
\eeq
for all  $x\in \Omega\cap \left(B_{\rho k_{\eps,i+1}}(0)\setminus
\overline{B}_{R\kei}(0)\right)$ and all $\eps>0$ small enough. The proof
of the second inequality of (\ref{ineq:LG:43}) goes similarly (see Step
\ref{sec:fe:1}.1 for details). This proves (\ref{ineq:LG:43}).

\medskip\noindent{\it Step \ref{sec:fe:1}.2.2:} It follows from point (A4)
of Proposition \ref{prop:exhaust} that there exists $C_1(R)>0$ and
$C_2(\rho)$ such that

$$\begin{array}{ll}
|\ue(x)|\leq C_1(R) \mu_{\eps,i}^{-\frac{n}{2}}d(x,\partial\Omega) &
\hbox{ for all } x\in\Omega\cap\partial B_{R\kei}(0)\\
&\\
|\ue(x)|\leq C_2(\rho) \mu_{\eps,i+1}^{-\frac{n}{2}}d(x,\partial\Omega) &
\hbox{ for all } x\in\Omega\cap\partial B_{\rho k_{\eps,i+1}}(0).
\end{array}$$
We let
$$D_{\eps,R,\rho}:=\left(B_{\rho k_{\eps,i+1}}(0)\setminus
\overline{B}_{R\kei}(0)\right)\cap\Omega.$$
We let
$$\alpha_\eps:=2C_1(R)C_1^{1-\nu}R^{n-\nu(n-1)}\alpha_i^{-\frac{n-\nu(n-1)}{\crit-2}}\mu_{\eps,i}^{\frac{n}{2}-\nu(n-1)}$$
and
$$\beta_\eps:=2C_2(\rho)C_2^{1-\nu}\rho^{\nu}\alpha_{i+1}^{-\frac{\nu}{\crit-2}}\mu_{\eps,i+1}^{-\frac{n}{2}+\nu},$$
and
$$\varphi_\eps(x):=\alpha_\eps H^{1-\nu}(x)+\beta_\eps
\psi^{1-\nu}(x)$$
for all $x\in \overline{D_{\eps,R,\delta}}$ and all $\eps>0$. Here, the
$\alpha_i$'s are as in Point (A7) of Proposition \ref{prop:exhaust}.
Similarly to what was done in Step \ref{sec:fe:1}.1, we then get that
\bequa\label{control:ue:bord:43}
|\ue(x)|\leq \varphi_\eps(x)
\eequa
for all $\eps>0$ and all $x\in \partial D_{\eps,R,\rho}$. The operator $L_\eps$ verifies the
comparison principle on $D_{\eps,R,\rho}$ as in Step \ref{sec:fe:1}.1.3.
It then follows that

$$|\ue(x)|\leq \varphi_\eps(x)$$
for all $x\in D_{\eps,R,\rho}$. With (\ref{ineq:Ge:1:43}), we then get
that (\ref{estim:nu:1:bis}) holds on $D_{\eps,R,\rho}$ for $R$ large and 
$\rho$ small. It follows from this last assertion and point
(A4) of Proposition \ref{prop:exhaust} that (\ref{estim:nu:1:bis}) holds 
on $\left(B_{\rho k_{\eps,i+1}}(0)\setminus
\overline{B}_{R\kei}(0)\right)\cap\Omega$ for all $R,\rho>0$.\hfill$\Box$

\medskip\noindent{\bf Step \ref{sec:fe:1}.3:} As easily checked, it
follows from (\ref{estim:nu:1}), (\ref{estim:nu:1:bis}) and Proposition
\ref{prop:exhaust} that for any $\nu\in (0,1)$, there exists $C_\nu>0$
such that
\bequa\label{eq:est:nu:global}
|\ue(x)|\leq
C_\nu\sum_{i=1}^N\frac{\mei^{\frac{n}{2}-(n-1)\nu}|x|^{1-\nu}}{\left(\mei^2+|x|^2\right)^{\frac{n}{2}(1-\nu)}}+C_\nu 
|x|^{1-\nu}
\eequa
for all $x\in\Omega$ and all $\eps>0$. Note that we have used that
$d(x,\partial\Omega)\leq |x-0|=|x|$ for all $x\in\Omega$. We let $G$ be
the Green's function of $\Delta$ on $\Omega$ with Dirichlet boundary 
condition. It follows from Green's
representation formula and (\ref{eq:est:nu:global}) that

\beqn
|\ue(x)|&=&\left|\int_\Omega
G(x,y)\left(\frac{|\ue(y)|^{\crit-2-\pe}\ue(y)}{|y|^s}-\ae(y)\ue(y)\right)\, 
dy\right|\label{eq:green:ue}\\
&\leq & C\int_\Omega
G(x,y)\left(\frac{|\ue(y)|^{\crit-1-\pe}}{|y|^s}+1\right)\, dy\nonumber\\
&\leq & C_\nu\sum_{i=1}^N\int_\Omega
\frac{G(x,y)}{|y|^s}\left(\frac{\mei^{\frac{n}{2}-(n-1)\nu}|y|^{1-\nu}}{\left(\mei^2+|y|^2\right)^{\frac{n}{2}(1-\nu)}}\right)^{\crit-1-\pe}\,
dy\nonumber\\
&&+C_\nu\int_\Omega G(x,y)\left(|y|^{(1-\nu)(\crit-1-\pe)-s}+1\right)\,
dy\label{ineq:ue:ef2}
\eeqn

\medskip\noindent{\bf Step \ref{sec:fe:1}.4:} We claim that there exists
$C>0$ such that
\bequa\label{ineq:He:alpha}
\int_\Omega G(x,y)\left(|y|^{(1-\nu)(\crit-1-\pe)-s}+1\right)\,
dy\leq C |x|
\eequa
\noindent{\it Proof of the Claim:} Indeed, we let $\psi_\eps\in
H_{1,0}^p(\Omega)$ ($1<p<\frac{n}{s}$) such that

$$\Delta\psi_\eps=|y|^{(1-\nu)(\crit-1-\pe)-s}+1\hbox{ in
}{\mathcal D}'(\Omega).$$
Here, $H_{1,0}^p(\Omega)$ denote the completion of $C^\infty_c(\Omega)$
for the norm $\Vert\cdot\Vert:=\Vert\nabla\cdot\Vert_p$. Since $s\in (0,2)$, it follows from standard elliptic theory 
that for $\nu>0$ small,
$\psi_\eps\in C^1(\overline{\Omega})$ and that there exists $C>0$ such
that
$$\Vert\psi_\eps\Vert_{C^{1}(\overline{\Omega})}\leq C.$$
Since $\psi_\eps(0)=0$, we get that
$$|\psi_\eps(x)|\leq C|x|$$
for all $x\in\Omega$. Moreover, since $s\in (0,2)$, we get with Green's representation formula that
$$\psi_\eps(x)=\int_\Omega G(x,y)\left(|y|^{(1-\nu)(\crit-1-\pe)-s}+1\right)\, dy$$
for all $x\in\Omega$ and all $\eps>0$. Inequation (\ref{ineq:He:alpha}) then
follows.\hfill$\Box$

\medskip\noindent{\bf Step \ref{sec:fe:1}.5:} We let $i\in \{1,...,N\}$.
We claim that there exists $C>0$ such that

\beqn\label{ineq:bubble:1}
&&\int_\Omega
\frac{G(x,y)}{|y|^s}\left(\frac{\mei^{\frac{n}{2}-(n-1)\nu}|y|^{1-\nu}}{\left(\mei^2+|y|^2\right)^{\frac{n}{2}(1-\nu)}}\right)^{\crit-1-\pe}\,
dy\nonumber\\
&&\leq
C\frac{\mei^{\frac{n}{2}}|x|}{\left(\mei^2+|x|^2\right)^{\frac{n}{2}}}
\eeqn
for all $x\in\Omega$ such that $|x|\geq \mei$.

\smallskip\noindent{\it Proof of the Claim:} Indeed, with point (G6) of
Theorem \ref{th:green:1} on the Green's function, we get that
\beq
&&\int_\Omega
\frac{G(x,y)}{|y|^s}\left(\frac{\mei^{\frac{n}{2}-(n-1)\nu}|y|^{1-\nu}}{\left(\mei^2+|y|^2\right)^{\frac{n}{2}(1-\nu)}}\right)^{\crit-1-\pe}\,
dy\\
&&\leq C\int_\Omega
\frac{|y|}{|x-y|^{n-1}|y|^s}\left(\frac{\mei^{\frac{n}{2}-(n-1)\nu}|y|^{1-\nu}}{\left(\mei^2+|y|^2\right)^{\frac{n}{2}(1-\nu)}}\right)^{\crit-1-\pe}\,
dy\\
&&\leq I_{1,\eps}(x)+I_{2,\eps}(x).
\eeq
Here,
$$I_{i,\eps}(x):=C\int_{\Omega_i(x)}
\frac{|y|}{|x-y|^{n-1}|y|^s}\left(\frac{\mei^{\frac{n}{2}-(n-1)\nu}|y|^{1-\nu}}{\left(\mei^2+|y|^2\right)^{\frac{n}{2}(1-\nu)}}\right)^{\crit-1-\pe}\,
dy$$
where
$$\Omega_1(x)=\Omega\cap\{|x-y|>|x|/2\}\hbox{ and
}\Omega_2(x)=\Omega\cap\{|x-y|<|x|/2\}.$$
We compute these two integrals separately. We let $R>0$ such that
$\Omega\subset B_R(0)$. We have that

\beqn
I_{1,\eps}(x)&\leq& C|x|^{1-n}\int_{B_R(0)}
|y|^{1-s}\left(\frac{\mei^{\frac{n}{2}-(n-1)\nu}|y|^{1-\nu}}{\left(\mei^2+|y|^2\right)^{\frac{n}{2}(1-\nu)}}\right)^{\crit-1-\pe}\,
dy\nonumber\\
&\leq&  C|x|^{1-n}\mei^{\frac{n}{2}}\int_{B_{\frac{R}{\mei}}(0)}
|y|^{1-s}\left(\frac{|y|^{1-\nu}}{\left(1+|y|^2\right)^{\frac{n}{2}(1-\nu)}}\right)^{\crit-1-\pe}\,
dy\nonumber\\
&\leq&  C'|x|^{1-n}\mei^{\frac{n}{2}}\label{ineq:last:1}
\eeqn
since $s\in (0,2)$ and up to taking $\nu>0$ small enough. Note that we 
have used here point (A7) of Proposition \ref{prop:exhaust}.\par

\smallskip\noindent We deal with the second integral. Note that when
$|x-y|\leq |x|/2$, we have that

$$\frac{|x|}{2}\leq |y|\leq \frac{3|x|}{2}.$$
Taking $\nu>0$ small enough, we then get that

\beqn
I_{2,\eps}(x)&\leq&
C|x|^{1-s}\left(\frac{\mei^{\frac{n}{2}-(n-1)\nu}|x|^{1-\nu}}{|x|^{n(1-\nu)}}\right)^{\crit-1-\pe}\int_{\{|x-y|\leq
|x|/2\}}|x-y|^{1-n}\, dy\nonumber\\
&\leq & C|x|^{1-n}\mei^{\frac{n}{2}}\cdot
|x|^{n+1-s}\mei^{-\frac{n}{2}}\left(\frac{\mei^{\frac{n}{2}-(n-1)\nu}|x|^{1-\nu}}{|x|^{n(1-\nu)}}\right)^{\crit-1-\pe}\nonumber\\
&\leq & C'|x|^{1-n}\mei^{\frac{n}{2}}\label{ineq:last:2}
\eeqn
since $|x|\geq \mei$. Plugging together (\ref{ineq:last:1}) and
(\ref{ineq:last:2}), we get that

$$\int_\Omega
\frac{G(x,y)}{|y|^s}\left(\frac{\mei^{\frac{n}{2}-(n-1)\nu}|y|^{1-\nu}}{\left(\mei^2+|y|^2\right)^{\frac{n}{2}(1-\nu)}}\right)^{\crit-1-\pe}\,
dy\leq  C|x|^{1-n}\mei^{\frac{n}{2}}$$
Since $|x|\geq \mei$, we get (\ref{ineq:bubble:1}).\hfill$\Box$

\medskip\noindent{\bf Step \ref{sec:fe:1}.6:} We let $i\in \{1,...,N\}$.
We claim that there exists $C>0$ such that
\beqn\label{ineq:bubble:2}
&&\int_\Omega
\frac{G(x,y)}{|y|^s}\left(\frac{\mei^{\frac{n}{2}-(n-1)\nu}|y|^{1-\nu}}{\left(\mei^2+|y|^2\right)^{\frac{n}{2}(1-\nu)}}\right)^{\crit-1-\pe}\,
dy\nonumber\\
&&\leq
C\frac{\mei^{\frac{n}{2}}|x|}{\left(\mei^2+|x|^2\right)^{\frac{n}{2}}}
\eeqn
for all $x\in\Omega$ such that $|x|\leq \mei$.

\smallskip\noindent{\it Proof of the Claim:} Indeed, let $p\in (1,n/s)$.
We let $\varphi_{\eps,i}\in H_{1,0}^p(\Omega)$ such that
\bequa\label{eq:phi:eps:i}
\Delta\varphi_{\eps,i}=\frac{1}{|x|^s}\left(\frac{\mei^{\frac{n}{2}-(n-1)\nu}|x|^{1-\nu}}{\left(\mei^2+|x|^2\right)^{\frac{n}{2}(1-\nu)}}\right)^{\crit-1-\pe}\hbox{
in }{\mathcal D}'(\Omega).
\eequa
We let $\varphi:U\to V$ defined in (\ref{def:vphi}) with $y_0=0$. We let
$$\tilde{\varphi}_{\eps,i}(x)=\mei^{\frac{n-2}{2}}\varphi_{\eps,i}\circ\varphi(\mei 
x)$$
for all $x\in \frac{U}{\mei}\cap\rnm$. We let $R>0$ such that
$\Omega\subset B_R(0)$. It follows from Green's representation formula and
the estimate (G5) on the Green's function that for any $x\in
\frac{U}{\mei}\cap\rnm$, we have that
\beq
|\tilde{\varphi}_{\eps,i}(x)|&\leq & \mei^{\frac{n-2}{2}}\int_\Omega
\frac{G(\varphi(\mei
x),y)}{|y|^s}\left(\frac{\mei^{\frac{n}{2}-(n-1)\nu}|y|^{1-\nu}}{\left(\mei^2+|y|^2\right)^{\frac{n}{2}(1-\nu)}}\right)^{\crit-1-\pe}\,
dy\\
&\leq & C\mei^{\frac{n-2}{2}} \int_\Omega \frac{1}{|\varphi(\mei
x)-y|^{n-2}|y|^s}\left(\frac{\mei^{\frac{n}{2}-(n-1)\nu}|y|^{1-\nu}}{\left(\mei^2+|y|^2\right)^{\frac{n}{2}(1-\nu)}}\right)^{\crit-1-\pe}\,
dy\\
&\leq & C\int_{B_{R/\mei}(0)} \frac{1}{\left|\frac{\varphi(\mei
x)}{\mei}-y\right|^{n-2}|y|^s}\left(\frac{|y|^{1-\nu}}{(1+|y|^2)^{\frac{n}{2}(1-\nu)}}\right)^{\crit-1-\pe}\,
dy.
\eeq
Since $s\in (0,2)$ and with the properties (\ref{def:vphi}) of $\varphi$,
we get that there exists $C>0$ such that
\bequa\label{bnd:tvarphi}
|\tilde{\varphi}_{\eps,i}(x)|\leq C
\eequa
for all $x\in B_3(0)\cap\rnm$ and all $\eps>0$. We let the metric
$(\tge)_{kl}=(\partial_k\varphi,\partial_l\varphi)(\mei x)$ for
$k,l=1,...,n$. Equation (\ref{eq:phi:eps:i}) rewrites as

$$\Delta_{\tge}\tilde{\varphi}_{\eps,i}=c_{i,\eps}\frac{\beta_{\eps,i}(x)^{(1-\nu)(\crit-1-\pe)-s}}{\left(1+\beta_{\eps,i}(x)^2\right)^{\frac{n}{2}(1-\nu)(\crit-1-\pe)}}\hbox{ 
in }{\mathcal D}'(B_3(0)\cap\rnm),$$
where $c_{i,\eps}\in\rr$ for all $\eps>0$ and $\lim_{\eps\to
0}c_{i,\eps}=c_i>0$.
$$\beta_{\eps,i}(x):=\left|\frac{\varphi(\mei x)}{\mei}\right|$$
for all $x\in B_3(0)\cap\rnm$. In particular, there exists $C>0$ such that
$$\frac{|x|}{C}\leq\beta_{\eps,i}(x)\leq C |x|$$
for all $x\in B_3(0)\cap\rnm$. Since (\ref{bnd:tvarphi}) holds, $s\in 
(0,2)$ and $\tilde{\varphi}_{\eps,i}\equiv 0$ on $\{x_1=0\}$, it follows 
from standard elliptic theory and the equation satisfies by 
$\tilde{\varphi}_{\eps,i}$ that there exists $C>0$ such that
$$\Vert\tilde{\varphi}_{\eps,i}\Vert_{C^{1}(B_2(0)\cap\overline{\rnm})}\leq
C$$
for all $\eps>0$. Since $\tilde{\varphi}_{\eps,i}(0)=0$, we get that
$$|\tilde{\varphi}_\eps^i(x)|\leq C |x|$$
for all $x\in B_2(0)\cap\rnm$ and all $\eps>0$. Coming back to the
definition of $\tilde{\varphi}_{\eps,i}$, we then get that there exists
$C>0$ such that
$$|\varphi_{\eps,i}(x)|\leq C
\frac{\mei^{\frac{n}{2}}|x|}{\left(\mei^2+|x|^2\right)^{\frac{n}{2}}}$$
for all $x\in\Omega\cap B_{\mei}(0)$. Inequality (\ref{ineq:bubble:2})
then follows from Green's representation formula.\hfill$\Box$

\medskip\noindent{\bf Step \ref{sec:fe:1}.7:} Plugging together
(\ref{ineq:He:alpha}), (\ref{ineq:bubble:1}) and (\ref{ineq:bubble:2}) 
into (\ref{ineq:ue:ef2}), we get that
$$|\ue(x)|\leq
C\sum_{i=1}^N\frac{\mei^{\frac{n}{2}}|x|}{\left(\mei^2+|x|^2\right)^{\frac{n}{2}}}+C |x|$$
for all $x\in\Omega$ and all $\eps>0$. This proves
(\ref{eq:fund:est:part1}). \hfill$\Box$

\section{Strong pointwise estimates, Part 2}\label{sec:fe:2}
This section is devoted to a refinement and a derivation of Proposition
\ref{prop:fund:est:part1}:

\begin{prop}\label{prop:fund:est:part2}
Let $\Omega$ be a smooth bounded domain of $\rn$, $n\geq 3$. We let $s\in 
(0,2)$. We let $(\pe)_{\eps>0}$ such that $\pe\in [0,\crit-2)$ for all 
$\eps>0$ and (\ref{lim:pe}) holds. We consider $(\ue)_{\eps>0}\in \huno$ 
such that (\ref{hyp:ae}), $(E_\eps)$ and (\ref{bnd:ue}) hold. We assume 
that blow-up occurs, that is

$$\lim_{\eps\to 0}\Vert\ue\Vert_{L^\infty(\Omega)}=+\infty.$$
We let $\mu_{\eps,1},...,\mu_{\eps,N}$ as in Proposition 
\ref{prop:exhaust}. Then, there exists $C>0$ such that

\beqn
&&|\ue(x)|\leq
C\sum_{i=1}^N\frac{\mei^{\frac{n}{2}}|x|}{\left(\mei^2+|x|^2\right)^{\frac{n}{2}}}+C|x|\label{est:co}\\
&&|\nabla\ue(x)|\leq
C\sum_{i=1}^N\frac{\mei^{\frac{n}{2}}}{\left(\mei^2+|x|^2\right)^{\frac{n}{2}}}+C\label{est:c1}
\eeqn
for all $\eps>0$ and all $x\in\Omega$.
\end{prop}

\noindent Inequality (\ref{est:co}) was proved in Proposition
\ref{prop:fund:est:part1}. We prove inequality (\ref{est:c1}). We let $G$ 
be the Green's function for the operator $\Delta$ on $\Omega$ with 
Dirichlet boundary condiction. Derivating Green's representation formula (\ref{eq:green:ue}) 
that

$$\nabla\ue(x)=\int_\Omega 
\nabla_xG(x,y)\left(\frac{|\ue(y)|^{\crit-2-\pe}\ue(y)}{|y|^s}-\ae(y)\ue(y)\right)\, 
dy$$
for all $x\in\Omega$ and all $\eps>0$. It then follows from (\ref{est:co}) 
that

\beqn
&&|\nabla\ue(x)|\leq C\int_\Omega 
|\nabla_xG(x,y)|\left(\frac{|\ue(y)|^{\crit-1-\pe}}{|y|^s}+1\right)\, 
dy\nonumber\\
&&\leq  C\sum_{i=1}^N\int_\Omega
\frac{|\nabla_xG(x,y)|}{|y|^s}\left(\frac{\mei^{\frac{n}{2}}|y|}{\left(\mei^2+|y|^2\right)^{\frac{n}{2}}}\right)^{\crit-1-\pe}\,
dy\nonumber\\
&&+C\int_\Omega |\nabla_xG(x,y)|\cdot(|y|^{\crit-1-s-\pe}+1)\, 
dy\label{proof:est:c1}
\eeqn
for all $x\in\Omega$ and all $\eps>0$.

\medskip\noindent{\bf Step \ref{sec:fe:2}.1:} We claim that there exists
$C>0$ such that

\bequa\label{sup:bubble:u0}
\int_\Omega |\nabla_xG(x,y)|\cdot\left(|y|^{\crit-1-s-\pe}+1\right)\, dy 
\leq C
\eequa
for all $x\in\Omega$ and all $\eps>0$.

\smallskip\noindent{\it Proof of the Claim:} Indeed, it follows from
property (G7) of Theorem \ref{th:green:1} that there exists $C>0$ such
that

\bequa\label{est:nablaG:sup}
|\nabla_xG(x,y)|\leq C |x-y|^{1-n}
\eequa
for all $x,y\in\Omega$ such that $x\neq y$. Since $s\in (0,2)$, we
then obtain that there exists $C>0$ such that

$$\int_\Omega |\nabla_xG(x,y)|\cdot\left(|y|^{\crit-1-s-\pe}+1\right)\, 
dy\leq C\int_\Omega
|x-y|^{1-n}\cdot\left(|y|^{\crit-1-s-\pe}+1\right)\, dy \leq C$$
for all $x\in\Omega$ and all $\eps>0$. This proves 
(\ref{sup:bubble:u0}).\hfill$\Box$

\medskip\noindent{\bf Step \ref{sec:fe:2}.2:} We let $i\in \{1,...,N\}$.
We claim that there exists $C>0$ such that
\bequa\label{sup:bubble:nablaG:1}
\int_\Omega
\frac{|\nabla_xG(x,y)|}{|y|^s}\left(\frac{\mei^{\frac{n}{2}}|y|}{\left(\mei^2+|y|^2\right)^{\frac{n}{2}}}\right)^{\crit-1-\pe}\,
dy\leq C\frac{\mei^{\frac{n}{2}}}{(\mei^2+|x|^2)^{\frac{n}{2}}}
\eequa
for all $x\in\Omega$ such that $|x|\leq\mei$ and all $\eps>0$.\par

\smallskip\noindent{\it Proof of the Claim:} We let
$\theta_\eps:=\frac{x}{\mei}$. Note that with our assumption, we have that
$|\theta_\eps|\leq 1$. We let $R>0$ such that $\Omega\subset B_R(0)$. With
(\ref{est:nablaG:sup}), and a change of variables, we get that

\beq
&&\int_\Omega
\frac{|\nabla_xG(x,y)|}{|y|^s}\left(\frac{\mei^{\frac{n}{2}}|y|}{\left(\mei^2+|y|^2\right)^{\frac{n}{2}}}\right)^{\crit-1-\pe}\,
dy\\
&&\leq
C\int_{B_R(0)}|x-y|^{1-n}|y|^{-s}\left(\frac{\mei^{\frac{n}{2}}|y|}{\left(\mei^2+|y|^2\right)^{\frac{n}{2}}}\right)^{\crit-1-\pe}\,
dy\\
&&\leq
C\mei^{-\frac{n}{2}}\int_{B_{\frac{R}{\mei}}(0)}\frac{|z|^{\crit-1-s-\pe}}{|\theta_\eps-z|^{n-1}(1+|z|^2)^{\frac{n}{2}(\crit-1-\pe)}}\,
dz
\eeq
Since $s\in (0,2)$ and $|\theta_\eps|\leq 1$, we get that there exists
$C>0$ such that

$$\int_\Omega
\frac{|\nabla_xG(x,y)|}{|y|^s}\left(\frac{\mei^{\frac{n}{2}}|y|}{\left(\mei^2+|y|^2\right)^{\frac{n}{2}}}\right)^{\crit-1-\pe}\,
dy\leq C\mei^{-\frac{n}{2}}.
$$
Since $|x|\leq\mei$, inequality (\ref{sup:bubble:nablaG:1})
follows.\hfill$\Box$

\medskip\noindent{\bf Step \ref{sec:fe:2}.3:} We let $i\in \{1,...,N\}$.
We claim that there exists $C>0$ such that

\bequa\label{sup:bubble:nablaG:2}
\int_\Omega
\frac{|\nabla_xG(x,y)|}{|y|^s}\left(\frac{\mei^{\frac{n}{2}}|y|}{\left(\mei^2+|y|^2\right)^{\frac{n}{2}}}\right)^{\crit-1-\pe}\,
dy\leq C\frac{\mei^{\frac{n}{2}}}{(\mei^2+|x|^2)^{\frac{n}{2}}}
\eequa
for all $x\in\Omega$ such that $|x|\geq\mei$ and all $\eps>0$.\par

\smallskip\noindent{\it Proof of the Claim:} We split the integral in two
parts:

\bequa\label{ineq:I}
\int_\Omega
\frac{|\nabla_xG(x,y)|}{|y|^s}\left(\frac{\mei^{\frac{n}{2}}|y|}{\left(\mei^2+|y|^2\right)^{\frac{n}{2}}}\right)^{\crit-1-\pe}\,
dy=I_{\eps,1}(x)+I_{\eps,2}(x)
\eequa
where
$$I_{\eps,j}(x)=\int_{\Omega_{\eps,j}(x)}
\frac{|\nabla_xG(x,y)|}{|y|^s}\left(\frac{\mei^{\frac{n}{2}}|y|}{\left(\mei^2+|y|^2\right)^{\frac{n}{2}}}\right)^{\crit-1-\pe}\,
dy$$
and
$$\Omega_{\eps,1}(x)=\Omega\cap\left\{|x-y|\geq\frac{|x|}{2}\right\}\hbox{
and }\Omega_{\eps,1}(x)=\Omega\cap \left\{|x-y|<\frac{|x|}{2}\right\}$$

\medskip\noindent{\it Step \ref{sec:fe:2}.3.1:} We deal with
$I_{\eps,1}(x)$. It follows from point (G8) of Theorem \ref{th:green:1}
that there exists $C>0$ such that

$$|\nabla_x G(x,y)|\leq C\frac{d(y,\partial\Omega)}{|x-y|^n}\leq
C\frac{|y|}{|x-y|^n}$$
for all $x,y\in\Omega$, $x\neq y$. We let $R>0$ such that $\Omega\subset
B_R(0)$. With a change of variable, we get that

\beq
I_{\eps,1}(x)&\leq &
C\int_{\Omega\cap\left\{|x-y|\geq\frac{|x|}{2}\right\}}\frac{|y|}{|x-y|^n|y|^s}\left(\frac{\mei^{\frac{n}{2}}|y|}{\left(\mei^2+|y|^2\right)^{\frac{n}{2}}}\right)^{\crit-1-\pe}\,
dy\\
&\leq &C 
|x|^{-n}\int_{B_R(0)}|y|^{1-s}\left(\frac{\mei^{\frac{n}{2}}|y|}{\left(\mei^2+|y|^2\right)^{\frac{n}{2}}}\right)^{\crit-1-\pe}\,
dy\\
&\leq & C 
|x|^{-n}\mei^{\frac{n}{2}}\int_{B_{\frac{R}{\mei}}(0)}\frac{|z|^{\crit-s-\pe}}{(1+|z|^2)^{\frac{n}{2}(\crit-1-\pe)}}\,
dz.
\eeq
Since $|x|\geq\mei$ and $s\in (0,2)$, we then get that

\bequa\label{ineq:I1}
I_{\eps,1}(x)\leq C |x|^{-n}\mei^{\frac{n}{2}}\leq C'
\frac{\mei^{\frac{n}{2}}}{(\mei^2+|x|^2)^{\frac{n}{2}}}.
\eequa

\medskip\noindent{\it Step \ref{sec:fe:2}.3.2:} We deal with
$I_{\eps,2}(x)$. As easily checked, we have that

\bequa\label{bnd:y:1}
\frac{|x|}{2}\leq |y|\leq \frac{3|x|}{2}
\eequa
for all $y\in \Omega_{\eps,2}(x)$. With (\ref{est:nablaG:sup}) and
(\ref{bnd:y:1}), we get that

\beq
I_{\eps,2}(x)&\leq &
\frac{C}{|x|^s}\left(\frac{\mei^{\frac{n}{2}}}{|x|^{n-1}}\right)^{\crit-1-\pe}\int_{\{|x-y|<\frac{|x|}{2}\}}|x-y|^{1-n}\,
dy\\
&\leq & C
|x|^{1-s}\left(\frac{\mei^{\frac{n}{2}}}{|x|^{n-1}}\right)^{\crit-1-\pe}\\
&\leq & C
|x|^{-n}\mei^{\frac{n}{2}}\frac{\mei^{\frac{n}{2}(\crit-2-\pe)}}{|x|^{n(\crit-1-\pe)-(\crit-1-\pe)-n-1+s}}.
\eeq
Since $|x|\geq\mei$ and $s\in (0,2)$, we then get that

\bequa\label{ineq:I2}
I_{\eps,2}(x)\leq C |x|^{-n}\mei^{\frac{n}{2}}\leq
C'\frac{\mei^{\frac{n}{2}}}{(\mei^2+|x|^2)^{\frac{n}{2}}}.
\eequa

\medskip\noindent Plugging (\ref{ineq:I1}) and (\ref{ineq:I2}) into
(\ref{ineq:I}), we get (\ref{sup:bubble:nablaG:2}).\hfill$\Box$

\medskip\noindent{\bf Step \ref{sec:fe:2}.4:} Plugging
(\ref{sup:bubble:u0}), (\ref{sup:bubble:nablaG:1}) and
(\ref{sup:bubble:nablaG:2}) into (\ref{proof:est:c1}), we get inequality
(\ref{est:c1}).

\section{Pohozaev identity and proof of compactness}\label{sec:poho}
This section is mainly devoted to the proof of the following proposition:

\begin{prop}\label{prop:poho} Let $\Omega$ be a smooth bounded domain of
$\rn$, $n\geq 3$, such that $0\in\partial\Omega$. We let $(\ue)$, $(\ae)$
and $(\pe)$ such that $(E_\eps)$, (\ref{hyp:ae}), (\ref{lim:pe}) and (\ref{bnd:ue}) hold. 
We assume that blow-up occurs, that is

$$\lim_{\eps\to 0}\Vert\ue\Vert_{L^\infty(\Omega)}=+\infty.$$
Then we have that

$$\lim_{\eps\to
0}\frac{\pe}{\meN}=\frac{(n-s)\displaystyle{\int_{\partial\rnm}}II_{0}(x,x)|\nabla\tu_N|^2\,
dx}{\displaystyle{(n-2)^2\alpha_N^{\frac{(n-1)(n-2)}{2(2-s)}}\sum_{i=1}^N}\alpha_{i}^{-\frac{(n-2)^2}{2(2-s)}}\displaystyle{\int_{\rnm}}|\nabla\tui|^2\,
dx}$$
when $n\geq 3$. In this expression, $II_0$ is the second fondamental form
at $0$ of the oriented boundary $\partial\Omega$ and $\partial\rnm$ is the 
oriented tangent space of $\partial\Omega$ at $0$. The sequences and 
families $\meN>0$, $\alpha_i$, $\tui$, $i\in \{1,...,N\}$ are as in 
Proposition \ref{prop:exhaust}. In addition, if $\ue\geq 0$ for all 
$\eps\geq 0$, we have that
$$\lim_{\eps\to
0}\frac{\pe}{\meN}=\frac{\displaystyle{(n-s)\int_{\partial\rnm}}|x|^2|\nabla\tu_N|^2\,
dx}{\displaystyle{n(n-2)^2\alpha_N^{\frac{(n-1)(n-2)}{2(2-s)}}\sum_{i=1}^N}\alpha_{i}^{-\frac{(n-2)^2}{2(2-s)}}\displaystyle{\int_{\rnm}}|\nabla\tui|^2\,
dx}\cdot H(0)$$
when $n\geq 3$. In this expression, $H(0)$ is the mean curvature at $0$ of 
the oriented boundary $\partial\Omega$.
\end{prop}

\noindent We prove the proposition in Steps \ref{sec:poho}.1 to
\ref{sec:poho}.3. We prove Theorem \ref{th:cpct} in Step \ref{sec:poho}.4.

\medskip\noindent{\bf Step \ref{sec:poho}.1:} We provide a Pohozaev-type
identity for $\ue$. It follows from Proposition \ref{prop:app} that
$\ue\in C^1(\overline{\Omega})$ and that $\Delta\ue\in L^p(\Omega)$ for
all $p\in (1,\frac{n}{s})$. We let

\bequa\label{def:re}
W_\eps:=\Omega\cap \varphi(B_{\re}(0)),\hbox{ where }\re=\sqrt{\meN}.
\eequa
In the sequel, we denote by $\nu(x)$ the outward normal vector at
$x\in\partial W_\eps$ of the oriented hypersurface $\partial W_\eps$
(oriented as the boundary of $W_\eps$). Integrating by parts, we get that

\beq
&&\int_{W_\eps} x^i\partial_i\ue\Delta\ue\, dx\\
&&=-\int_{\partial W_\eps}x^i\partial_i\ue\partial_{\nu}\ue\,
d\sigma+\int_{W_\eps}\partial_j(x^i\partial_i\ue)\partial_j\ue\, dx\\
&&=-\int_{\partial W_\eps}x^i\partial_i\ue\partial_{\nu}\ue\,
d\sigma+\int_{W_\eps}|\nabla\ue|^2\, dx+\int_{W_\eps}
x^i\partial_i\frac{|\nabla\ue|^2}{2}\, dx\\
&&=\left(1-\frac{n}{2}\right)\int_{W_\eps}|\nabla\ue|^2\,
dx+\int_{\partial
W_\eps}\left((x,\nu)\frac{|\nabla\ue|^2}{2}-x^i\partial_i\ue\partial_{\nu}\ue\right)\,
d\sigma\\
&&=\left(1-\frac{n}{2}\right)\left(\int_{\partial
W_\eps}\ue\partial_{\nu}\ue\, d\sigma+\int_{W_\eps}\ue\Delta\ue\,
dx\right)\\
&&+\int_{\partial
W_\eps}\left((x,\nu)\frac{|\nabla\ue|^2}{2}-x^i\partial_i\ue\partial_{\nu}\ue\right)\,
d\sigma.
\eeq
Using the equation $(E_\eps)$ in the RHS, we get that

\beqn
&&\int_{W_\eps} x^i\partial_i\ue\Delta\ue\,
dx=\left(1-\frac{n}{2}\right)\left(\int_{W_\eps}\frac{|\ue|^{\crit-\pe}}{|x|^s}\,
dx-\int_{W_\eps}\ae\ue^2\, dx\right)\nonumber\\
&&+\int_{\partial
W_\eps}\left(\left(1-\frac{n}{2}\right)\ue\partial_{\nu}\ue+(x,\nu)\frac{|\nabla\ue|^2}{2}-x^i\partial_i\ue\partial_{\nu}\ue\right)\,
d\sigma.\label{eq:poho:1}
\eeqn
On the other hand, using the equation $(E_\eps)$ satisfied by $\ue$, we
get that

\beqn
&&\int_{W_\eps} x^i\partial_i\ue\Delta\ue\, dx=\int_{W_\eps}
x^i\partial_i\ue\frac{|\ue|^{\crit-2-\eps}\ue}{|x|^s}\, dx-\int_{W_\eps}
x^i\partial_i\ue\ae\ue\, dx\nonumber\\
&&=\int_{W_\eps}
x^i|x|^{-s}\partial_i\left(\frac{|\ue|^{\crit-\pe}}{\crit-\pe}\right)\,
dx-\int_{W_\eps} x^i\partial_i\ue\ae\ue\, dx\nonumber\\
&&
=-\int_{W_\eps}\left(\partial_i(x^i|x|^{-s})\frac{|\ue|^{\crit-\pe}}{\crit-\pe}+x^i\partial_i\ue\ae\ue\right)\,
dx\nonumber\\
&&+\int_{\partial W_\eps}
\frac{(x,\nu)}{\crit-\pe}\cdot\frac{|\ue|^{\crit-\pe}}{|x|^s}\,
d\sigma\nonumber\\
&&
=-\int_{W_\eps}\frac{n-s}{|x|^s}\cdot\frac{|\ue|^{\crit-\pe}}{\crit-\pe}\,
dx+\frac{1}{2}\int_{W_\eps} (n\ae +x^i\partial_i\ae)\ue^2\, dx\nonumber\\
&&+\int_{\partial W_\eps}
\frac{(x,\nu)}{\crit-\pe}\cdot\frac{|\ue|^{\crit-\pe}}{|x|^s}\,
d\sigma-\int_{\partial W_\eps}\frac{(x,\nu)}{2}\ae\ue^2 \,
d\sigma.\label{eq:poho:2}
\eeqn
Plugging together (\ref{eq:poho:1}) and (\ref{eq:poho:2}), we get that

\beqn
&&\left(\frac{n-2}{2}-\frac{n-s}{\crit-\pe}\right)\int_{W_\eps}\frac{|\ue|^{\crit-\pe}}{|x|^s}\,
dx+\int_{W_\eps}\left(\ae+\frac{(x,\nabla\ae)}{2}\right)\ue^2\,
dx\nonumber\\
&&=\int_{\partial
W_\eps}\left(-\frac{n-2}{2}\ue\partial_\nu\ue+(x,\nu)\frac{|\nabla\ue|^2}{2}\right.\nonumber\\
&&\left.-x^i\partial_i\ue\partial_{\nu}\ue-\frac{(x,\nu)}{\crit-\pe}\cdot\frac{|\ue|^{\crit-\pe}}{|x|^s}\right)\,
d\sigma+\int_{\partial W_\eps}\frac{(x,\nu)}{2}\ae\ue^2\, dx\nonumber
\eeqn
for all $\eps>0$. Since

$$\partial W_\eps=[\varphi(B_{\re}(0))\cap\partial\Omega]\cup [\Omega\cap
\varphi(\partial B_{\re}(0))]$$
and since $\ue\equiv 0$ on $\partial\Omega$, we get that

\beqn
&&\frac{(n-2)\pe}{2\cdot(\crit-\pe)}\int_{\varphi(B_{\re}(0))\cap\Omega}\frac{|\ue|^{\crit-\pe}}{|x|^s}\,
dx-\int_{\varphi(B_{\re}(0))\cap\Omega}\left(\ae+\frac{(x,\nabla\ae)}{2}\right)\ue^2\,
dx\nonumber\\
&&=\frac{1}{2}\int_{\varphi(B_{\re}(0))\cap \partial
\Omega}(x,\nu)|\nabla\ue|^2\, d\sigma\label{eq:poho:3}\\
&&-\int_{\Omega\cap \varphi( \partial
B_{\re}(0))}\left(-\frac{n-2}{2}\ue\partial_\nu\ue+(x,\nu)\frac{|\nabla\ue|^2}{2}-x^i\partial_i\ue\partial_{\nu}\ue\right.\nonumber\\
&&\left.-\frac{(x,\nu)}{\crit-\pe}\cdot\frac{|\ue|^{\crit-\pe}}{|x|^s}+\frac{(x,\nu)}{2}\ae\ue^2\right)\,
d\sigma.\nonumber
\eeqn
It follows from (\ref{est:co}) and (\ref{est:c1}) that there exists $C>0$
such that

$$|\ue(x)|\leq C \re\hbox{ and }|\nabla\ue(x)|\leq C$$
for all $x\in \Omega\cap\varphi(\partial B_{\re}(0))$ (recall that
$\re=\sqrt{\meN}$). We then get that
\beqn
&&\int_{\Omega\cap\varphi(\partial
B_{\re}(0))}\left(-\frac{n-2}{2}\ue\partial_\nu\ue+(x,\nu)\frac{|\nabla\ue|^2}{2}-x^i\partial_i\ue\partial_{\nu}\ue\right.\nonumber\\
&&\left.-\frac{(x,\nu)}{\crit-\pe}\cdot\frac{|\ue|^{\crit-\pe}}{|x|^s}-\frac{(x,\nu)}{2}\ae\ue^2\right)\,
d\sigma=O(\meN^{\frac{n}{2}})=o(\meN)\label{eq:boundary:term}
\eeqn
when $\eps\to 0$ since $n\geq 3$. With (\ref{est:co}) and Proposition
\ref{prop:exhaust}, we get that

\beqn
&&\left|\int_{\varphi(B_{\re}(0))\cap\Omega}\left(\ae+\frac{(x,\nabla\ae)}{2}\right)\ue^2\,
dx\right|\leq C\int_{\varphi(B_{\re}(0))\cap\Omega}\ue^2\, dx\nonumber\\
&&\leq
C\sum_{i=1}^N\int_{\varphi(B_{\re}(0))\cap\Omega}\frac{\mei^{n}}{(\mei^2+|x|^2)^{n-1}}\,
dx+C\int_{\varphi(B_{\re}(0))\cap\Omega}|x|^2\, dx\nonumber\\
&&\leq C\sum_{1=1}^N\mei^2\int_{\rn}\frac{dx}{(1+|x|^2)^{n-1}}\, dx+C \re^{n+2}\nonumber\\
&&=o(\meN)\label{eq:full:term}
\eeqn
when $\eps\to 0$ since $n\geq 3$. Plugging (\ref{eq:boundary:term})
and (\ref{eq:full:term}) in (\ref{eq:poho:3}), we get that

\bequa\frac{(n-2)\pe}{2\cdot(\crit-\pe)}\int_{\varphi(B_{\re}(0))\cap\Omega}\frac{|\ue|^{\crit-\pe}}{|x|^s}\,
dx=\frac{1}{2}\int_{\varphi(B_{\re}(0))\cap \partial
\Omega}(x,\nu)|\nabla\ue|^2\, d\sigma+o(\meN)\label{eq:poho:5}
\eequa
when $\eps\to 0$ and $n\geq 3$.

\medskip\noindent{\bf Step \ref{sec:poho}.2:} We deal with the LHS of
(\ref{eq:poho:5}). We let $\varphi$ as in (\ref{def:vphi}). Since

$$\lim_{\eps\to 0}\frac{\re}{\meN}=+\infty$$
(see (\ref{def:re})), with a change of variables, we get for any
$R>\alpha>0$ that

\beqn\label{lim:0}
&&\int_{\varphi(B_{\re}(0))\cap\Omega}\frac{|\ue|^{\crit-\pe}}{|x|^s}\,
dx=\int_{\varphi(B_{\re}(0)\cap\rnm)}\frac{|\ue|^{\crit-\pe}}{|x|^s}\,
dx\\
&&=\int_{B_{\re}(0)\cap\rnm}
\frac{|\ue\circ\varphi(x)|^{\crit-\pe}}{|\varphi(x)|^s}\cdot|\hbox{Jac
}\varphi(x)|\, dx\nonumber\\
&&= \int_{B_{R\keun}(0)\cap\rnm}
\frac{|\ue\circ\varphi(x)|^{\crit-\pe}}{|\varphi(x)|^s}\cdot|\hbox{Jac
}\varphi(x)|\, dx\nonumber\\
&&+ \sum_{i=1}^{N-1}\left[\int_{(B_{R k_{\eps,i+1}}(0)\setminus
\overline{B}_{\alpha k_{\eps,i+1}}(0))\cap\rnm}
\frac{|\ue\circ\varphi(x)|^{\crit-\pe}}{|\varphi(x)|^s}\cdot|\hbox{Jac
}\varphi(x)|\, dx\right.\nonumber\\
&&\left.\qquad\qquad+\int_{(B_{\alpha k_{\eps,i+1}}(0)\setminus
\overline{B}_{ R k_{\eps,i}}(0))\cap\rnm}
\frac{|\ue\circ\varphi(x)|^{\crit-\pe}}{|\varphi(x)|^s}\cdot|\hbox{Jac
}\varphi(x)|\, dx\right]\nonumber\\
&&+\int_{(B_{\re}(0)\setminus \overline{B}_{R\keN}(0))\cap\rnm}
\frac{|\ue\circ\varphi(x)|^{\crit-\pe}}{|\varphi(x)|^s}\cdot|\hbox{Jac
}\varphi(x)|\, dx\nonumber
\eeqn
It follows from Proposition \ref{prop:exhaust} that

\bequa\label{lim:1}
\lim_{R\to +\infty}\lim_{\eps\to 0}\int_{B_{R\keun}(0)\cap\rnm}
\frac{|\ue\circ\varphi(x)|^{\crit-\pe}}{|\varphi(x)|^s}\cdot|\hbox{Jac
}\varphi(x)|\,
dx=\alpha_1^{-\frac{(n-2)^2}{2(2-s)}}\int_{\rnm}\frac{|\tu_1|^{\crit}}{|x|^s}\, 
dx
\eequa
and for any $i\in\{1,...,N-1\}$ that
\beqn
&&\lim_{R\to +\infty}\lim_{\alpha\to 0}\lim_{\eps\to 0}\int_{(B_{R
k_{\eps,i+1}}(0)\setminus \overline{B}_{\alpha k_{\eps,i+1}}(0))\cap\rnm}
\frac{|\ue\circ\varphi(x)|^{\crit-\pe}}{|\varphi(x)|^s}\cdot|\hbox{Jac
}\varphi(x)|\, dx\nonumber\\
&&=\alpha_{i+1}^{-\frac{(n-2)^2}{2(2-s)}}\int_{\rnm}\frac{|\tu_{i+1}|^{\crit}}{|x|^s}\,
dx.\label{lim:2}
\eeqn
It follows from the pointwise estimate (\ref{est:co}) that there exists
$C>0$ such that
$$|\ue(x)|\leq C\meN^{\frac{n}{2}}|x|^{1-n}+C|x|$$
for all $x\in\Omega$. It then follows that there exists $C>0$ independant
of $R>1$ such that

\beq
&&\int_{(B_{\re}(0)\setminus \overline{B}_{R\keN}(0))\cap\rnm}
\frac{|\ue\circ\varphi(x)|^{\crit-\pe}}{|\varphi(x)|^s}\cdot|\hbox{Jac
}\varphi(x)|\, dx\\
&&\leq C\int_{B_{\re}(0)\setminus
\overline{B}_{R\keN}(0)}\frac{1}{|y|^s}\left(|y|+\frac{\meN^{\frac{n}{2}}}{|y|^{n-1}}\right)^{\crit-\pe}\,
dy\\
&&\leq C\int_{B_{\re}(0)}|y|^{\crit-s-\pe}\, dy+
C\meN^{\frac{n}{2}\crit}\int_{B_{\re}(0)\setminus
\overline{B}_{R\keN}(0)}|y|^{-(n-1)(\crit-\pe)-s}\, dy\\
&&\leq C \re^{n}+\frac{C}{R^{(n-1)(\crit-\pe)-n+s}}.
\eeq
Since $\lim_{\eps\to 0}\re=0$, we get that

\bequa\label{lim:3}
\lim_{R\to +\infty}\lim_{\eps\to 0}\int_{(B_{\re}(0)\setminus
\overline{B}_{R\keN}(0))\cap\rnm}
\frac{|\ue\circ\varphi(x)|^{\crit-\pe}}{|\varphi(x)|^s}\cdot|\hbox{Jac
}\varphi(x)|\, dx=0.
\eequa
We let $i\in \{1,...,N-1\}$. Using the pointwise estimate (\ref{est:co}),
we get that

$$|\ue(x)|\leq C \frac{\mei^{\frac{n}{2}}}{|x|^{n-1}}+C \mu_{\eps,i+1}^{1-\frac{n}{2}}$$
for all $x\in\Omega$ and all $\eps>0$. With computations similar to the
ones provided for the proof of (\ref{lim:3}), we get that

\bequa\label{lim:4}
\lim_{R\to +\infty}\lim_{\alpha\to 0}\lim_{\eps\to 0}\int_{(B_{\alpha
k_{\eps,i+1}}(0)\setminus \overline{B}_{ R k_{\eps,i}}(0))\cap\rnm}
\frac{|\ue\circ\varphi(x)|^{\crit-\pe}}{|\varphi(x)|^s}\cdot|\hbox{Jac
}\varphi(x)|\, dx=0.
\eequa
Plugging together (\ref{lim:1}), (\ref{lim:2}), (\ref{lim:3}) and
(\ref{lim:4}) in (\ref{lim:0}), using point (A4) of Proposition
\ref{prop:exhaust}, we get that

\beqn
\lim_{\eps\to
0}\int_{\varphi(B_{\re}(0))\cap\Omega}\frac{|\ue|^{\crit-\pe}}{|x|^s}\,
dx&=&\sum_{i=1}^N\alpha_{i}^{-\frac{(n-2)^2}{2(2-s)}}\int_{\rnm}\frac{|\tui|^{\crit}}{|x|^s}\,
dx\nonumber\\
&=&\sum_{i=1}^N\alpha_{i}^{-\frac{(n-2)^2}{2(2-s)}}\int_{\rnm}|\nabla\tui|^2\,
dx\label{lim:lhs}
\eeqn

\medskip\noindent{\bf Step \ref{sec:poho}.3:} We deal with the RHS of
(\ref{eq:poho:5}). We have that

\beqn
&&\int_{\varphi(B_{\re}(0))\cap \partial \Omega}(x,\nu)|\nabla\ue|^2\,
d\sigma=
\int_{\varphi(B_{R\keun}(0))\cap\partial\Omega}(x,\nu)|\nabla\ue|^2\,
d\sigma\nonumber\\
&&+\sum_{i=1}^{N-2}\int_{\varphi(B_{R k_{\eps,i+1}}(0)\setminus
\overline{B}_{R k_{\eps,i}}(0))\cap\partial\Omega}(x,\nu)|\nabla\ue|^2\,
d\sigma\nonumber\\
&&+\int_{\varphi(B_{\alpha k_{\eps,N}}(0)\setminus \overline{B}_{R
k_{\eps,N-1}}(0))\cap\partial\Omega}(x,\nu)|\nabla\ue|^2\,
d\sigma\nonumber\\
&&+\int_{\varphi(B_{R k_{\eps,N}}(0)\setminus \overline{B}_{\alpha
k_{\eps,N}}(0))\cap\partial\Omega}(x,\nu)|\nabla\ue|^2\,
d\sigma\nonumber\\
&&+\int_{\varphi(B_{\re}(0)\setminus \overline{B}_{R
k_{\eps,N}}(0))\cap\partial\Omega}(x,\nu)|\nabla\ue|^2\,
d\sigma\label{est:bord:1}
\eeqn
Using the expression of $\varphi$ (see (\ref{def:vphi})), we get that
$$\nu(\varphi(x))=\frac{(1,-\partial_2\varphi_0(x), ...,
-\partial_n\varphi_0(x))}{\sqrt{1+\sum_{i=2}^n(\partial_i\varphi_0(x))^2}}$$
for all $x\in U\cap\{x_1=0\}$. With the expression of $\varphi$, we then
get that

\bequa\label{est:ps:1}
(\nu\circ\varphi(x),\varphi(x))=(1+O(1)|x|^2)\cdot\left(\varphi_0(x)-\sum_{i=2}^n
x^i\partial_i\varphi_0(x)\right)
\eequa
for all $x\in U\cap\{x_1=0\}$. In this expression, there exists $C>0$ such
that $|O(1)|\leq C$ for all $x\in U\cap\{x_1=0\}$. Since $\varphi_0(0)=0$
and $\nabla\varphi_0(0)=0$ (see (\ref{def:vphi})), we then get that there
exists $C>0$ such that

\bequa\label{est:ps:2}
|(\varphi(x),\nu\circ\varphi(x))|\leq C|x|^2
\eequa
for all $x\in U\cap\{x_1=0\}$.

\medskip\noindent{\it Step \ref{sec:poho}.3.1:} We deal with the second
term in the RHS of (\ref{est:bord:1}). We let $i\in\{1,...,N-2\}$. It
follows from the pointwise estimate (\ref{est:c1}) that

\bequa\label{est:gradient:bis}
|\nabla\ue(x)|\leq
C\mei^{\frac{n}{2}}|x|^{-n}+C\mu_{\eps,i+1}^{-\frac{n}{2}}
\eequa
for all $x\in\Omega$. With (\ref{est:ps:2}) and (\ref{est:gradient:bis}),
we get that

\beqn\label{final:est:1}
&&\int_{\varphi(B_{R k_{\eps,i+1}}(0)\setminus \overline{B}_{R
k_{\eps,i}}(0))\cap\partial\Omega}(x,\nu)|\nabla\ue|^2\,
d\sigma\nonumber\\
&&\leq C\int_{B_{2R k_{\eps,i+1}}(0)\setminus \overline{B}_{R
k_{\eps,i}/2}(0)\cap\{x_1=0\}}|x|^2\left(\mei^{n}|x|^{-2n}+\mu_{\eps,i+1}^{-n}\right)\,
dx\nonumber\\
&&\leq  C\mei+C\mu_{\eps,i+1}=o(\meN)
\eeqn
when $\eps\to 0$ when $n\geq 3$. Here, we have used that $i+1<N$ and point
(A3) of Proposition \ref{prop:exhaust}. With the same type of arguments,
we get that

\bequa\label{final:est:2}
\int_{\varphi(B_{R\keun}(0))\cap\partial\Omega}(x,\nu)|\nabla\ue|^2\,
d\sigma=o(\meN)
\eequa
when $\eps\to 0$ as soon as $N\geq 2$.

\medskip\noindent{\it Step \ref{sec:poho}.3.2:} We deal with the third
term of the RHS of (\ref{est:bord:1}). It follows from the pointwise
estimate (\ref{est:c1}) that

\bequa\label{est:gradient:ter}
|\nabla\ue(x)|\leq
C\mu_{\eps,N-1}^{\frac{n}{2}}|x|^{-n}+C\mu_{\eps,N}^{-\frac{n}{2}}
\eequa
for all $x\in\Omega$. With (\ref{est:ps:2}) and (\ref{est:gradient:ter}),
we get that

\beq
&&\int_{\varphi(B_{\alpha k_{\eps,N}}(0)\setminus \overline{B}_{R
k_{\eps,N-1}}(0))\cap\partial\Omega}(x,\nu)|\nabla\ue|^2\, d\sigma\\
&&\leq C\int_{B_{2\alpha k_{\eps,N}}(0)\setminus \overline{B}_{R
k_{\eps,N-1}/2}(0)\cap\{x_1=0\}}|x|^2\left(\mu_{\eps,N-1}^{n}|x|^{-2n}+\mu_{\eps,N}^{-n}\right)\,
dx\\
&&\leq  C\mu_{\eps,N-1}+C\alpha^{n+1} \mu_{\eps,N}
\eeq
since $n\geq 3$ and where $C>0$ is independant of $\alpha$ and $\eps>0$. With point
(A3) of Proposition \ref{prop:exhaust}, we get that

\bequa\label{final:est:3}
\lim_{\alpha\to 0}\lim_{\eps\to 0}\mu_{\eps,N}^{-1}\int_{\varphi(B_{\alpha
k_{\eps,N}}(0)\setminus \overline{B}_{R
k_{\eps,N-1}}(0))\cap\partial\Omega}(x,\nu)|\nabla\ue|^2\, d\sigma=0.
\eequa

\medskip\noindent{\it Step \ref{sec:poho}.3.3:} We deal with the fifth
term of the RHS of (\ref{est:bord:1}). It follows from the pointwise
estimate (\ref{est:c1}) that

\bequa\label{est:gradient:quat}
|\nabla\ue(x)|\leq C\mu_{\eps,N}^{\frac{n}{2}}|x|^{-n}+C
\eequa
for all $x\in\Omega$. With (\ref{est:ps:2}) and (\ref{est:gradient:quat}),
we get that

\beq
&&\int_{\varphi(B_{\re}(0)\setminus \overline{B}_{R
k_{\eps,N}}(0))\cap\partial\Omega}(x,\nu)|\nabla\ue|^2\, d\sigma\\
&&\leq C\int_{B_{2\re}(0)\setminus \overline{B}_{R
k_{\eps,N}/2}(0)\cap\{x_1=0\}}|x|^2\left(\mu_{\eps,N}^{n}|x|^{-2n}+C\right)\,
dx\\
&&\leq C R^{1-n} \mu_{\eps,N}+C\re^{n+1}
\eeq
since $n\geq 3$ and where $C>0$ is independant of $R$ and $\eps>0$. With the definition
(\ref{def:re}) of $\re$, we get that $\re^{n+1}=o(\meN)$ when $\eps\to 0$.
It then follows from point (A3) of Proposition \ref{prop:exhaust} that

\bequa\label{final:est:4}
\lim_{R\to +\infty}\lim_{\eps\to
0}\mu_{\eps,N}^{-1}\int_{\varphi(B_{\re}(0)\setminus \overline{B}_{R
k_{\eps,N}}(0))\cap\partial\Omega}(x,\nu)|\nabla\ue|^2\, d\sigma=0
\eequa
when $n\geq 3$.

\medskip\noindent{\it Step \ref{sec:poho}.3.4:} We deal with the fourth
term of the RHS of (\ref{est:bord:1}). Since $\varphi_0(0)=0$ and $\nabla\varphi_0(0)=0$, it follows from the definition (\ref{def:vphi}) of $\varphi$ and (\ref{est:ps:1})
that

\beqn
&&(\varphi(\keN x),\nu\circ\varphi(\keN x))\nonumber\\
&&=(1+O(\keN^2|x|^2))\left(\varphi_0(\keN
x)-\ke\sum_{i=2}^nx^i\partial_i\varphi_0(\keN x)\right)\nonumber\\
&&=-\frac{1}{2}\keN^2\sum_{i,j=2}^n\partial_{ij}\varphi_0(0)x^ix^j+\theta_{\eps,R}(x)\keN^2,\label{dev:PS}
\eeqn
for all $\eps>0$ and all $x\in B_R(0)\cap\{x_1=0\}$ and where
$\lim_{\eps\to 0}\sup_{B_R(0)\cap\{x_1=0\}}|\theta_{\eps,R}|=0$ for any
$R>0$. With a change of variable, (\ref{dev:PS}) and the definition of 
$\tilde{u}_{\eps,N}$ (see Proposition \ref{prop:exhaust}), we have that

\beq
&&\meN^{-1}\int_{\varphi(B_{R k_{\eps,N}}(0)\setminus \overline{B}_{\alpha
k_{\eps,N}}(0))\cap\partial\Omega}(x,\nu)|\nabla\ue|^2\, d\sigma\\
&&=\left(\frac{\keN}{\meN}\right)^{n-1}\left(-\frac{1}{2}\int_{(B_R(0)\setminus\overline{B}_\alpha(0))\cap\{x_1=0\}}\sum_{i,j=2}^nx^ix^j\partial_{ij}\varphi_0(x)|\nabla\tu_{\eps,N}|_{\tge}^2\,
dv_{\tge}\right)\\
&&+o(1)
\eeq
when $\eps\to 0$. In this expression,
$(\tge)_{ij}=(\partial_i\varphi,\partial_j\varphi)(\keN x)$ for all
$i,j=2,...,n$. It follows from (\ref{est:gradient:quat}) and the 
definition of $\tilde{u}_{\eps,N}$ that there exists $C>0$ such that
\bequa\label{upper:grnd:states}
|\nabla\tu_N(x)|\leq \frac{C}{1+|x|^n}
\eequa
for all $x\in\rnm$. With points (A4) and (A7) of Proposition 
\ref{prop:exhaust} and inequality (\ref{upper:grnd:states}), we get that
\beqn
&&\lim_{R\to +\infty}\lim_{\alpha\to 0}\lim_{\eps\to
0}\meN^{-1}\int_{\varphi(B_{R k_{\eps,N}}(0)\setminus \overline{B}_{\alpha
k_{\eps,N}}(0))\cap\partial\Omega}(x,\nu)|\nabla\ue|^2\,
d\sigma\nonumber\\
&&=-\frac{\alpha_N^{-\frac{n-1}{\crit-2}}}{2}\int_{\partial\rnm}\sum_{i,j=2}^nx^ix^j\partial_{ij}\varphi_0(x)|\nabla\tu_N|^2\,
dx\label{final:est:5}
\eeqn
when $n\geq 3$. Plugging (\ref{final:est:1}), (\ref{final:est:2}),
(\ref{final:est:3}), (\ref{final:est:4}) and (\ref{final:est:5}) in
(\ref{est:bord:1}), we get that

\bequa\label{lim:rhs}
\lim_{\eps\to 0}\meN^{-1}\int_{\varphi(B_{\re}(0))\cap \partial
\Omega}(x,\nu)|\nabla\ue|^2\,
d\sigma=-\frac{\alpha_N^{-\frac{n-1}{\crit-2}}}{2}\int_{\partial\rnm}\sum_{i,j=2}^nx^ix^j\partial_{ij}\varphi_0(x)|\nabla\tu_N|^2\,
dx
\eequa
We consider the second fondamental form associated to $\partial\Omega$, 
namely

$$II_p(x,y)=(d\nu_px,y)$$
for all $p\in\partial\Omega$ and all $x,y\in T_{p}\partial\Omega$ (recall
that $\nu$ is the outward normal vector at the hypersurface
$\partial\Omega$). In the canonical basis of
$\partial\rnm=T_0\partial\Omega$, the matrix of the bilinear form $II_{0}$
is $-D^2_0\varphi_0$, where $D^2_0\varphi_0$ is the Hessian matrix of
$\varphi_0$ at $0$. With this remark, plugging (\ref{lim:lhs}) and
(\ref{lim:rhs}) into (\ref{eq:poho:5}), we get that

\bequa\label{lim:pe:me}
\lim_{\eps\to
0}\frac{\pe}{\meN}=\frac{n-s}{(n-2)^2}\cdot\alpha_N^{-\frac{(n-1)(n-2)}{2(2-s)}}\cdot\frac{\displaystyle{\int_{\partial\rnm}}II_{0}(x,x)|\nabla\tu_N|^2\,
dx}{\displaystyle{\sum_{i=1}^N}\alpha_{i}^{-\frac{(n-2)^2}{2(2-s)}}\displaystyle{\int_{\rnm}}|\nabla\tui|^2\,
dx}
\eequa
when $n\geq 3$. This proves the first part of Proposition \ref{prop:poho}.

\medskip\noindent We prove the second part of the Proposition and assume 
that $\ue\geq 0$ for all $\eps$. It follows that the limit function 
$\tu_N$ is nonnegative, and then positive on $\rnm$. Moreover, we have 
that
$$\Delta\tu_N=\frac{\tu_N^{\crit-1}}{|x|^s}$$
in $\rnm$. It follows from (\ref{est:co}) that there exists $C>0$ such that
$$|\tu_N(x)|\leq\frac{C}{1+|x|^{n-1}}$$
for all $x\in\rnm$. It then follows from Proposition 
\ref{prop:sym} of Appendix C that there exists $v\in 
C^2(\rr_{-}^\star\times \rr)$ such that $\tu_N(x_1,x')=v(x_1,|x'|)$ for 
all $(x_1,x')\in \rr_-^\star\times \rr^{n-1}$. In particular, 
$|\nabla\tu_N|(0,x')$ is radially symmetrical wrt $x'\in\partial\rnm$. 
Since we have chosen a chart $\varphi$ that is Euclidean at $0$, we get 
that
\beq
\int_{\partial\rnm}II_{0}(x,x)|\nabla\tu_N|^2\,
dx&=&\frac{\sum_{i=2}^n(II_0)^{ii}}{n}\int_{\partial\rnm}|x|^2|\nabla\tu_N|^2\,
dx\\
&=&\frac{H(0)}{n}\int_{\partial\rnm}|x|^2|\nabla\tu_N|^2\,
dx.
\eeq
Note that we have used here that in the chart $\varphi$ defined in 
(\ref{def:vphi}), the matrix of the first fundamental form at $0$ is the 
identity. The second part of the Proposition then follows.

\medskip\noindent{\bf Step \ref{sec:poho}.4: Proof of Theorem
\ref{th:cpct}:} We let $(\ue)$, $(\ae)$ and $(\pe)$ such that $(E_\eps)$,
(\ref{hyp:ae}), (\ref{lim:pe}) and (\ref{bnd:ue}) hold. Assume that

\bequa\label{proof:th}
\lim_{\eps\to 0}\Vert\ue\Vert_{L^\infty(\Omega)}=+\infty.
\eequa
Then we can apply Proposition \ref{prop:poho}, and (\ref{lim:pe:me})
holds. Since the principal curvatures of $\partial\Omega$ at $0$ are
nonpositive, but do not all vanish, we have that $II_0(x,x)\leq 0$ for all
$x\in\partial\rnm$, but $II_0\not\equiv 0$. In particular, the RHS of
(\ref{lim:pe:me}) is negative. A contradiction since $\pe\geq 0$, and then
the LHS of (\ref{lim:pe:me}) is nonnegative. Then (\ref{proof:th}) does
not hold, and there exists $C>0$ such that $|\ue(x)|\leq C$ for all
$\eps>0$ and all $x\in\Omega$. The first part of Theorem \ref{th:cpct} then follows from
Proposition \ref{prop:bound}. In the case $\ue\geq 0$ for all $\eps>0$, we apply the second part of Proposition \ref{prop:poho} to recover compactness as soon as $H(0)<0$, and the second part of Theorem \ref{th:cpct} is proved.

\section{Proof of existence and multiplicity}\label{sec:proof:th}

\subsection{Proof of Theorem 1.1}
For any subcritical $p$, i.e.,  $2<p<{\crit}(s)$ we define the corresponding
best constant 
\bequa\label{def:musp}
\mu_{s,p} (\Omega) := \inf \left \{ \int_{\Omega}| \nabla u|^2
dx;\, u \in \huno \hbox{ and }  \int_{\Omega} \frac {|u|^p}{|x|^s}\, dx
=1\right\}.
\eequa
Because of the compactness of the embedding $\huno$ into $L^p(\Omega;
|x|^{-s}dx)$,  the infimum $\mu_{s,p} (\Omega)$ is attained at a positive extremal $v_p$
satisfying
\bequa\label{def:ext}
\left\{\begin{array}{ll}
\Delta u=\frac{u^{p-1}}{|x|^s}& \hbox{ in }{\mathcal D}'(\Omega)\\
u>0&\hbox{ in }\Omega\\
u=0&\hbox{ on }\partial\Omega.
\end{array}\right.
\eequa
Moreover, the family $(v_p)$ is uniformly bounded in $\huno$ when $p\to\crit$. Part 2 of the main compactness Theorem \ref{th:cpct} for positive 
sequences now yields a nontrivial limit $v$ that is an extremal for
$\mu_{s} (\Omega)$.

\subsection{Proof of Theorem 1.2}
For each  $2<p\leq {\crit}(s)$, consider the $C^2$-functional
\begin{equation}
I_p(u)=\frac{1}{2}\int_\Omega\vert \nabla u\vert^2\, dx-
\frac{1}{p}\int_{\Omega} \frac {|u|^p}{|x|^s}\, dx
\end{equation}
on $\huno$ whose critical points are the weak solutions of
\begin{equation}
\label{main}
\left\{ \begin{array}{ll}
\Delta u= \frac{|u|^{p-2}u}{|x|^s}  \hspace{5mm}& {\rm on} \ \Omega \\
  \ \    u  =   0 & {\rm on} \ \partial\Omega.
\end{array} \right.
\end{equation}
First note that for a fixed $u\in \huno$, we have since
$$I_p(\lambda u)=\frac{\lambda^2}{2}\int_\Omega\vert \nabla u\vert^2\, dx-
\frac{\lambda^p}{p}\int_{\Omega} \frac {|u|^p}{|x|^s}\, dx$$
that ${\rm limit}_{\lambda \to \infty}I(\lambda u)=-\infty$, which means 
that for each finite dimensional subspace $E_k \subset E:=\huno$, there 
exists $R_k>0$ such that
\begin{equation}
\label{neg}
\sup \{I_p(u); u\in E_k, \|u\| >R_k\} <0
\end{equation}
when $p\to\crit$. Let $(E_k)^\infty_{k=1}$ be an increasing sequence of subspaces of
$\huno$ such that
$\dim E_k=k $ and $\overline{\cup^\infty_{k=1} E_k}=E:=\huno$ and define
the min-max values:
$$c_{p,k} = {\ds \inf_{h \in {\bf H}_k} \sup_{x\in E_k} I_p(h(x))},$$
where
$$ {\bf H}_k=\{h \in C(E,E); \hbox{ $h$ is odd and $h(v)=v$
   for $\|v\| > R_k$ for some $R_k>0 \}$}.$$
\begin{prop}  With the above notation and assuming  $n\ge 3$, we have:
  \begin{enumerate}
  \item For each $k\in \nn$, $c_{p,k}>0$ and  $\lim\limits_{p\to \crit}
c_{p,k}=c_{\crit,k}:=c_k.$
   \item If $2<p<{\crit}$, there exists for each $k$, functions  $u_{p,k}\in
\huno$ such that $I'_p(u_{p,k})=0$, and $I_p(u_{p,k})=c_{p,k}$.
   \item  For each  $2<p < {\crit}$, we have $c_{p,k}$ satisfy
$ c_{p,k}\ge D_{n,p} k^{\frac{p+1}{p-1}\frac{2}{n}}$ where $D_{n,p}>0$ is 
such that
$\lim_{p\to  {\crit}}D_{n,p}=0$.
\item $\lim\limits_{k\to \infty} c_k=\lim\limits_{k\to \infty} 
c_{\crit,k}=+\infty$.

\end{enumerate}
\end{prop}
\noindent {\bf Proof:} (1) First note that in view of the Hardy-Sobolev 
inequality, we have
\[
I_p(u)\geq \frac{1}{2}\|\nabla u\|_2^2-C\|\nabla u\|_2^p =\|\nabla 
u\|_2^2\left(\frac{1}{2}-C\|\nabla u\|_2^{p-2}\right)\geq \alpha >0
\]
provided   $\|u\|_{\huno}= \rho$ for some $\rho>0$ small enough.  A 
standard intersection lemma gives that
the sphere $S_\rho=\{u\in E; \|u\|_{\huno}= \rho\}$ must intersect every 
image $h(E_k)$  by an odd continuous function $h$. It follows that
\[
c_{p,k}\geq \inf \{I_p(u); u\in S_\rho \} \geq \alpha >0.
\]
In view of (\ref{neg}), it follows that for each $h\in {\bf H}_k$, we have 
that
\[
\sup\limits_{x\in
E_k}I_{p_i}(h(x))=\sup\limits_{x\in D_k}I_p(h(x))
\]
  where $D_k$ denotes the
ball in $E_k$ of radius $R_k$.
  Consider now a sequence $p_i \to \crit$ and note
first that for each $u\in E$,  we have that $I_{p_i}(u) \to I_{\crit}
(u)$.   Since $h(D_k)$ is compact and the family of
functionals $(I_p)_p$ is equicontinuous, it follows that
   $\sup\limits_{x\in E_k}I_p(h(x))\to \sup\limits_{x\in E_k}I_{\crit
}(h(x))$, from which follows that
   $\limsup\limits_{i\in \nn}c_{p_i,k}\leq \sup\limits_{x\in E_k}I_{\crit}
(h(x))$. Since this holds for any $h\in {\bf H}_k$, it follows that
   \[
    \limsup\limits_{i\in \nn}c_{p_i,k}\leq c_{\crit, k}=c_k.
   \]
On the other hand, the function $f(r)=\frac{1}{p}r^p-\frac{1}{2^*}r^{2^*}$ 
attains its maximum on $[0, +\infty)$ at $r=1$ and therefore
$h(r) \leq \frac{1}{p}-\frac{1}{2^*}$ for all $r>0$. It follows
   \[
   I_{\crit}(u) = I_{p}(u)+\int_\Omega \frac{1}{|x|^s}
\left(\frac{1}{p}|u(x)|^p-\frac{1}{\crit}|u(x)|^{\crit }\right)\, dx\leq
I_{p}(u)+\int_\Omega \frac{1}{|x|^s} \left(\frac{1}{p}-\frac{1}{\crit}\right)\, dx
   \]
from which follows that $c_k\leq  \liminf\limits_{i\in \nn}c_{p_i,k}$, and 
claim (1) is proved. \\
If now $p< {\crit}$, we are in the subcritical case, that is we have
compactness in the Sobolev embedding $\huno \to L^p(\Omega; |x|^{-s}dx)$ 
and
therefore $I_p$ has the Palais-Smale condition. It is then standard to
find critical points $u_{p,k}$ for $I_p$ at each level $c_{p,k}$ (see for
example \cite{Gh}).
Now there are many ways to establish growth estimates for $c_{p,k}$ as
$k\to +\infty$, and we shall use here the one based on the  Morse indices
of these variationally obtained solutions, a method first used by
by Bahri-Lions \cite{bl} and independently by Tanaka \cite{Ta}. We need
the following key estimate of Li-Yau \cite{LY}.
    \begin{lem}   Let  $V \in L^{n/2}(\Omega)$ and denote by $m^*\,(V)$ the
number of non-positive eigenvalues of the following eigenvalue problem:
$$
\left\{ \begin{array}{ll}
     \Delta u -Vu=\lambda u   \hspace{5mm}& {\rm on } 
       \Omega \\
     \ \ \ \ \quad \quad \quad    u  =   0 & {\rm on} \ \partial\Omega. \\
\end{array} \right.
$$
If $n\ge 3$, then there is a constant $C_n>0$ such that
$m^*\,(V)\le C_n\| V\|^{n/2}_{n/2}$.
\end{lem}
To prove the growth estimates on the critical values $c_{p,k}$, one can
follow \cite{Ta} (see also \cite{Gh}) and identify a cohomotopic family of
sets ${\bf F_k}$ of dimension $k$ in such a way that if $D_k$ denotes the 
ball in $E_k$ of radius $R_k$ and if $\gamma\in {\bf H}_k$, then
$\gamma (D_k) \in {\bf F_k}$. It then follows that there exists 
$v_{p,k}\in
\huno$ such that $I_p(v_{p,k})\le c_{p,k}$, $I^\prime(v_{p,k})=0$ and
$m^*(v_{p,k})\ge k $, where $m^*(v_{p,k})$ is the augmented Morse index of 
$I_p$ at $v_{p,k}$. In other words, since
\begin{eqnarray*}
I_p''(v)(h, h) &= \int\limits_\Omega |\nabla h|^2\, dx-(p-1)
\int\limits_\Omega \frac{|v|^{p-2}}{|x|^s} h^2\, dx\\
&= \langle (\Delta -(p-1)\frac{|v|^{p-2}}{|x|^s} )h, h\rangle,
\end{eqnarray*} 
in $H^{-1}(\Omega)$, this means that the operator $(\Delta 
-(p-1)\frac{|v_{p,k}|^{p-2}}{|x|^s} )$
 possesses at least $k$ non-positive eigenvalues. Applying the above lemma, 
we get that the number of these non-positive eigenvalues  is bounded above 
by
$C_n\int\limits_\Omega
{[(p-1)\frac{|v_{p,k}|^{p-2}}{|x|^s}]}^{\frac{n}{2}}\, dx)$. \\
Since $p<\frac{2n}{n-2}$, we have $q:=\frac{2p}{n(p-2)}>1$, as well as its
conjugate $q'$. Moreover, since $p<\frac{2(n-s)}{n-2}$, we have that 
$\frac{2sn}{2p-np+2n}<n$. It then follows from Holder's inequality that:
\begin{eqnarray}
k &\leq& C_n\int\limits_\Omega
|p-1|^{\frac{n}{2}}\frac{|v_{p,k}|^{(p-2)\frac{n}{2}}}{|x|^{s\frac{n}{2}}}
\, dx\label{ineq:topo}\\
&\leq & C_n |p-1|^\frac{n}{2}\left(\int\limits_\Omega\frac{1}{|x|^{\frac{2sn}{2p-np+2n}}}\, dx\right)^{\frac{1}{q'}}{\left(\int\limits_\Omega \frac{|v_{p,k}|^p}{|x|^s}\,
dx\right)}^{\frac{n(p-2)}{2p}}  \nonumber\\
&\leq& C_{n,p} \left(\int\limits_\Omega\frac{|v_{p,k}|^p}{|x|^s}\,
dx\right)}^{\frac{n(p-2)}{2p}\nonumber
\end{eqnarray}
where $C_{n,p}= C_n |p-1|^\frac{n}{2}\left(\int\limits_\Omega
\frac{1}{|x|^{\frac{2sn}{2p-np+2n}}}\, dx\right)^{\frac{1}{q'}}$. \\
Since $\langle I'(v_{p,k}), v_{p,k}\rangle=0$, it follows that
$
\int\limits_\Omega |\nabla v_{p,k}|^2\, dx=\int\limits_\Omega
\frac{|v_{p,k}|^p}{|x|^s}\, dx,
$
which finally implies that
\begin{eqnarray*}
c_{p,k}\geq I(v_{p,k}) &=& \frac{1}{2}\int\limits_\Omega |\nabla
v_{p,k}|^2\, dx-\frac{1}{p}\int\limits_\Omega \frac{|v_{p,k}|^p}{|x|^s}\,
dx\\
&=&
\left(\frac{1}{2}-\frac{1}{p}\right)\int\limits_\Omega\frac{|v_{p,k}|^p}{|x|^s}\,
dx\\
&\geq& D_{n,p}k^{\frac{2p}{n(p-2)}},
\end{eqnarray*}
where $D_{n,p}=(\frac{1}{2}-\frac{1}{p})C_{n,p}^{-{\frac{2p}{n(p-2)}}  }.$ 
\\

To prove 4) we proceed by contradiction and assume that $(c_k)_k$ is 
bounded so that a subsequence
of which converges to some real number $c$. Using  the first claim of the 
proposition, there exists for each $k \in \nn$, $2<p_k<\crit$ such that
$|c_{p_{k}, k}-c_k|<\frac{1}{k}$ in such a way that $\lim_{k\to +\infty}p_k=\crit$ and 
\begin{equation}
\label{critical.values.limit}
\lim\limits_{k\to +\infty}c_{p_{k},k}=\lim\limits_{k\to +\infty}c_k=c.
\end{equation}
As above, there exists $v_{p_k,k}\in
\huno$ such that $I_{p_k}(v_{p_k,k})\le c_{p_k,k}$, $I_{p_k}^\prime(v_{p_k,k})=0$ 
and
$m^*(v_{p_k,k})\ge k $, where $m^*(v_{p_k,k})$ is the augmented Morse 
index of $I_{p_k}$ at $v_{p_k,k}$.

But  (\ref{critical.values.limit}) gives that the energies of 
$(v_{p_k,k})_k$  are uniformly bounded and therefore $(v_{p_k,k})_k$ is 
bounded in $\huno$.  It follows from Proposition \ref{prop:app} and the 
compactness Theorem \ref{th:cpct} that they converge to a solution $v$ of 
(\ref{main}) with energy below level $c$. In particular, there exists $C>0$ such that 
\bequa\label{est:pt:vk}
|v_{p_k,k}(x)|\leq C
\eequa
for all $x\in\Omega$ and all $k\in\nn$. With (\ref{ineq:topo}) applied to $v_{p_k,k}$, we get that
$$k \leq C_n\int\limits_\Omega
|p_k-1|^{\frac{n}{2}}\frac{|v_{p_k,k}|^{(p-2)\frac{n}{2}}}{|x|^{s\frac{n}{2}}}
\, dx.$$
With (\ref{est:pt:vk}), we get that there exists a constant $C>0$ independant of $k$ such that
$$k\leq C\int_{\Omega}\frac{dx}{|x|^{\frac{sn}{2}}}.$$
In particular, since $s\in (0,2)$, the integral is finite and there existe $C>0$ such that $k\leq C$ for all $k\in\nn$. A contradiction, and we are done with the 
proposition. \\

To complete the proof of Theorem 1.3,   notice that since for each $k$, we
have $ \lim\limits_{p_i\to \crit}I_{p_i}(u_{p_i,k})=\lim\limits_{p_i\to\crit}c_{p_i,k}=c_k$, it follows that the sequence  $(u_{p_i,k})_i $ is
uniformly bounded in $\huno$. Moreover, since $I_{p_i}'(u_{p_i,k})=0$, it 
follows from Proposition \ref{prop:app} and the compactness Theorem 
\ref{th:cpct} that by letting $p_i\to \crit$, we
get a solution $u_k$ of (\ref{main}) in such a way that
$I_{\crit}(u_k)=\lim\limits_{p\to \crit}I_p(u_{p,k})=\lim\limits_{p\to
\crit}c_{p,k}=c_k$. Since the latter sequence goes to infinity, it follows 
that   (\ref{main}) has  an infinite number of critical levels. The result 
for the equation $\Delta u+au=\frac{|u|^{\crit-2}u}{|x|^s}$ when 
$\Delta+a$ is coercive goes the same way, and Theorem \ref{th:multi} is 
proved.

\section{Appendix A: Regularity of weak solutions}\label{sec:app}
In this appendix, we prove the following regularity result. Note that such a $C^1-$regularity was first proved out by Egnell \cite{eg}. We include the proof for completeness.

\begin{prop}\label{prop:app}
Let $\Omega$ be a smooth domain of $\rn$, $n\geq 3$. We assume that either
$\Omega$ is bounded, or $\Omega=\rnm$. We let $s\in (0,2)$ and $a\in
C^0(\overline{\Omega})$. We let $\eps\in [0,\crit-2)$ and consider
$u\in\huno$ a weak solution of
$$\Delta u+au=\frac{|u|^{\crit-2-\eps}u}{|x|^s}  \hbox{ in }{\mathcal
D}'(\Omega).$$
Then $u\in C^{1,\theta}(\overline{\Omega})$ for all 
$\theta\in(0,\min\{1,\crit-\eps-s\})$ if $\Omega$ is bounded, and $u\in 
C^{1,\theta}_{loc}(\overline{\rnm})$ for all 
$\theta\in(0,\min\{1,\crit-\eps-s\})$ if $\Omega=\rnm$. In addition, in 
all the cases, we have that $u\in C^2(\overline{\Omega}\setminus\{0\})$ if 
$a\in C^{0,\alpha}(\overline{\Omega})$ for some $\alpha\in (0,1)$.
\end{prop}

\begin{proof} We prove the result when $\Omega$ is bounded. The arguments
and the results are basically local, and the proof goes the same way when
$\Omega=\rnm$.

\medskip\noindent{\bf Step \ref{sec:app}.1:} We follow the strategy
developed by Trudinger (\cite{trudinger}, and \cite{hebeybook1} for an 
exposition in book form). Let $\beta\geq 1$, and $L>0$. We let

$$G_L(t)=\left\{\begin{array}{ll}
|t|^{\beta-1}t& \hbox{ if }|t|\leq L\\
\beta L^{\beta-1}(t-L)+L^\beta & \hbox{ if }t\geq L\\
\beta L^{\beta-1}(t+L)-L^\beta & \hbox{ if }t\leq -L
\end{array}\right.$$
and
$$H_L(t)=\left\{\begin{array}{ll}
|t|^{\frac{\beta-1}{2}}t& \hbox{ if }|t|\leq L\\
\frac{\beta+1}{2}L^{\frac{\beta-1}{2}}(t-L)+L^{\frac{\beta+1}{2}} & \hbox{
if }t\geq L\\
\frac{\beta+1}{2}L^{\frac{\beta-1}{2}}(t+L)-L^{\frac{\beta+1}{2}} & \hbox{
if }t\leq -L
\end{array}\right.$$
As easily checked,

$$0\leq t G_L(t)\leq H_L(t)^2\hbox{ and
}G_L'(t)=\frac{4\beta}{(\beta+1)^2}(H_L'(t))^2$$
for all $t\in\rr$ and all $L>0$.
Let $\eta\in C^\infty_c(\rn)$. As easily checked, $\eta^2G_L(u),\eta
H_L(u)\in\huno$. With the equation verified by $u$, we get that
\bequa\label{eq:C}
\int_\Omega\nabla u\nabla (\eta^2 G_L(u))\, dx=\int_\Omega
\frac{|u|^{\crit-2-\eps}}{|x|^s}\eta^2 u G_L(u)\, dx-\int_\Omega a\eta^2
uG_L(u)\, dx.
\eequa
We let $J_L(t)=\int_0^t G_L(\tau)\, d\tau$ for all $t\in\rr$. Integrating
by parts, we get that

\beqn
&&\int_\Omega\nabla u\nabla (\eta^2 G_L(u))\, dx=\int_\Omega\eta^2
G_L'(u)|\nabla u|^2\, dx+\int_\Omega \nabla\eta^2\nabla J_L(u)\,
dx\nonumber\\
&&=\frac{4\beta}{(\beta+1)^2}\int_\Omega \eta^2 |\nabla H_L(u)|^2\,
dx+\int_\Omega(\Delta\eta^2)J_L(u)\, dx\nonumber\\
&&= \frac{4\beta}{(\beta+1)^2}\int_\Omega |\nabla (\eta H_L(u))|^2\,
dx-\frac{4\beta}{(\beta+1)^2}\int_\Omega\eta\Delta\eta |H_L(u)|^2\,
dx\nonumber\\
&&+\int_\Omega(\Delta\eta^2)J_L(u)\, dx\label{eq:A}
\eeqn
On the other hand, with H\"older's inequality and the definition of
$\mu_{s}(\rn)$, we get that

\beqn
&&\int_\Omega \left(\frac{|u|^{\crit-2-\eps}}{|x|^s}-a\right)\cdot \eta^2
u G_L(u)\, dx\leq \int_\Omega\left( |a|+
\frac{|u|^{\crit-2-\eps}}{|x|^s}\right)\cdot(\eta H_L(u))^2\,
dx\nonumber\\
&&\leq \left(\int_{\Omega\cap\hbox{Supp }\eta}\frac{(|a|\cdot
|x|^s+|u|^{\crit-2-\eps})^{\frac{\crit-\eps}{\crit-2-\eps}}}{|x|^s}\, dx\right)^{1-\frac{2}{\crit-\eps}}\nonumber\\
&&\times \left(\int_{\Omega}\frac{|\eta
H_L(u)|^{\crit}}{|x|^s}\, dx\right)^{\frac{2}{\crit}}\times\left(\int_{\Omega\cap\hbox{Supp }\eta}\frac{dx}{|x|^s}\right)^{\frac{2\eps}{\crit(\crit-\eps)}}\nonumber\\
&&\leq \alpha\int_{\Omega}|\nabla(\eta H_L(u))|^2\, dx\label{eq:B}
\eeqn
where
\begin{eqnarray*}
\alpha&:=&\left(\int_{\Omega\cap\hbox{Supp }\eta}\frac{(|a|\cdot
|x|^s+|u|^{\crit-2-\eps})^{\frac{\crit-\eps}{\crit-2-\eps}}}{|x|^s}\,
dx\right)^{1-\frac{2}{\crit-\eps}}\\
&&\times\mu_{s}(\rn)^{-1}\left(\int_{\Omega\cap\hbox{Supp 
}\eta}\frac{dx}{|x|^s}\right)^{\frac{2\eps}{\crit(\crit-\eps)}}
\end{eqnarray*}
Plugging (\ref{eq:A}) and (\ref{eq:B}) into (\ref{eq:C}), we get that

\bequa\label{eq:D}
A\cdot \int_\Omega |\nabla (\eta H_L(u))|^2\, dx\leq
\frac{4\beta}{(\beta+1)^2}\int_\Omega|\eta\Delta\eta| |H_L(u)|^2\,
dx+\int_\Omega|\Delta(\eta^2)J_L(u)|\, dx
\eequa
where
\begin{eqnarray*}
A&:=&\frac{4\beta}{(\beta+1)^2} -\left(\int_{\Omega\cap\hbox{Supp
}\eta}\frac{(|a|\cdot
|x|^s+|u|^{\crit-2-\eps})^{\frac{\crit-\eps}{\crit-2-\eps}}}{|x|^s}\,
dx\right)^{1-\frac{2}{\crit-\eps}}\\
&&\times
\mu_{s}(\rn)^{-1}\left(\int_{\Omega\cap\hbox{Supp 
}\eta}\frac{dx}{|x|^s}\right)^{\frac{2\eps}{\crit(\crit-\eps)}}
\end{eqnarray*}

\medskip\noindent{\bf Step \ref{sec:app}.2:} We let

$$p_0=\sup\{p\geq 1/\, u\in L^p(\Omega)\}.$$
It follows from Sobolev's embedding theorem that $p_0\geq\frac{2n}{n-2}$.
We claim that

$$p_0=+\infty.$$
We proceed by contradiction and assume that

$$p_0<\infty.$$
Let $p\in (2,p_0)$. It follows from the definition of $p_0$ that $u\in
L^p(\Omega)$. Let $\beta=p-1>1$. For any $x\in\overline{\Omega}$, we let
$\delta_x>0$ such that

\begin{eqnarray}
&&\left(\int_{\Omega\cap B_{2\delta_x}(x)}\frac{(|a|\cdot
|y|^s+|u|^{\crit-2-\eps})^{\frac{\crit-\eps}{\crit-2-\eps}}}{|y|^s}\,
dy\right)^{1-\frac{2}{\crit-\eps}}\mu_{s}(\rn)^{-1}\nonumber\\
&&\times
\left(\int_{\Omega\cap B_{2\delta_x}(x)}\frac{dy}{|y|^s}\right)^{\frac{2\eps}{\crit(\crit-\eps)}}\leq
\frac{2\beta}{(\beta+1)^2}.\label{eq:E}
\end{eqnarray}
Since $\overline{\Omega}$ is compact, we get that there exists
$x_1,...,x_N\in\overline{\Omega}$ such that
$$\Omega\subset \bigcup_{i=1}^N B_{\delta_{x_i}}(x_i).$$
We fix $i\in\{1,...,N\}$ and let $\eta\in
C^\infty(B_{2\delta_{x_i}}(x_i))$ such that $\eta(x)=1$ for all $x\in
B_{\delta_{x_i}}(x_i)$. We then get with (\ref{eq:D}) and (\ref{eq:E})
that

\beqn
&&\frac{2\beta}{(\beta+1)^2}\int_\Omega |\nabla (\eta H_L(u))|^2\,
dx\nonumber\\
&&\leq \frac{4\beta}{(\beta+1)^2}\int_\Omega|\eta\Delta\eta| |H_L(u)|^2\,
dx+\int_\Omega|\Delta\eta^2| |J_L(u)|\, dx.\label{eq:F}
\eeqn
Recall that it follows from Sobolev's inequality that there exists
$K(n,2)>0$ that depends only on $n$ such that

\bequa\label{ineq:sob}
\left(\int_{\rn}|f|^{\frac{2n}{n-2}}\, dx\right)^{\frac{n-2}{n}}\leq
K(n,2)\int_{\rn}|\nabla f|^2\, dx
\eequa
for all $f\in \hunrn$. It follows from (\ref{eq:F}) and (\ref{ineq:sob})
that

\beq
&&\frac{2\beta}{(\beta+1)^2}K(n,2)^{-1} \left(\int_\Omega |\eta
H_L(u)|^{\frac{2n}{n-2}}\, dx\right)^{\frac{n-2}{n}}\\
&&\leq \frac{4\beta}{(\beta+1)^2}\int_\Omega|\eta\Delta\eta| |H_L(u)|^2\,
dx+\int_\Omega|\Delta\eta^2|\cdot |J_L(u)|\, dx
\eeq
for all $L>0$. As easily checked, there exists $C_0>0$ such that
$|J_L(t)|\leq C_0 |t|^{\beta+1}$ for all $t\in\rr$ and all $L>0$.
Since $u\in L^{\beta+1}(\Omega)$, we get that there exists a constant
$C=C(\eta,u,\beta,\Omega)$ independant of $L$ such that

$$\int_{\Omega\cap B_{\delta_{x_i}}(x_i)} |H_L(u)|^{\frac{2n}{n-2}}\,
dx\leq \int_\Omega |\eta H_L(u))|^{\frac{2n}{n-2}}\, dx\leq C$$
for all $L>0$. Letting $L\to +\infty$, we get that
$$\int_{\Omega\cap B_{\delta_{x_i}}(x_i)} |u|^{\frac{n}{n-2}(\beta+1)}\,
dx<+\infty,$$
for all $i=1...N$. We then get that $u\in
L^{\frac{n}{n-2}(\beta+1)}(\Omega)=L^{\frac{n}{n-2}p}(\Omega)$. And then,
$\frac{n}{n-2}p\leq p_0$ for all $p\in (2,p_0)$. Letting $p\to p_0$, we
get a contradiction. Then $p_0=+\infty$ and $u\in L^p(\Omega)$ for all
$p\geq 1$. This ends Step \ref{sec:app}.2.

\medskip\noindent{\bf Step \ref{sec:app}.3:} We claim that

$$u\in C^{0,\alpha}(\overline{\Omega})$$
for all $\alpha\in (0,1)$. Indeed, it follows from Step \ref{sec:app}.2
and the assumption $0<s<2$ that there exists $p>\frac{n}{2}$ such that

$$\fe:=\frac{|u|^{\crit-2-\eps}u}{|x|^s}-au\in L^p(\Omega).$$
It follows from standard elliptic theory that, in this case, $u\in
C^{0,\alpha}(\overline{\Omega})$ for all $\alpha\in (0,\min\{2-s,1\})$. We
let

$$\alpha_0=\sup\{\alpha\in (0,1)/\, u\in
C^{0,\alpha}(\overline{\Omega})\}.$$
Note that it follows from the preceding remark that $\alpha_0>0$. We let $\alpha\in (0,\alpha_0)$. Then $u\in
C^{0,\alpha}(\overline{\Omega})$. Since $u(0)=0$, we then get that

\bequa\label{ineq:app:1}
|u(x)|\leq |u(x)-u(0)|\leq C |x|^\alpha.
\eequa
We then get with (\ref{ineq:app:1}) that

$$\left|\fe(x)\right|=\left|\frac{|u(x)|^{\crit-1-\eps}u}{|x|^s}-au\right|\leq
\frac{C}{|x|^{s-(\crit-1-\eps)\alpha}}+C$$
for all $x\in\Omega$. We distinguish 2 cases:

\smallskip\noindent{\it Case \ref{sec:app}.3.1:}
$s-(\crit-1-\eps)\alpha_0\leq 0$. In this case, for any $p>1$, up to
taking $\alpha$ close enough to $\alpha_0$, we get that

$$\fe\in L^p(\Omega).$$
Since $\Delta u=\fe$ and $u\in\huno$, it follows from standard elliptic
theory that for any $\theta\in (0,1)$, we have that $u\in
C^{1,\theta}(\overline{\Omega})$. It follows that $\alpha_0=1$. This
proves the claim in Case \ref{sec:app}.3.1.

\smallskip\noindent{\it Case \ref{sec:app}.3.2:}
$s-(\crit-1-\eps)\alpha_0>0$. In this case, for any
$p<\frac{n}{s-(\crit-1-\eps)\alpha_0}$, up to taking $\alpha$ close enough
to $\alpha_0$, we get that

$$\fe\in L^p(\Omega).$$
We distinguish 3 subcases.

\smallskip\noindent{\it Case \ref{sec:app}.3.2.1:}
$s-(\crit-1-\eps)\alpha_0<1$. In this case, up to taking $\alpha$ close
enough to $\alpha_0$, there exists $p>n$ such that

$$\fe\in L^p(\Omega).$$
Since $\Delta u=\fe$ and $u\in\huno$, it follows from standard elliptic
theory that there exist exists $\theta\in (0,1)$ such that $u\in
C^{1,\theta}(\overline{\Omega})$. It follows that $\alpha_0=1$. This
proves the claim in Case \ref{sec:app}.3.2.1.

\smallskip\noindent{\it Case \ref{sec:app}.3.2.2:}
$s-(\crit-1-\eps)\alpha_0=1$. In this case, for any $p<n$, up to taking
$\alpha$ close enough to $\alpha_0$, we get that

$$\fe\in L^p(\Omega).$$
Since $\Delta u=\fe$ and $u\in\huno$, it follows from standard elliptic
theory that $u\in C^{0,\tilde{\alpha}}(\overline{\Omega})$ for all
$\tilde{\alpha}\in (0,1)$. It follows that $\alpha_0=1$. This proves the
claim in Case \ref{sec:app}.3.2.2.

\smallskip\noindent{\it Case \ref{sec:app}.3.2.3:}
$s-(\crit-1-\eps)\alpha_0>1$. In this case, it follows from standard
elliptic theory that $u\in C^{0,\tilde{\alpha}}(\overline{\Omega})$ for
all

$$\tilde{\alpha}\leq 2-(s-(\crit-1-\eps)\alpha_0).$$
It follows from the definition of $\alpha_0$ that

$$\alpha_0\geq  2-(s-(\crit-1-\eps)\alpha_0),$$
and then

$$0\geq 2-s+\left(\crit-2-\eps\right)\alpha_0>0,$$
a contradiction since $s<2$ and $\eps<\crit-2$. This proves that Case
\ref{sec:app}.3.2.3 does not occur, and we are back to the other cases.

\smallskip\noindent Clearly, theses cases end Step \ref{sec:app}.3.

\medskip\noindent{\bf Step \ref{sec:app}.4:} We claim that

$$u\in C^{1,\theta}(\overline{\Omega})$$
for all $\theta\in (0,\min\{1,\crit-\eps-s\})$. We proceed as in Step
\ref{sec:app}.3. We let $\alpha\in (0,1)$ (note that $\alpha_0=1$). We
then get that

$$\left|\fe(x)\right|=\left|\frac{|u|^{\crit-1-\eps}u}{|x|^s}-au\right|\leq
\frac{C}{|x|^{s-(\crit-1-\eps)\alpha}}+C$$
for all $x\in\Omega$. We distinguish 2 cases:

\smallskip\noindent{\it Case \ref{sec:app}.4.1:} $s-(\crit-1-\eps)\leq 0$.
In this case, for any $p>1$, up to taking $\alpha$ close enough to
$\alpha_0$, we get that

$$\fe\in L^p(\Omega).$$
Since $\Delta u=\fe$ and $u\in\huno$, it follows from standard elliptic
theory that $u\in C^{1,\theta}(\overline{\Omega})$ for all $\theta\in
(0,1)$. It follows that $\alpha_0=1$. This proves the claim in Case
\ref{sec:app}.4.1.

\smallskip\noindent{\it Case \ref{sec:app}.4.2:} $s-(\crit-1-\eps)>0$. In
this case, for any $p<\frac{n}{s-(\crit-1-\eps)}$, up to taking $\alpha$
close enough to $1$, we get that

$$\fe\in L^p(\Omega).$$
As easily checked,

$$1-(s-(\crit-1-\eps))=2-s+(\crit-1-\eps)-1>\crit-2-\eps> 0$$
We then get that there exists $p>n$ such that $\fe\in L^p(\Omega)$. Since
$\Delta u=\fe$ and $u\in\huno$, it follows from standard elliptic
theory that $u\in C^{1,\theta}(\overline{\Omega})$ for all $\theta\in
(0,\min\{1,\crit-\eps-s\})$. This proves the claim in Case
\ref{sec:app}.4.2.

\medskip\noindent Combining Case \ref{sec:app}.4.1 and Case
\ref{sec:app}.4.2, we obtain Step \ref{sec:app}.4. Proposition
\ref{prop:app} follows from Step \ref{sec:app}.4.
\end{proof}

\section{Appendix B: Properties of the Green's function}\label{sec:green}

This section is devoted to the proof of some useful properties of the
Green's function for a coercive operator. Concerning notations, for any
function $F: X\times Y\to \rr$ and any $x\in X$, we let $F_x:Y\to \rr$
such that $F_x(y)=F(x,y)$ for all $y\in Y$. We prove the following:

\begin{thm}\label{th:green:1}
Let $\Omega$ be a bounded domain of $\rn$, $n\geq 3$. Let $K,\lambda>0$. 
Let $\theta\in (0,1)$ and $a\in C^{0,\theta}(\overline{\Omega})$ such that

\bequa\label{co:green}
|a(x)|\leq K\hbox{ and }|a(x)-a(y)|\leq K |x-y|^\theta
\eequa
for all $x,y\in\overline{\Omega}$ and

\bequa\label{coerc:green}
\int_\Omega(|\nabla\varphi|^2+a\varphi^2)\, dx\geq
\lambda\int_\Omega\varphi^2\, dx
\eequa
for all $\varphi\in C^\infty_c(\Omega)$. Then there exists
$G:\overline{\Omega}\times\overline{\Omega}\setminus\{(x,x)/\,
x\in\overline{\Omega}\}\to \rr$ such that

\smallskip\noindent{\bf (G1)} For any $x\in\Omega$, $G_x\in L^1(\Omega)$
and $G_x\in C^{2,\theta}(\overline{\Omega}\setminus\{x\})$.

\smallskip\noindent{\bf (G2)} For any $x\in\Omega$, $G_x>0$ in
$\Omega\setminus\{x\}$ and $G_x=0$ on $\partial\Omega$.

\smallskip\noindent{\bf (G3)} For any $\varphi\in C^2(\overline{\Omega})$
such that $\varphi\equiv 0$ on $\partial\Omega$, we have that

$$\varphi(x)=\int_\Omega G(x,y)(\Delta\varphi+a\varphi)(y)\, dy$$
for all $x\in\Omega$.\par

\smallskip\noindent{\bf (G4)} $G(x,y)=G(y,x)$ for all $x,y\in\Omega$,
$x\neq y$.\par

\smallskip\noindent{\bf (G5)} There exists $C=C(\Omega,K,\lambda)>0$ such
that

$$|x-y|^{n-2}|G(x,y)|\leq C(\Omega,K,\lambda)$$
for all $x,y\in\Omega$, $x\neq y$.

\smallskip\noindent{\bf (G6)} There exists $C=C(\Omega,K,\lambda)>0$ such
that

$$|x-y|^{n-1}|G(x,y)|\leq C(\Omega,K,\lambda) d(y,\partial\Omega)$$
for all $x,y\in\Omega$, $x\neq y$.

\smallskip\noindent{\bf (G7)} There exists $C=C(\Omega,K,\lambda)>0$ such
that

$$|x-y|^{n-1}|\nabla G_x(y)|\leq C(\Omega,K,\lambda)$$
for all $x,y\in\Omega$, $x\neq y$.

\smallskip\noindent{\bf (G8)} There exists $C=C(\Omega,K,\lambda)>0$ such
that

$$|x-y|^{n}|\nabla_y G_x(y)|\leq C(\Omega,K,\lambda) d(x,\partial\Omega)$$
for all $x,y\in\Omega$, $x\neq y$.
\end{thm}

\noindent Some similar properties are available for the normal derivative
of $G$ at the boundary. Namely,

\begin{thm}\label{th:green:2}
Let $\Omega$ be a bounded domain of $\rn$, $n\geq 3$. We assume
$0\in\partial\Omega$. Let $K,\lambda>0$. Let $\theta\in
(0,1)$ and $a\in C^{0,\theta}(\overline{\Omega})$ such that
(\ref{co:green}) and (\ref{coerc:green}) hold. We let $G$ as in Theorem
\ref{th:green:1}. We let $H(x)=-\partial_\nu G_x(0)$ for all
$x\in\overline{\Omega}\setminus\{x\}$. Then the following assertions hold:

\smallskip\noindent{\bf (G9)} $H\in C^2(\overline{\Omega}\setminus\{0\})$,
$H>0$ in $\Omega$ and $H\equiv 0$ on $\partial\Omega\setminus\{0\}$,\par

\smallskip\noindent{\bf (G10)} $\Delta H+aH=0$ in $\Omega$,\par

\smallskip\noindent{\bf (G11)} There exists $C=C(\Omega,K,\lambda)>0$ such
that

$$\frac{d(x,\partial\Omega)}{C|x|^n}\leq H(x)\leq \frac{C d(x,\partial\Omega)}{|x|^n}$$
for all $x\in\Omega$.\par

\smallskip\noindent{\bf (G12)} There exists $C=C(\Omega,K,\lambda)>0$  and
$\delta=\delta(\Omega,K,\lambda)>0$ such that
$$\frac{1}{C|x|^n}\leq |\nabla H(x)|\leq \frac{C}{|x|^n}$$
for all $x\in B_\delta(0)\cap \Omega$.
\end{thm}

The proof of Theorem \ref{th:green:1} is very close to the proof of the
existence of the Green's function on a compact manifold without boundary
provided in \cite{dhr}. We just give the main steps of the proof and
outline the difference with \cite{dhr} when necessary. We prove Theorem
\ref{th:green:2} in details.

\medskip\noindent{\bf Step \ref{sec:green}.1:} This Step is devoted to the
proof of points (G1)-(G5) of Theorem \ref{th:green:1}. We only sketch the
proof. Details are available in \cite{dhr}. We define

$${\mathcal H}(x,y)=\frac{1}{(n-2)\omega_{n-1}|x-y|^{n-2}}$$
for all $x,y\in\rn$ such that $x\neq y$. In this expression,
$\omega_{n-1}$ denotes the volume of the standard $(n-1)-$sphere. The
function ${\mathcal H}$ is the standard Green kernel of the Laplacian in
$\rn$. We define the functions $\Gamma_i$'s by induction. Given
$x,y\in\overline{\Omega}$, $x\neq y$, we let

$$\begin{array}{ll}
\Gamma_1(x,y)=-a(y){\mathcal H}(x,y)&\\
\Gamma_{i+1}(x,y)=\int_\Omega \Gamma_i(x,z)\Gamma_1(z,y)\, dz &\hbox{ for
all }i\geq 1.
\end{array}$$
As easily checked, $\Gamma_i\in
C^0(\overline{\Omega}\times\overline{\Omega}\setminus\{(x,x)/\,
x\in\overline{\Omega}\})$ for all $i\geq 1$. Standard computations yield
that there exists $C(\Omega,n,K)>0$ such that

$$\begin{array}{ll}
|\Gamma_i(x,y)|\leq C(\Omega,n,K)|x-y|^{2i-n} & \hbox{ if }2i<n\\
|\Gamma_i(x,y)|\leq C(\Omega,n,K)\left(1+\left|\ln
|x-y|\right|\right) & \hbox{ if }2i=n\\
|\Gamma_i(x,y)|\leq C(\Omega,n,K) & \hbox{ if }2i>n,\; i\leq n.
\end{array}$$
for all $x,y\in\overline{\Omega}$, $x\neq y$. In addition, $\Gamma_i$ can
be extended to a continuous function in
$\overline{\Omega}\times\overline{\Omega}$ for all $i>n/2$. We let
$x\in\Omega$. We let $U_x \in H_{1,0}^2(\Omega)$ such that

$$\Delta U_x+aU_x=\Gamma_{n+1}(x,\cdot)\hbox{ in }{\mathcal D}'(\Omega).$$
Since $\Gamma_{n+1}$ is uniformly bounded in $L^\infty$, it follows from
standard elliptic theory that $U_x\in H_2^p(\Omega)$ for all $p>1$ and
that there exists $C(\Omega,K,\lambda)>0$ such that

$$\Vert U_x\Vert_{C^1(\overline{\Omega})}\leq C(\Omega,K,\lambda)$$
for all $x\in\Omega$. We let $V_x\in H_1^2(\Omega)$ such that
$$\left\{\begin{array}{ll}
\Delta V_x+aV_x=0 &\hbox{ in }{\mathcal D}'(\Omega)\\
V_x(y)=-{\mathcal H}(x,y)-\sum_{i=1}^n\int_\Omega \Gamma_i(x,z){\mathcal
H}(z,y)\, dz& \hbox{ for all }y\in\partial\Omega.
\end{array}\right.$$
It follows from standard elliptic theory that for any $x\in\Omega$,
$V_x\in C^1(\Omega)$. Moreover, it follows from the explicit expression of
${\mathcal H}$ and the $\Gamma_i$'s that there exists
$C(\Omega,K,\lambda)'>0$ such that $V_x(y)\leq C(\Omega,K,\lambda)'$ for 
all $x\in\Omega$ and all $y\in\partial\Omega$. Since $\Delta+a$ is 
coercive, it follows from the comparison principle that there exists 
$C(\Omega,K,\lambda)>0$ such that
$$V_x(y)\leq C(\Omega,K,\lambda)$$
for all $x\in\Omega$ and all $y\in\Omega$. We let

$$G_x(y):={\mathcal H}(x,y)+\sum_{i=1}^n\int_\Omega \Gamma_i(x,z){\mathcal
H}(z,y)\, dz+ U_x(y)+V_x(y)$$
for all $y\in\Omega$. It follows from the
construction of $G$ that there exists $C(\Omega,K,\lambda)>0$ such that
$$G(x,y)\leq C(\Omega,K,\lambda)\cdot |x-y|^{2-n}$$
for all $x,y\in\Omega$, $x\neq y$ and that $G_x$ vanishes on $\partial\Omega$ for all $x\in\Omega$. This prove point (G5). We let $\varphi\in C^2(\overline{\Omega})$ such that
$\varphi\equiv 0$ on $\partial\Omega$. Noting that
$$\varphi(z)=\int_\Omega {\mathcal H}(z,y)\Delta\varphi(y)\,
dy+\int_{\partial\Omega}{\mathcal H}(x,y)\partial_\nu\varphi(y)\,
d\sigma(y)$$
for all $z\in\Omega$, we get with some integrations by parts that
$$\varphi(x)=\int_\Omega G(x,y)(\Delta\varphi+a\varphi)(y)\, dy.$$
This proves point (G3). It then follows that
$$\Delta G_x+aG_x=0\hbox{ in }{\mathcal D}'(\Omega\setminus\{x\}).$$
Since $G_x\equiv 0$ on $\partial\Omega$, we get that $G_x\in
C^{2,\theta}_{loc}(\overline{\Omega}\setminus\{x\})$. It the follows from
the construction and the maximum principle that $G_x>0$ in
$\Omega\setminus\{x\}$. This proves points (G2) and (G1). Point (G4) is standard, we refer to \cite{aubin} or \cite{dhr}.

\medskip\noindent{\bf Step \ref{sec:green}.2:} We prove points (G6) and
(G7) of Theorem \ref{th:green:1}. We proceed by contradiction and assume
that there exists a sequence $(a_k)_{k\in\nn}\in
C^{0,\theta}(\overline{\Omega})$ and sequences
$(x_k)_{k\in\nn},(y_k)_{k\in\nn}\in\Omega$ such that (\ref{co:green}) and
(\ref{coerc:green}) hold and

\bequa\label{hyp:G:c1}
\lim_{k\to +\infty}\left[|x_k-y_k|^{n-1}|\nabla
G_{x_k}(y_k)|+\frac{|x_k-y_k|^{n-1}G_{x_k}(y_k)}{d(y_k,\partial\Omega)}\right]=+\infty
\eequa
where $G_{x_k}$ is the Green's function for $\Delta+a_k$ at $x_k$. We let
$x_\infty=\lim_{k\to +\infty}x_k$ and $y_\infty=\lim_{k\to +\infty}y_k$
(these limits exist up to a subsequence).

\medskip\noindent{\it Case 1:} $x_\infty\neq y_\infty$. We let
$0<\delta<|x_\infty-y_\infty|/4$. It follows from point (G5) that there
exists $C>0$ independant of $k$ such that $|G_{x_k}(y)|\leq C$ for all
$y\in\Omega\cap B_{y_\infty}(2\delta)$. Since $\Delta G_{x_k}+a_k
G_{x_k}=0$ and $G_{x_k}=0$ on $\partial\Omega$, it follows from standard
elliptic theory that

$$\Vert G_{x_k}\Vert_{C^1(\overline{\Omega}\cap
B_{y_\infty}(\delta))}=O(1)$$
when $k\to +\infty$. Since $G_{x_k}$ vanishes on $\partial\Omega$, we get
that there exists $C>0$ such that

$$|G_{x_k}(y)|\leq C d(y,\partial\Omega)\hbox{ and }|\nabla
G_{x_k}(y)|\leq C$$
for all $y\in \overline{\Omega}\cap B_{y_\infty}(\delta)$ and all
$\eps>0$. A contradiction with (\ref{hyp:G:c1}).

\medskip\noindent{\it Case 2:} $x_\infty=y_\infty$.\par

\smallskip\noindent{\it Case 2.1:} We assume that

\bequa\label{ineq:hyp:case:1}
d(x_k,\partial\Omega)\geq 2|y_k-x_k|
\eequa
up to a subsequence. We let

$$\tilde{G}_k(z)=|y_k-x_k|^{n-2}G(x_k, x_k+|y_k-x_k|z)$$
for all $z\in B_{3/2}(0)$. With our assumption, this is well defined. It
follows from (G5) that there exists $C>0$ such that

$$|\tilde{G}_k(z)|\leq C$$
for all $z\in B_{3/2}(0)\setminus \overline{B}_{1/4}(0)$. Moreover,
$\tilde{G}_k$ verifies the equation

$$\Delta\tilde{G}_k+|y_k-x_k|^2a_k(x_k+|y_k-x_k|z)\tilde{G}_k(z)=0 $$
in $B_{3/2}(0)\setminus \overline{B}_{1/4}(0)$. It follows from standard
elliptic theory that

$$\Vert\tilde{G}_k\Vert_{C^1(B_{5/4}(0)\setminus
\overline{B}_{1/2}(0))}=O(1)$$
when $k\to +\infty$. Taking $z=\frac{y_k-x_k}{|y_k-x_k|}$ and coming back to $G_{x_k}$, we get that
\bequa\label{ineq:A}
|x_k-y_k|^{n-1}|\nabla G_{x_k}(y_k)|=O(1)
\eequa
when $k\to +\infty$. Moreover, it follows from point (G5) of Theorem 
\ref{th:green:1} and (\ref{ineq:hyp:case:1}) that there exists $C>0$ such 
that
\bequa\label{ineq:B}
|x_k-y_k|^{n-1} G_{x_k}(y_k)\leq C d(y_k,\partial\Omega)
\eequa
when $k\to +\infty$. Inequations (\ref{ineq:A}) and (\ref{ineq:B}) 
contradict (\ref{hyp:G:c1}).

\smallskip\noindent{\it Case 2.2:} We assume that

\bequa\label{hyp:22}
d(x_k,\partial\Omega)\leq 2|y_k-x_k|
\eequa
up to a subsequence. In particular, $x_\infty\in\partial\Omega$. We let a
chart $\varphi: U\to V$ as in (\ref{def:vphi}) with $y_0=x_\infty$ and
where $U,V$ are open neighborhoods of $0$ and $x_\infty$ respectively. We
let $\tilde{x}_k,\tilde{y}_k\in U\cap \{x_1<0\}$ such that
$x_k=\varphi(\tilde{x}_k)$ and $y_k=\varphi(\tilde{y}_k)$. As a remark,
$\lim_{k\to +\infty}\tilde{x}_k=\lim_{k\to +\infty}\tilde{y}_k=0$. We let
$\tilde{x}_{k,1}<0$ be the first coordinate of $\tilde{x}_k$. As in Step
\ref{sec:exh}.2, we have that
$d(x_k,\partial\Omega)=(1+o(1))|\tilde{x}_{k,1}|$ when $k\to +\infty$. We
then get with (\ref{hyp:22}) that
$\tilde{x}_{k,1}=O(|\tilde{y}_k-\tilde{x}_k|)$ when $k\to +\infty$. We let

$$\rho_k=\frac{\tilde{x}_{k,1}}{|\tilde{y}_k-\tilde{x}_k|}\hbox{ and
}\rho_\infty=\lim_{k\to+\infty}\rho_k$$
(this limit exists up to a subsequence). We let $R>0$ and we let
$$\tilde{G}_k(z)=|\tilde{y}_k-\tilde{x}_k|^{n-2}G(x_k,
\varphi\left(\tilde{x}_k+|\tilde{y}_k-\tilde{x}_k|\left(z-\rho_k\vec{e}_1\right)\right))$$
for all $k$ and all $z\in B_R(0)\cap \{z_1\leq 0\}$. Here $\vec{e}_1$
denotes the first vector of the canonical basis of $\rn$. Note that
$\tilde{G}_k$ vanishes on $B_R(0)\cap \{z_1= 0\}$. It follows of the
pointwise estimate (G5) that for any $R,\delta>0$, there exists
$C(R,\delta)>0$ such that

$$|\tilde{G}_k(z)|\leq C(R,\delta)$$
for all $z\in [B_R(0)\setminus
\overline{B}_\delta((\rho_\infty,0,...,0))]\cap \{z_1\leq 0\}$. The
function $\tilde{G}_k$ verifies the equation

$$\Delta_{g_k}\tilde{G}_k+|\tilde{y}_k-\tilde{x}_k|^2a_k(\varphi\left(\tilde{x}_k+|\tilde{y}_k-\tilde{x}_k|\left(z-(\rho_k,0,..,0)\right)\right))\tilde{G}_k=0$$
in $[B_R(0)\setminus \overline{B}_\delta((\rho_\infty,0,...,0))]\cap
\{z_1\leq 0\}$. It then follows from standard elliptic theory that
$\Vert\tilde{G}_k\Vert_{C^1([B_{R/2}(0)\setminus
\overline{B}_{2\delta}(\rho_\infty,...,0)]\cap \{z_1\leq 0\})}=O(1)$ when 
$k\to +\infty$. As
in Case 2.1, we get that
\bequa\label{ineq:A:bis}
|x_k-y_k|^{n-1}|\nabla G_{x_k}(y_k)|=O(1)
\eequa
when $k\to +\infty$. Moreover, since $\tilde{G}_k$ vanishes on 
$\partial\rnm$, there exists $C>0$ such that
$$|\tilde{G}_k(z)|\leq C|z_1|$$
for all $z\in [B_{R/2}(0)\setminus
\overline{B}_{2\delta}(\rho_\infty,...,0)]\cap \{z_1\leq 0\}$. Taking 
$z=(\rho_k,...,0)+\frac{\tilde{y}_k-\tilde{x}_k}{|\tilde{y}_k-\tilde{x}_k|}$, 
we get that
\bequa\label{ineq:B:bis}
|x_k-y_k|^{n-1}G_{x_k}(y_k)\leq C d(y_k,\partial\Omega)
\eequa
for all $k$ large enough. A contradiction with (\ref{hyp:G:c1}).

\smallskip\noindent In all the cases, we have contradicted
(\ref{hyp:G:c1}). This proves points (G6) and (G7) of Theorem
\ref{th:green:1}.

\medskip\noindent{\bf Step \ref{sec:green}.3:} We prove point (G8) of
Theorem \ref{th:green:1}. More precisely, we claim that there exists
$C=C(\Omega,K,\lambda)>0$ such that

\bequa\label{ineq:proof:H}
|x-y|^{n}G(x,y)\leq C d(y,\partial\Omega)d(x,\partial\Omega)
\eequa
and
$$|x-y|^n|\nabla_y G(x,y)|\leq C d(x,\partial\Omega)$$
for all $x,y\in\Omega$, $x\neq \Omega$. Indeed we proceed as in the proof
of points (G6) and (G7). We proceed by contradiction and assume that there
exist a sequence $(a_k)_{k\in\nn}\in C^{0,\theta}(\overline{\Omega})$ and
sequences $(x_k)_{k\in\nn},(y_k)_{k\in\nn}\in\Omega$ such that
(\ref{co:green}) and (\ref{coerc:green}) hold and

\bequa\label{hyp:G:co}
\lim_{k\to
+\infty}|x_k-y_k|^{n}\frac{|G(x_k,y_k)|}{d(x_k,\partial\Omega)
d(y_k,\partial\Omega)}+|x_k-y_k|^{n}\frac{|\nabla
G_{x_k}(y_k)|}{d(x_k,\partial\Omega)}=+\infty
\eequa
where $G_{x_k}$ is the Green's function for $\Delta+a_k$ at $x_k$. We let
$x_\infty=\lim_{k\to +\infty}x_k$ and $y_\infty=\lim_{k\to +\infty}y_k$
(these limits exist up to a subsequence).

\medskip\noindent{\it Case 1:} $x_\infty\neq y_\infty$. We let
$0<\delta<|x_\infty-y_\infty|/4$. We let

$$\tilde{G}_k(z)=\frac{G_k(x_k,z)}{d(x_k,\partial\Omega)}$$
for all $z\in\Omega$. As in Case 1 of the proof of (G6)-(G7), using (G6),
we get that

$$\Vert \tilde{G}_k\Vert_{C^1(\overline{\Omega}\cap
B_{y_\infty}(\delta))}=O(1)$$
when $k\to +\infty$. It then follows that

$$\tilde{G}_k(y_k)\leq C d(y_k,\partial\Omega)\hbox{ and
}|\nabla\tilde{G}_k(y_k)|\leq C$$
when $k\to +\infty$. A contradiction with (\ref{hyp:G:co}).

\medskip\noindent{\it Case 2:} $x_\infty=y_\infty$.\par

\smallskip\noindent{\it Case 2.1:} We assume that \

$$d(x_k,\partial\Omega)\geq 2|y_k-x_k|$$
up to a subsequence. We then obtain that $|x_k-y_k|\leq
d(y_k,\partial\Omega)$. This inequality and (G6)-(G7) yield to a
contradiction with (\ref{hyp:G:co}).

\smallskip\noindent{\it Case 2.2:} We assume that

$$d(x_k,\partial\Omega)\leq 2|y_k-x_k|$$
up to a subsequence. In particular, $x_\infty\in\partial\Omega$. We let a
chart $\varphi: U\to V$ as in (\ref{def:vphi}) with $y_0=x_\infty$ and
where $U,V$ are open neighborhoods of $0$ and $x_\infty$ respectively. We
let $\tilde{x}_k,\tilde{y}_k\in U\cap \{x_1<0\}$ such that
$x_k=\varphi(\tilde{x}_k)$ and $y_k=\varphi(\tilde{y}_k)$. We let

$$\tilde{G}_k(z)=|\tilde{y}_k-\tilde{x}_k|^{n-1}\frac{G\left[x_k,
\varphi\left(\tilde{x}_k+|\tilde{y}_k-\tilde{x}_k|\left(z-(\frac{\tilde{x}_{k,1}}{|\tilde{y}_k-\tilde{x}_k|},0,..,0)\right)\right)\right]}{d(x_k,\partial\Omega)}$$
for all $z\in [B_{R}(0)\setminus
\overline{B}_{\delta}(\rho_\infty,0,...,0)]\cap \{z_1\leq 0\}$.
As in Case 2.2 of the proof of (G6)-(G7), we get with (G6) that for any
$R>4\delta>0$, we have that

$$\Vert\tilde{G}_k\Vert_{C^1([B_{R/2}(0)\setminus
\overline{B}_{2\delta}(\rho_\infty,0,...,0)]\cap \{z_1\leq 0\})}=O(1)$$
when $k\to +\infty$, where $\rho_\infty=\lim_{k\to +\infty} 
\frac{\tilde{x}_{k,1}}{|\tilde{y}_k-\tilde{x}_k|}$. Since $\tilde{G}_k$ 
vanishes on $\{z_1=0\}$, it then
follows that there exists $C>0$ such that $|\tilde{G}_k(z)|\leq C|z_1|$
for all $z\in [B_{R/2}(0)\setminus 
\overline{B}_{2\delta}(\rho_\infty,0,...,0)]\cap
\{z_1\leq 0\}$. Coming back to the definition of $\tilde{G}_k$ and noting
that $d(y_k,\partial\Omega)=(1+o(1))|\tilde{y}_{k,1}|$ when $k\to
+\infty$, we get a contradiction with (\ref{hyp:G:co}) as in Case 2.2 of
Step \ref{sec:green}.2.

\smallskip\noindent In all the cases, we have contradicted
(\ref{hyp:G:co}). This proves the claim and ends Step \ref{sec:green}.3.

\medskip\noindent The proof of Theorem \ref{th:green:1} is complete. We
prove Theorem \ref{th:green:2}.

\medskip\noindent{\bf Step \ref{sec:green}.4:} We let $H(x)=-\partial_\nu
G_x(0)$ for any $x\in\overline{\Omega}\setminus\{0\}$. It follows from
(\ref{ineq:proof:H}) that  there exists $C=C(\Omega,K,\lambda)>0$ such
that

\bequa\label{upper:bnd:H}
0\leq H(x)\leq \frac{C d(x,\partial\Omega)}{|x|^{n}}\leq
\frac{C}{|x|^{n-1}}
\eequa
for all $x\in\Omega$. Since $\Delta G_x+aG_x=0$ in $\Omega\setminus\{x\}$,
using the symetry (G4) of $G$ and (\ref{upper:bnd:H}), we get that $H\in
C^2(\overline{\Omega}\setminus\{0\})$ and that $\Delta H+aH=0$ in $\Omega$
and $H(x)=0$ for all $x\in\partial\Omega\setminus\{0\}$. Derivating (G3), 
we get that

\bequa\label{eq:dist:H}
\partial_\nu\varphi(0)=-\int_\Omega H(x)(\Delta\varphi+a\varphi)(x)\, dx
\eequa
for all $\varphi\in C^2(\overline{\Omega})$ such that $\varphi\equiv 0$ on
$\partial\Omega$.

\medskip\noindent{\bf Step \ref{sec:green}.5:} Assume that there exists a
sequence $(a_k)_{k>0}\in C^{0,\theta}(\overline{\Omega})$ such that
(\ref{coerc:green}) and (\ref{co:green}) hold, that there exists a
sequence $(r_k)_{k>0}\in\rr$ such that $r_k>0$, $\lim_{k\to +\infty}r_k=0$
and

$$\lim_{k\to +\infty}
\sup_{|x|=r_k}\frac{H_k(x)|x|^n}{d(x,\partial\Omega)}=0,$$
where $H_k$ comes from the Green's function of $\Delta+a_k$. We claim that in this situation, we have that
\bequa\label{lim:global}
\lim_{k\to +\infty}\sup_{\frac{1}{2}r_k\leq |x|\leq
3r_k}\left(\frac{H_k(x)|x|^n}{d(x,\partial\Omega)}+ |x|^n|\nabla
H_k(x)|\right)=0.
\eequa
Indeed, we let $\varphi:U\to V$ as in (\ref{def:vphi}) where $U,V$ are
open neighborhoods of $0$. We let

$$\tilde{H}_k(x)=r_k^{n-1}H_k(\varphi(r_k x))$$
for all $x\in \frac{U}{r_k}\cap\{x_1\leq 0\}$. It follows from
(\ref{upper:bnd:H}) that for any $R>\delta>0$, there exists
$C(R,\delta)>0$ such that $|\tilde{H}_k(x)|\leq C(R,\delta)$ for all $x\in
[B_R(0)\setminus B_\delta(0)]\cap \{x_1\leq 0\}$. In addition $\tilde{H}_k$ vanishes when $x_1=0$. Moreover, we have that

$$\Delta_{g_k}\tilde{H}_k+r_k^2 a_k(\varphi(r_k x))\tilde{H}_k=0,$$
where $(g_k)_{ij}=(\partial\varphi,\partial_j\varphi)(r_k x)$ for
$i,j\in\{1,...,n\}$. It then follows from standard elliptic theory that
there exists $\tilde{H}\in C^2(\overline{\rnm}\setminus\{0\})$ such that
$\Delta\tilde{H}=0$ in $\rnm\setminus\{0\}$ and

$$\lim_{k\to +\infty}\tilde{H}_k=\tilde{H}$$
in $C^2_{loc}(\overline{\rnm}\setminus\{0\})$. As easily checked, we have
that

\beqn
&&\lim_{k\to +\infty}\sup_{\frac{1}{2}r_k\leq |x|\leq
3r_k}\left(\frac{H_k(x)|x|^n}{d(x,\partial\Omega)}+ |x|^n|\nabla
H_k(x)|\right)\nonumber\\
&&=\sup_{\frac{1}{2}\leq |x|\leq 3}\left(\frac{\tilde{H}(x)|x|^n}{|x_1|}+
|x|^n|\nabla \tilde{H}(x)|\right)\label{lim:global:bis}
\eeqn
and
\bequa\label{eq:lim:inf}
0=\lim_{k\to +\infty}\sup_{
|x|=r_k}\frac{H_k(x)|x|^n}{d(x,\partial\Omega)}=\sup_{|x|=1}\left(\frac{\tilde{H}(x)|x|^n}{|x_1|}\right).
\eequa
Assume that $\tilde{H}\not\equiv 0$. Then, since $\tilde{H}\geq 0$
vanishes on $\partial\rnm$, we have that $\tilde{H}>0$ in $\rnm$ and
$\partial_1\tilde{H}<0$ on $\partial\rnm\setminus\{0\}$. It then follows
that the RHS of (\ref{eq:lim:inf}) is positive. A contradiction, since the
LHS is $0$. Then $\tilde{H}\equiv 0$, and (\ref{lim:global}) follows from
(\ref{lim:global:bis}). This ends Step \ref{sec:green}.5.

\medskip\noindent{\bf Step \ref{sec:green}.6:} We claim that there exists
$\eps(\Omega,K,\lambda)>0$ such that

\bequa\label{liminf:H}
\liminf_{r\to 0}\sup_{|x|=r}\frac{H(x)|x|^n}{d(x,\partial\Omega)}\geq
\eps(\Omega,K,\lambda).
\eequa
Indeed, we argue by contradiction and assume that there exists a sequence
$(a_k)_{k>0}\in C^{0,\theta}(\overline{\Omega})$ such that
(\ref{coerc:green}) and (\ref{co:green}) hold, that there exists a
sequence $(r_k)_{k>0}\in\rr$ such that $r_k>0$, $\lim_{k\to +\infty}r_k=0$
and

$$\lim_{k\to +\infty}
\sup_{|x|=r_k}\frac{H_k(x)|x|^n}{d(x,\partial\Omega)}=0,$$
where $H_k$ comes from the Green's function of $\Delta+a_k$. It then
follows from Step \ref{sec:green}.5. that

\bequa\label{lim:global:proof}
\lim_{k\to +\infty}m_k=0.
\eequa
where $$m_k:=\sup_{\frac{1}{2}r_k\leq |x|\leq
3r_k}\left(\frac{H_k(x)|x|^n}{d(x,\partial\Omega)}+ |x|^n|\nabla
H_k(x)|\right).$$
We let $\tilde{\eta}\in C^\infty(\rn)$ such $\tilde{\eta}\equiv 0$ in
$B_1(0)$ and $\tilde{\eta}\equiv 1$ in $\rn\setminus B_2(0)$. We let
$\eta_k(x)=\tilde{\eta}(x/r_k)$ for all $x\in\rn$ and all $k>0$. We let
$\varphi_k\in C^2(\overline{\Omega})$ such that

$$\Delta\varphi_k+a_k\varphi_k=1\hbox{ in }\Omega\hbox{ and
}\varphi_k\equiv 0\hbox{ on }\partial\Omega.$$
It follows from standard elliptic theory that $\lim_{k\to
+\infty}\varphi_k=\varphi\not\equiv 0$ in $C^2(\overline{\Omega})$. It then follows
from Hopf's maximum principle that

\bequa\label{der:hopf}
\partial_\nu\varphi(0)<0.
\eequa
Integrating by parts and using that $\Delta H_k+a_kH_k=0$, we obtain that

\beq
\int_\Omega H_k(x)(\Delta\varphi_k+a_k\varphi_k)(x)\, dx&=& \int_\Omega
(\eta_kH_k)(x)(\Delta\varphi_k+a_k\varphi_k)(x)\, dx+o(1)\\
&=& \int_\Omega (\Delta(\eta_kH_k)+a\eta_kH_k)\varphi_k\, dx+o(1)\\
&=& \int_\Omega ((\Delta\eta_k)H_k-2\nabla\eta_k\nabla H_k)\varphi_k\,
dx+o(1)\\
&=& \int_{\Omega\cap B_{2r_k}(0)\setminus B_{r_k}(0)}
((\Delta\eta_k)H_k-2\nabla\eta_k\nabla H_k)\varphi_k\, dx\\
&&+o(1)
\eeq
where $\lim_{k\to +\infty}o(1)=0$. Since $\varphi_k(0)=0$ and $\lim_{k\to
+\infty}\varphi_k=\varphi$ in $C^1(\overline{\Omega})$, using the
definition of $m_k$ we get that

\beq
\int_\Omega H_k(x)(\Delta\varphi_k+a_k\varphi_k)(x)\,
dx&=&O\left(r_k^n(m_k r_k^{-2}r_k^{1-n}r_k)\right)+o(1)=O(m_k)+o(1).
\eeq
With (\ref{lim:global:proof}), letting $k\to +\infty$, and using
(\ref{eq:dist:H}) we get that

$$\partial_\nu\varphi(0)=\partial_\nu\varphi_k(0)+o(1)=-\int_\Omega
H_k(x)(\Delta\varphi_k+a_k\varphi_k)(x)\, dx+o(1)=0.$$
A contradiction with (\ref{der:hopf}), and the claim is proved.

\medskip\noindent{\bf Step \ref{sec:green}.7:} We claim that there exists
$\eps(\Omega,K,\lambda)>0$ such that

\bequa\label{liminf:infH}
\liminf_{r\to 0}\inf_{|x|=r}\frac{H(x)|x|^n}{d(x,\partial\Omega)}\geq
\eps(\Omega,K,\lambda).
\eequa
Indeed, we argue by contradiction and assume that there exists a sequence
$(a_k)_{k>0}\in C^{0,\theta}(\overline{\Omega})$ such that
(\ref{coerc:green}) and (\ref{co:green}) hold, that there exists a
sequence $(r_k)_{k>0}\in\rr$ such that $r_k>0$, $\lim_{k\to +\infty}r_k=0$
and

$$\lim_{k\to +\infty}
\inf_{|x|=r_k}\frac{H_k(x)|x|^n}{d(x,\partial\Omega)}=0,$$
where $H_k$ comes from the Green's function of $\Delta+a_k$. Mimicking the
proof of Step \ref{sec:green}.5, we obtain that
$\tilde{H}_k(x):=r_k^{n-1}H_k(\varphi(r_kx))$ converges to $\tilde{H}$ in
$C^1_{loc}(\overline{\rnm}\setminus\{0\})$. We get that

$$\inf_{|x|=1}\frac{\tilde{H}(x)|x|^n}{|x_1|}=\lim_{k\to +\infty}
\inf_{|x|=r_k}\frac{H_k(x)|x|^n}{d(x,\partial\Omega)}=0.$$
Since $\tilde{H}\geq 0$ is harmonic and vanishes on
$\partial\rnm\setminus\{0\}$, it follows from Hopf's maximum principle
that $\tilde{H}\equiv 0$. We then get that
$$\lim_{k\to +\infty}
\sup_{|x|=r_k}\frac{H_k(x)|x|^n}{d(x,\partial\Omega)}=\sup_{|x|=1}\frac{\tilde{H}(x)|x|^n}{|x_1|}=0.$$
A contradiction with Step \ref{sec:green}.6. This proves the claim.

\medskip\noindent{\bf Step \ref{sec:green}.8:} We claim that there exists
$C=C(\Omega,K,\lambda)>0$ such that

$$\frac{d(x,\partial\Omega)}{C|x|^n}\leq H(x)\leq
\frac{C d(x,\partial\Omega)}{|x|^n}$$
for all $x\in\Omega\setminus\{0\}$. Indeed, this claim is a consequence of
(\ref{upper:bnd:H}),  Step \ref{sec:green}.7 and standard elliptic theory.
This proves point (G11).

\medskip\noindent{\bf Step \ref{sec:green}.9:} We claim that there exists
$C(\Omega,K,\lambda)>0$ such that

\bequa\label{claim:step9}
|x|^n|\nabla H(x)|\leq C(\Omega,K,\lambda)
\eequa
for all $x\in\Omega\setminus\{0\}$. We proceed by contradiction and assume
that that there exists a sequence $(a_k)_{k>0}\in
C^{0,\theta}(\overline{\Omega})$ such that (\ref{coerc:green}) and
(\ref{co:green}) hold, that there exists a sequence $(x_k)_{k>0}\in\Omega$
such that

\bequa\label{step9}
\lim_{k\to+\infty}|x_k|^n|\nabla H_k(x_k)|=+\infty,
\eequa
where $H_k$ comes from the Green's function of $\Delta+a_k$.

\smallskip\noindent{\it Case 1:} $\lim_{k\to +\infty}x_k\neq 0$. In this
case, since $\Delta H_k+a_kH_k=0$, it follows from (\ref{upper:bnd:H}) and
standard elliptic theory that $|\nabla H_k(x_k)|=O(1)$ when $k\to
+\infty$.

\smallskip\noindent{\it Case 2:} $\lim_{k\to +\infty}x_k=0$. We consider
$\varphi:U\to V$ as in (\ref{def:vphi}) with $y_0=0$ and $U,V$ are open
neighborhoods of $0$. We let $x_k=\varphi(\tilde{x}_k)$. We let

$$\tilde{H}_k(x)=|\tilde{x}_k|^{n-1}H_k(\varphi(|\tilde{x}_k|x))$$
for all $x\in \frac{U}{|\tilde{x}_k|}\cap\{x_1\leq 0\}$. As in Step
\ref{sec:green}.5, we get that there exists $C>0$ such that

$$\Vert\tilde{H}_k\Vert_{C^1(\{x_1\leq 0\}\cap B_2(0)\setminus
B_{1/2}(0))}\leq C.$$
Estimating the gradient at $\tilde{x}_k/|\tilde{x}_k|$, we get that
$$|x_k|^n|\nabla H_k(x_k)|=O(1)$$
when $k\to +\infty$.

\smallskip\noindent In both cases, we have contradicted (\ref{step9}).
This proves (\ref{claim:step9}).

\medskip\noindent{\bf Step \ref{sec:green}.10:} We claim that there exists
$\delta(\Omega,K,\lambda), C(\Omega,K,\lambda)>0$ such that

\bequa\label{claim:step10}
|x|^n|\nabla H(x)|\geq C(\Omega,K,\lambda)
\eequa
for all $x\in\Omega\setminus\{0\}$ such that $|x|\leq
\delta(\Omega,K,\lambda)$. We proceed by contradiction and assume that
that there exists a sequence $(a_k)_{k>0}\in
C^{0,\theta}(\overline{\Omega})$ such that (\ref{coerc:green}) and
(\ref{co:green}) hold, that there exists a sequence $(x_k)_{k>0}\in\Omega$
such that $\lim_{k\to +\infty}x_k=0$ and

\bequa\label{step10}
\lim_{k\to+\infty}|x_k|^n|\nabla H_k(x_k)|=0,
\eequa
where $H_k$ comes from the Green's function of $\Delta+a_k$. We let
$x_k=\varphi(\tilde{x}_k)$ and $y_k=\varphi(\tilde{y}_k)$. We let

$$\tilde{H}_k(x)=|\tilde{x}_k|^{n-1}H_k(\varphi(|\tilde{x}_k|x))$$
for all $x\in \frac{U}{|\tilde{x}_k|}\cap\{x_1\leq 0\}$. Mimicking the
proof of Steps \ref{sec:green}.5 and \ref{sec:green}.9, we get that there
exists $\tilde{H}\in C^2(\overline{\rnm}\setminus\{0\})$ such that

\bequa\label{lim:Hk}
\lim_{k\to +\infty}\tilde{H}_k=\tilde{H}
\eequa
in $C^2_{loc}(\overline{\rnm}\setminus\{0\})$. In particular, we have that
$\tilde{H}$ is harmonic. It follows from Step \ref{sec:green}.8  and
(\ref{lim:Hk}) that there exists $C>0$ such that

$$\frac{|x_1|}{C|x|^n}\leq \tilde{H}(x)\leq
\frac{C|x_1|}{|x|^n}$$
for all $x\in\rnm\setminus\{0\}$. It then follows from the rigidity
Property \ref{prop:rigidity} below that $\nabla\tilde{H}(x)\neq 0$ for all
$x\in\overline{\rnm}\setminus\{0\}$. It follows from (\ref{step10}) and
(\ref{lim:Hk}) that there exists $\hat{x}\in\overline{\rnm}\setminus\{0\}$
such that $\nabla\tilde{H}(\hat{x})=0$. A contradiction. This proves
(\ref{claim:step10}).

\medskip\noindent Clearly Theorem \ref{th:green:2} is a consequence of
Steps \ref{sec:green}.4 to \ref{sec:green}.10.

\medskip\noindent{\bf Step \ref{sec:green}.11:} Our last step is the proof
of the following rigidity result:

\begin{prop}\label{prop:rigidity}
Let $h\in C^2(\overline{\rnm}\setminus\{0\})$. We assume that $h$ is
nonnegative in a neighborhood of $0$, harmonic and vanishes on
$\partial\rnm\setminus\{0\}$. We assume that there exists $C>0$ such that
$|h(x)|\leq C|x|^{1-n}$ for all $x\in\rnm\setminus\{0\}$. Then there
exists $\alpha\geq 0$ such that

$$h(x)=\alpha\frac{|x_1|}{|x|^n}$$
for all $x\in\rnm\setminus\{0\}$.
\end{prop}

\begin{proof} Up to rescaling, we assume that $h\geq 0$ in
$B_2(0)\setminus\{0\}$. We let

$$\alpha:=\max\left\{\lambda\geq 0/\, h(x)\geq
\lambda\frac{|x_1|}{|x|^n}\hbox{ for all }x\in\rnm\cap
\overline{B}_1(0)\right\}.$$
We let $\tilde{h}(x)=h(x)-\alpha\frac{|x_1|}{|x|^n}$ for all
$x\in\rnm$. The new function $\tilde{h}$ satisfies the hypothesis of
Proposition \ref{prop:rigidity}. In addition, it follows from the
definition of $\alpha$ and Hopf's maximum principle that

$$\liminf_{r\to 0}\inf_{|x|=r}\frac{\tilde{h}(x)|x|^n}{-x_1}=0.$$
Mimicking what was done in Steps \ref{sec:green}.5 and \ref{sec:green}.7,
we get that

$$\liminf_{r\to 0}\sup_{|x|=r,\,
x\in\rnm}\frac{\tilde{h}(x)|x|^n}{-x_1}=0.$$
We let
$$\hat{h}(x_1,\tilde{x}):=\left\{\begin{array}{ll}
\tilde{h}(x_1,\tilde{x}) & \hbox{ if }x_1\leq 0\hbox{ and
}(x_1,\tilde{x})\neq 0\\
-\tilde{h}(-x_1,\tilde{x}) & \hbox{ if }x_1>0.
\end{array}\right.$$
As easily checked, we have that $\hat{h}\in C^2(\rn\setminus\{0\})$ and
$\Delta\hat{h}=0$ in $\rn\setminus\{0\}$. With the definition of
$\hat{h}$, we immediately get that

$$\liminf_{r\to 0}\sup_{|x|=r}\frac{|\hat{h}(x)|\cdot |x|^n}{|x_1|}=0.$$
We let $(r_k)_{k>0}$ such that $\lim_{k\to +\infty}r_k=0$ and

$$\liminf_{k\to +\infty}\sup_{|x|=r_k}\frac{|\hat{h}(x)|\cdot
|x|^n}{|x_1|}=0.$$
We let $\tilde{\eta}\in C^\infty(\rn)$ such that $\tilde{\eta}\equiv 0$ in
$B_1(0)$ and $\tilde{\eta}\equiv 1$ in $\rn\setminus B_2(0)$. We let
$\eta_k(x):=\tilde{\eta}(x/r_k)$. Mimicking what was done in Step
\ref{sec:green}.6, we let $\varphi\in C^\infty_c(\rn)$ and get that

\beq
\int_{\rn}\hat{h}\Delta\varphi\, dx&=&
\int_{\rn}\eta_k\hat{h}\Delta\varphi\, dx+o(1)\\
&=& \int_{\rn}\Delta(\eta_k\hat{h})\cdot(\varphi-\varphi(0))\,
dx+\varphi(0) \int_{\rn}\Delta(\eta_k\hat{h})\, dx+o(1)\\
&=& o(1)+\varphi(0) \int_{\rn}\Delta(\eta_k\hat{h})\, dx
\eeq
We let $R>3$, and choose $k_0$ such that $0<r_k<1$ for $k>k_0$. We then
get that

\beq
\left|\int_{\rn}\Delta(\eta_k\hat{h})\,
dx\right|&=&\left|\int_{B_R(0)}\Delta(\eta_k\hat{h})\,
dx\right|=\left|\int_{\partial B_R(0)}\partial_\nu(\eta_k\hat{h})\,
d\sigma\right|\\
&=&\left|\int_{\partial B_R(0)}\partial_\nu\hat{h}\, d\sigma\right|\leq
C R^{n-1} R^{-n}\leq \frac{C}{R}
\eeq
Letting $R\to +\infty$, we get that $\int_{\rn}\Delta(\eta_k\hat{h})\,
dx=0$. We finally get that

$$\int_{\rn}\hat{h}\Delta\varphi\, dx=0$$
for all $\varphi\in C^\infty_c(\rn)$. As a consequence, $\Delta\hat{h}=0$
in ${\mathcal D}'(\rn)$, and $\hat{h}\in C^2(\rn)$. Since there exists
$C>0$ such that $|\hat{h}(x)|\leq C |x|^{1-n}$, we then get that
$\hat{h}$ is uniformly bounded on $\rn$. Since $\Delta\hat{h}=0$, we get
that $\hat{h}\equiv 0$. In particular,

$$h(x)=\alpha\frac{|x_1|}{|x|^n}$$
for all $x\in\rnm\setminus\{0\}$.
\end{proof}

\section{Appendix C: Symmetry of the positive solutions to the limit 
equation}\label{sec:sym}

This section is devoted to the proof of a symmetry property for the 
positive solutions to the limit equations involved in Proposition 
\ref{prop:exhaust}.

\begin{prop}\label{prop:sym}
Let $n\geq 3$ and $s\in (0,2)$. We let $u\in C^2(\rnm)\cap 
C^1(\overline{\rnm})$ such that
\bequa\label{sys:sym}\left\{\begin{array}{ll}
\Delta u=\frac{u^{\crit-1}}{|x|^s}& \hbox{ in }\rnm\\
u>0 & \hbox{ in }\rnm\\
u=0 & \hbox{ on }\partial\rnm,
\end{array}
\right.\eequa
where $\crit=\frac{2(n-s)}{n-2}$. We assume that there exists $C>0$ such 
that $u(x)\leq C(1+|x|)^{1-n}$ for all $x\in\rnm$. Then we have that 
$u\circ\sigma=u$ for all isometry of $\rn$ such that $\sigma(\rnm)=\rnm$. 
In particular, there exists $v\in C^2(\rr_-^\star\times \rr)\cap 
C^1(\rr_-\times \rr)$ such that for all $x_1<0$ and all $x'\in\rr^{n-1}$, 
we have that $u(x_1,x')=v(x_1,|x'|)$.
\end{prop}

We prove the Proposition in the sequel. We let $u\in C^2(\rnm)\cap 
C^1(\overline{\rnm})$ that verifies the system (\ref{sys:sym}) and such 
that there exists $C>0$ such that
\bequa\label{ineq:sym}
u(x)\leq \frac{C}{(1+|x|)^{n-1}}
\eequa
for all $x\in\rnm$. We $\vec{e}_1$ be the first vector of the canonical 
basis of $\rn$. We let the open ball
$$D:=B_{1/2}\left(-\frac{1}{2}\vec{e}_1\right).$$
We define
\bequa
v(x):=|x|^{2-n}u\left(\vec{e}_1+\frac{x}{|x|^2}\right)
\eequa
for all $x\in \overline{D}\setminus\{0\}$. We prolongate $v$ by $0$ at 
$0$. Clearly, this is well-defined.

\medskip\noindent{\bf Step \ref{sec:sym}.1:} We claim that
\bequa\label{ineq:hopf}
v\in C^2(D)\cap C^1(\overline{D})\hbox{ and }\frac{\partial v}{\partial 
\nu}<0\hbox{ on }\partial D
\eequa
where $\partial/\partial\nu$ denotes the outward normal derivative.
\begin{proof} It follows from the assumptions on $u$ that $v\in C^2(D)\cap 
C^1(\overline{D}\setminus\{0\})$. Moreover, $v(x)>0$ for all $x\in D$ and 
$v(x)=0$ for all $x\in\partial D\setminus\{0\}$. It follows from the 
estimate (\ref{ineq:sym}) that there exists $C>0$ such that
\bequa\label{ineq:v:1}
v(x)\leq C|x|
\eequa
for all $x\in \overline{D}\setminus \{0\}$. Since $v(0)=0$, we have that 
$v\in C^0(\overline{D})$. The function $v$ verifies the equation
\bequa\label{eq:v}
\Delta 
v=\frac{v^{\crit-1}}{|x+|x|^2\vec{e}_1|^s}=\frac{v^{\crit-1}}{|x|^s\left|x+\vec{e}_1\right|^s}
\eequa
in $D$. Since $-\vec{e}_1\in \partial D\setminus\{0\}$ and $v\in 
C^1(\overline{D}\setminus\{0\})\cap C^0(\overline{D})$, there exists $C>0$ 
such that
\bequa\label{ineq:v:2}
v(x)\leq C |x+\vec{e}_1|
\eequa
for all $x\in \overline{D}$. It then follows from (\ref{ineq:v:1}), 
(\ref{eq:v}), (\ref{ineq:v:2}) and standard elliptic theory that $v\in 
C^1(\overline{D})$. Since $v>0$ in $D$, it follows from Hopf's Lemma that 
$\frac{\partial v}{\partial \nu}<0$ on $\partial D$.\end{proof}

We prove the symmetry of $u$ by proving a symmetry property of $v$, which 
is defined on a ball. Our proof uses the moving plane method. We take 
largely inspiration in \cite{gnn} and \cite{cgs}. Classically, for any 
$\mu\geq 0$ and any $x=(x',x_n)\in \rn$ ($x'\in \rr^{n-1}$ and 
$x_n\in\rr$), we let
$$x_\mu=(x',2\mu-x_n)\hbox{ and }D_\mu=\{x\in D/\, x_\mu\in D\}.$$
It follows from Hopf's Lemma (see (\ref{ineq:hopf})) that there exists 
$\epsilon_0>0$ such that for any $\mu\in 
(\frac{1}{2}-\epsilon_0,\frac{1}{2})$, we have that $D_\mu\neq\emptyset$ 
and $v(x)\geq v(x_\mu)$ for all $x\in D_\mu$ such that $x_n\leq\mu$. We 
let $\mu\geq 0$. We say that $(P_\mu)$ holds if $D_{\mu}\neq \emptyset$ 
and
$$v(x)\geq v(x_{\mu})$$
for all $x\in D_{\mu}$ such that $x_n\leq\mu$. We let
\bequa\label{def:lambda}
\lambda:=\min\left\{\mu\geq 0/\, (P_{\nu})\hbox{ holds for all }\nu\in 
\left(\mu,\frac{1}{2}\right) \right\}.
\eequa

\medskip\noindent{\bf Step \ref{sec:sym}.2:} We claim that $\lambda=0$.
\begin{proof} We proceed by contradiction and assume that $\lambda>0$. We 
then get that $D_{\lambda}\neq\emptyset$ and that $(P_{\lambda})$ holds. 
We let
$$w(x):=v(x)-v(x_{\lambda})$$
for all $x\in D_{\lambda}\cap\{x_n<\lambda\}$. Since $(P_{\lambda})$ 
holds, we have that $w(x)\geq 0$ for all $x\in 
D_{\lambda}\cap\{x_n<\lambda\}$. With the equation (\ref{eq:v}) of $v$ and 
$(P_{\lambda})$, we get that
\beq
\Delta w&=& \frac{v(x)^{\crit-1}}{|x+|x|^2\vec{e}_1|^s}- 
\frac{v(x_{\lambda})^{\crit-1}}{|x_{\lambda}+|x_{\lambda}|^2\vec{e}_1|^s}\\
&\geq &v(x_{\lambda})^{\crit-1}\left(\frac{1}{|x+|x|^2\vec{e}_1|^s}- 
\frac{1}{|x_{\lambda}+|x_{\lambda}|^2\vec{e}_1|^s}\right)
\eeq
for all $x\in D_{\lambda}\cap\{x_n<\lambda\}$. With straightforward 
computations, we have that
\beq
&& |x_{\lambda}|^2-|x|^2=4\lambda(\lambda-x_n)\\
&& 
|x_{\lambda}+|x_{\lambda}|^2\vec{e}_1|^2-|x+|x|^2\vec{e}_1|^2=(|x_{\lambda}|^2-|x|^2)\left(1+|x_{\lambda}|^2+|x|^2+2x_1)\right)
\eeq
for all $x\in\rn$. It follows that $\Delta w(x)>0$ for all $x\in 
D_{\lambda}\cap\{x_n<\lambda\}$. Note that we have used that $\lambda>0$. 
It then follows from Hopf's Lemma and the strong comparison principle that
\bequa\label{ppty:w}
w>0\hbox{ in }D_\lambda\cap\{x_n<\lambda\}\hbox{ and }\frac{\partial 
w}{\partial \nu}<0\hbox{ on }D_\lambda\cap\{x_n=\lambda\}.
\eequa
By definition, there exists a sequence 
$(\lambda_i)_{i\in\mathbb{N}}\in\rr$ and a sequence 
$(x^i)_{i\in\mathbb{N}}\in D$ such that $\lambda_i<\lambda$, $x^i\in 
D_{\lambda_i}$, $(x^i)_n<\lambda_i$, $\lim_{i\to 
+\infty}\lambda_i=\lambda$ and
\bequa\label{ineq:sym:2}
v(x^i)<v((x^i)_{\lambda_i})
\eequa
for all $i\in\mathbb{N}$. Up to extraction a subsequence, we assume that 
there exists $x\in \overline{\left(D_{\lambda}\cap\{x_n< \lambda\}\right)}$ such 
that $\lim_{i\to +\infty}x^i=x$ with $x_n\leq \lambda$. Passing to the 
limit $i\to +\infty$ in (\ref{ineq:sym:2}), we get that $v(x)\leq 
v(x_{\lambda})$. It follows from this last inequality and (\ref{ppty:w}) 
that $v(x)-v(x_\lambda)=w(x)=0$, and then $x\in 
\partial(D_{\lambda}\cap\{x_n<\lambda\})$.

\smallskip\noindent{\it Case 1:} We assume that $x\in\partial D$. Then 
$v(x_{\lambda})=0$ and $x_{\lambda}\in\partial D$. Since $D$ is a ball and 
$\lambda>0$, we get that $x=x_{\lambda}\in\partial D$. Since $v$ is $C^1$, 
we get that there exists $\tau_i\in ((x^i)_n,2\lambda_i-(x^i)_n)$ such 
that
$$v(x^i)-v((x^i)_{\lambda_i})=\partial_n v((x')^i,\tau_i)\times 
2((x^i)_n-\lambda_i)$$
Letting $i\to +\infty$, using that $(x^i)_n<\lambda_i$ and 
(\ref{ineq:sym:2}), we get that $\partial_n v(x)\geq 0$. On the other 
hand, we have that
$$\partial_n v(x)=\frac{\partial v}{\partial\nu}(x)\cdot 
(\nu(x)|\vec{e}_n)=\frac{\lambda}{|x+\vec{e}_1/2|}\frac{\partial 
v}{\partial\nu}(x)<0.$$
A contradiction.

\smallskip\noindent{\it Case 2:} We assume that $x\in D$. Since 
$v(x_{\lambda})=v(x)$, we then get that $x_{\lambda}\in D$. Since $x\in 
\partial(D_{\lambda}\cap\{x_n<\lambda\})$, we then get that $x\in 
D\cap\{x_n=\lambda\}$. With the same argument as in the preceding step, we 
get that $\partial_n v(x)\geq 0$. On the other hand, since $x_n=\lambda$, we get with (\ref{ppty:w}) that $\partial_n v(x)=\frac{\partial_n w(x)}{2}<0$. A contradiction.

\smallskip\noindent In all the cases, we have obtained a contradiction. 
This proves that $\lambda=0$. \end{proof}

\medskip\noindent{\bf Step \ref{sec:sym}.3:} Here goes the final argument. 
Since $\lambda=0$, it follows from the definition (\ref{def:lambda}) of 
$\lambda$ that $v(x',x_n)\geq v(x',-x_n)$ for all $x\in D$ such that 
$x_n\leq 0$. With the same technique, we get the reverse inequality, and 
then, we get that
$$v(x',x_n)=v(x',-x_n)$$
for all $x=(x',x_n)\in D$. In other words, $v$ is symmetric with respect 
to the hyperplane $\{x_n=0\}$. The same analysis holds for any 
hyperplane containing $\vec{e_1}$. Coming back to the initial function 
$u$, this proves the Theorem.

\end{document}